\documentclass[reqno, 11pt]{amsart}
\usepackage{amsmath}
\usepackage{latexsym}
\usepackage{amsfonts}

\parskip=\smallskipamount
\newtheorem{Pa}{Paper}[section]
\newtheorem{Tm}[Pa]{{\bf Theorem}}
\newtheorem{La}[Pa]{{\bf Lemma}}
\newtheorem{Cy}[Pa]{{\bf Corollary}}
\newtheorem{Rk}[Pa]{{\bf Remark}}
\newtheorem{Pn}[Pa]{{\bf Proposition}}
\newtheorem{Pb}[Pa]{{\bf Problem}}
\newtheorem{Dn}[Pa]{{\bf Definition}}
\newtheorem{Ex}[Pa]{{\bf Example}}

\newcommand{\tH}{\widetilde{\mathcal {H}}}

\newcommand{\bs}{{\bf s}}
\newcommand{\br}{{\bf r}}
\newcommand{\Col}{\boldsymbol{\Sigma}}

\newcommand{\G}{\Gamma}

\newcommand{\C}{{\mathbb C}}
\newcommand{\B}{{\mathbb B}}

\newcommand{\D}{{\mathbb D}}
\newcommand{\cS}{\mathcal S}

\newcommand{\tU}{\widetilde{U}}
\newcommand{\td}{\widetilde{d}}

\newcommand{\T}{F}
\newcommand{\del}{T}
\newcommand{\hH}{\widehat{\mathcal H}}

\newcommand{\bU}{{\bf U}}

\newcommand{\cH}{{\mathcal{H}}}
\newcommand{\cG}{{\mathcal {G}}}
\newcommand{\cL}{{\mathcal {L}}}
\newcommand{\cK}{{\mathcal {K}}}
\newcommand{\cN}{{\mathcal {N}}}

\newcommand{\cF}{{\mathcal {F}}}

\newcommand{\cU}{{\mathcal U}}
\newcommand{\cY}{{\mathcal Y}}
\newcommand{\CE}{{\mathcal T}}
\newcommand{\sbm}[1]{\left[\begin{smallmatrix} #1
		\end{smallmatrix}\right]}

\newcommand{\K}{{\mathbb K}}
\begin{document}
\title[Interpolation in the noncommutative Schur-Agler class]
{Interpolation in the noncommutative Schur-Agler class}

\author[J. A. Ball]{Joseph A. Ball}
\address{Department of Mathematics,
Virginia Polytechnic Institute,
Blacksburg, VA 24061-0123, USA}
\email{ball@calvin.math.vt.edu}
\author[V. Bolotnikov]{Vladimir Bolotnikov}
\address{Department of Mathematics,
The College of William and Mary,
Williamsburg VA 23187-8795, USA}
\email{vladi@math.wm.edu}
\maketitle
\begin{quote}
The class of Schur-Agler functions over a domain ${\mathcal
D} \subset {\mathbb C}^{d}$ is defined as the class of holomorphic
operator-valued
functions on ${\mathcal D}$ for which a certain von Neumann inequality
is satisfied when a commuting tuple of operators satisfying a certain
polynomial norm inequality is plugged in for the variables.
Such functions are alternatively characterized as those
having a linear-fractional presentation which
identifies them as transfer functions of a certain type of
conservative structured multidimensional linear system.  There now has
been introduced a noncommutative version of the Schur-Agler class which
consists of formal power series in noncommuting indeterminants
satisfying a noncommutative version of the von Neumann inequality when a
tuple of operators (not necessarily commuting) coming from a
noncommutative operator ball are plugged in for the formal
indeterminants.  Formal power series in this noncommutative
Schur-Agler class in turn are characterized as those having a certain
linear-fractional presentation in noncommuting variables identifying
them as transfer functions of a recently introduced class of
conservative structure multidimensional linear systems having
evolution along a free semigroup rather than along an integer
lattice.  The purpose of this paper is to extend the previously
developed interpolation theory for the commutative Schur-Agler class
to this noncommutative setting.
   \end{quote}
\date{}

\section{Introduction} \label{S:intro}
\setcounter{equation}{0}

\textbf{The classical setting.}
By way of introduction we recall the
classical Schur class $\cS$ of analytic functions mapping the
unit disk $\D$ into the closed unit disk $\overline{\D}$. The
operator-valued {\em Schur class} $\cS(\cU, \cY)$ consists, by definition,
of analytic functions $F$ on $\D$ with values $F(z)$ equal to contraction
operators between two Hilbert spaces $\cU$ and $\cY$. In what follows, the
symbol $\cL(\cU,\cY)$ stands for the algebra of bounded linear
operators mapping $\cU$ into $\cY$, and we often abbreviate $\cL(\cU,
\cU)$ to $\cL(\cU)$. The class $\cS(\cU, \cY)$
admits several remarkable characterizations. In
particular any such function $F(z)$ can be realized in the form
\begin{equation}
F(z) = D + z C (I - zA)^{-1} B
\label{1.1}
\end{equation}
where the connecting operator (or {\em colligation})
$$
\bU = \begin{bmatrix} A & B \\ C & D \end{bmatrix} \colon
\begin{bmatrix} {\mathcal H} \\ \cU \end{bmatrix} \to
\begin{bmatrix} {\mathcal H} \\ \cY \end{bmatrix}
$$
is unitary, and where ${\mathcal H}$ is some auxiliary Hilbert space
(the {\em internal space} for the colligation). From the point of
view of system theory, the function \eqref{1.1} is the
{\em transfer function} of the linear system
$$
\boldsymbol{\Sigma} =\boldsymbol{\Sigma}(\bU) \colon \left\{
\begin{array}{rcl}
                x(n+1) & = & A x(n) + B u(n) \\
                y(n) & = & C x(n) + D u(n)
\end{array}.  \right.
$$
It is also well known that the Schur-class functions satisfy a
von Neumann inequality: {\em if $F \in \cS(\cU,\cY)$
and $T \in {\mathcal L}({\mathcal K})$ satisfies $\|T\|<1$, then
$F(T)$ is a contraction operator ($\|F(T)\| \le 1$), where $F(T)$ is defined by
$$
F(T) = \sum_{n=0}^{\infty} F_{n} \otimes T^n \in {\mathcal
L}(\cU \otimes {\mathcal K}, \cY \otimes{\mathcal K})
\quad\text{if}\quad F(z) = \sum_{n=0}^{\infty} F_{n}z^{n}.$$}
There is also a well-developed interpolation theory for the classical
Schur class.  One convenient formalism which encodes classical
Nevanlinna-Pick and Carath\'eodory-Fej\'er interpolation (see e.g.
\cite{BGR, FFGK}) proceeds as follows.
Making use of power series expansions one can introduce the
left and the right evaluation maps
\begin{equation}
F^{\wedge L}(T_{L})= \sum_{n=0}^\infty T_{L}^n
            F_n\quad\mbox{and}\quad F^{\wedge
            R}(T_{R})=\sum_{n=0}^\infty F_nT_{R}^n,
\label{1.2}
\end{equation}
which make sense for $F\in\cS(\cU,\cY)$ and for every choice of strictly
contractive operators $T_{L}\in {\mathcal L}(\cY)$ and $T_{R}\in
{\mathcal L}(\cU)$. One can then
formulate an interpolation problem with the data sets consisting of two
Hilbert spaces $\cK_{L}$ and $\cK_{R}$ and operators
\begin{align*}
& T_{L}\in\cL(\cK_{L}),\quad T_{R}\in\cL(\cK_{R}), \quad
X_L\in\cL(\cY,\cK_{L}), \\
&  Y_L\in\cL(\cU,\cK_{L}),\quad
X_R\in\cL(\cK_{R},\cY),\quad Y_R\in\cL(\cK_{R},\cU)
\end{align*}
as follows:
\begin{Pb}
Given the data as above, find necessary and sufficient conditions for
existence of a function $S\in\cS(\cU, \, \cY)$ such that
\begin{equation}
\left(X_L S\right)^{\wedge L}(T_{L})=Y_L \quad\text{and}\quad
\left(S Y_R\right)^{\wedge R}(T_{R})=X_R.
\label{1.3}
\end{equation}
\label{pb0}
\end{Pb}
The answer is well known: Problem \ref{pb0} has a solution if and only if
there exists a positive semidefinite operator
$P\in\cL(\cK_L\oplus\cK_R)$ subject to the Stein identity
$$
M^*PM-N^*PN=X^*X-Y^*Y
$$
where
\begin{equation}
M=\begin{bmatrix} I_{\cK_L} & 0 \\ 0 & T_R\end{bmatrix},\quad
N=\begin{bmatrix} T_L & 0 \\ 0 & I_{\cK_R} \end{bmatrix},\quad
X=\begin{bmatrix} X_L^* & X_R\end{bmatrix},\quad
Y=\begin{bmatrix} Y_L^* & Y_R\end{bmatrix}.
\label{1.4}
\end{equation}
\medskip
\textbf{Multivariable extensions.}
Multivariable generalizations of these and many other related results
have been obtained recently; one very general formulation introduced
(see \cite{at, AE, bb4}) proceeds as
follows.  Let
${\bf Q}$ be a $m\times k$ matrix-valued polynomial
\begin{equation}
{\bf Q}(z)=\begin{bmatrix}{\bf q}_{11}(z) & \ldots & {\bf q}_{1k}(z)\\
\vdots &&\vdots \\ {\bf q}_{m1}(z) & \ldots & {\bf
q}_{mk}(z)\end{bmatrix}: \; \C^d\to\C^{m\times k}
\label{1.5}
\end{equation}
and let ${\mathcal D}_{\bf Q}\in\C^d$ be the domain defined by
\begin{equation}
{\mathcal D}_{\bf Q}=\left\{z\in\C^d: \; \|{\bf Q}(z)\|<1\right\}.
\label{1.5a}
\end{equation}
For $\cU$ and $\cY$ two separable Hilbert spaces, in analogy with the
classical case it is natural to define the Schur class for the
domain ${\mathcal D}_{{\mathbf Q}}$ as the class ${\mathcal
S}_{{\mathbf Q}}(\cU, \cY)$  of holomorphic $\cL(\cU, \cY)$-valued
functions $S$ on ${\mathcal D}_{{\mathbf Q}}$ such that $\|S(z)\| \le
1$ for all $z \in {\mathcal D}_{{\mathbf Q}}$.  We say that an $S \in
{\mathcal S}_{{\mathbf Q}}(\cU, \cY)$ satisfies the
{\em ${\mathbf Q}$-von Neumann
inequality over ${\mathcal D}_{{\mathbf Q}}$} if $\|S(T_{1}, \dots,
T_{d})\| \le 1$ for all commuting tuples $(T_{1}, \dots, T_{d})$ of
operators on a Hilbert space $\cK$ with $\|{\mathbf Q}(T_{1}, \dots,
T_{d})\| < 1$.  (Here the fact that $\| {\mathbf Q}(T_{1}, \dots,
T_{d})\| < 1$ implies that the Taylor joint spectrum of $(T_{1},
\dots, T_{d})$ is contained in ${\mathcal D}_{{\mathbf Q}}$, so one
can use a tensored version of the Taylor functional calculus to
define $S(T_{1}, \dots, T_{d})$---see \cite{at}.) We define the
Schur-Agler class over ${\mathcal D}_{{\mathbf Q}}$, denoted by $\mathcal
{S A}_{{\mathbf Q}}(\cU, \cY)$, to consist of all $S \in {\mathcal
S}_{{\mathbf Q}}(\cU, \cY)$ which in addition satisfy the
${\mathbf Q}$-von Neumann inequality over ${\mathcal D}_{{\mathbf Q}}$.
As was first understood
for the tridisk case $({\mathbf Q}(z_{1}, z_{2}, z_{3}) = \sbm{ z_{1}&
0 & 0 \\ 0 & z_{2}& 0 \\ 0 & 0 & z_{3}}$), it can happen that the
containment $\mathcal {S A}_{{\mathbf Q}}(\cU, \cY) \subset {\mathcal
S}_{{\mathbf Q}}(\cU, \cY)$ is strict. It is this smaller class
$\mathcal{ SA}_{{\mathbf Q}}(\cU, \cY)$ which has a characterization
analogous to \eqref{1.1} and thereby can be interpreted as the set of
transfer functions of some type of conservative linear system, namely
(see \cite{bb4,AE}): {\em an $\cL(\cU, \cY)$-valued function  analytic on
${\mathcal D}_{\bf Q}$ belongs to the
class ${\mathcal{SA}}_{\bf Q}$ if and only if there exists
an auxiliary Hilbert space $\cH$  and a unitary operator
$$
U=\begin{bmatrix} A & B \\ C & D\end{bmatrix} \; : \;
\; \begin{bmatrix}\C^p\otimes{\cH} \\ \cU\end{bmatrix} \rightarrow
\begin{bmatrix}\C^q\otimes{\cH} \\ \cY\end{bmatrix}
$$
such that}
\begin{equation}
S(z)=D+C\left(I_{\C^p\otimes\cH}-({\bf Q}(z)\otimes\cH)A\right)^{-1}
({\bf Q}(z)\otimes\cH)B.
\label{1.6}
\end{equation}
Note that special choices of
\begin{equation}
{\bf Q}(z)={\rm diag} \, (z_1,\ldots,z_d)
\quad\mbox{and}\quad {\bf Q}(z)=\begin{bmatrix}z_1 &
z_2 & \ldots & z_d\end{bmatrix}
\label{1.7}
\end{equation}
lead to the unit polydisk ${\mathcal D}_{\bf Q}=\D^d$ and the unit ball
${\mathcal D}_{\bf Q}=\B^d$ of $\C^d$, respectively. The classes
${\mathcal {SA}}_{\bf Q}(\cU, \, \cY)$ for these two generic cases
have been known for a while. The polydisk setting was first
presented by J. Agler in \cite{agler-hellinger} and then extended to the
operator valued case in \cite{BLTT, bt}; see also \cite{aglmccar-poly,
bb3,BallMal}. The Schur-Agler functions on the unit ball appeared in
\cite{Drury} and later in \cite{agler2, quig, mccull2, aglmccar} in
connection with complete Nevanlinna-Pick kernels and in \cite{arv,
Popescu-Poisson} in connection with the study of commutative unitary
dilations of commutative row contractions; the Schur-Agler class for the
unit ball case has the extra structure that it can be
identified with the unit ball of the space of operator-valued
multipliers over the Arveson space (the reproducing kernel
Hilbert space over the unit ball ${\mathbb B}^{d} \subset {\mathbb C}^{d}$
with reproducing kernel $k_{d}(z,w) = \frac{1}{1 - \langle z, w
\rangle}$)---we refer to \cite{BTV} for a thorough review of the
operator-valued case.
The case when ${\mathcal D}_{\bf Q}$ is the Cartesian product of unit
balls (of arbitrary finite dimensions) was considered in \cite{tomerlin}.
Schur-Agler-class functions on $\D^d$ and $\B^d$ arise as the transfer
functions of {\em Givone-Roesser} (see \cite{Roesser, Kac}) and {\em
Fornasini-Marchesini} (see \cite{FM, Kac}) systems, respectively, which
satisfy an additional energy-balance relation (see \cite{BSV}).
In the general case, formula \eqref{1.6} can be interpreted as representing
$S$ as the
transfer function of a more general type of multidimensional conservative
linear system (see \cite[Section 4]{bb5} for more detail).

An interpolation problem similar to Problem \ref{pb0} has been studied
in \cite{bb5}. Interpolation conditions for this problem are the same as
in \eqref{1.3} but $T_{L}$ and $T_{R}$ are now commuting $d$-tuples
satisfying conditions
\begin{equation}
\|{\bf Q}(T_{L})\|<1\quad\mbox{and}\quad\|{\bf Q}(T_{R})\|<1
\label{1.8}
\end{equation}
and definitions of the left and the right evaluation maps are more
involved and rely on the Martinelli kernel (see \cite{vasilescu}) of the
Taylor functional calculus \cite{taylor1, taylor2}. Similarly to the one
variable case, the problem has a solution if and only if there is  a
positive semidefinite operator $P\in\cL((\cK_{L})^m\oplus(\cK_{R})^k)$
subject to the Stein identity
\begin{equation}
\sum_{j=1}^m M_j^*PM_j-\sum_{\ell=1}^k N_\ell^*PN_\ell=X^*X-Y^*Y
\label{1.9}
\end{equation}
where $X$ and $Y$ are the same as in \eqref{1.4} and $M_j$ and
$N_\ell$ are certain operators depending on $T_{L}$ and $T_{R}$
respectively (see \cite[Theorem 1.4]{bb5}).

\medskip
\textbf{The noncommutative setting.}
System theoretical aspects of the
above ideas has been extended recently \cite{Mal,BGM1,BGM2} to
noncommutative
multidimensional linear systems of a certain structure. These
systems, called {\em structured noncommutative multidimensional linear
systems} or SNMLSs in \cite{BGM1}) have
evolution along a free semigroup rather than along an integer lattice as
is usually taken in work in multidimensional linear system theory, and the
transfer function is a formal power series in noncommuting indeterminants
rather than an analytic function of several complex variables.
Furthermore, the transfer function of a {\em conservative} SNMLS
satisfies a certain von Neumann type inequality which leads to the
definition of a noncommutative Schur-Agler class associated with certain
noncommutative analogues of the domains ${\mathcal D}_{\bf Q}$ (but
where ${\bf Q}$ is restricted to be linear). The main
result \cite[Theorem 5.3]{BGM2} states that every noncommutative
Schur-Agler function admits a unitary realization similar to \eqref{1.6}.
The purpose of the present paper is to study related interpolation
problems of Nevanlinna-Pick type in the noncommutative Schur-Agler class.

The precise definitions and constructions involve a certain type
of graph (an ``admissible graph'' as defined below). Let $\G$ be a graph
consisting  of a  set of vertices $V = V(\G)$ and edges $E=E(\G)$.
An edge $e$ connects its {\em source vertex} $s$, denoted by $s = {\mathbf
s}(e) \in V$, to its
{\em range vertex} $r$, denoted by $r = {\mathbf r}(e) \in V$.
Following \cite{BGM1}, we say that $\G$ is {\em admissible} if it is a
finite ($V$ and $E$ are  finite sets)  bipartite graph such that each
connected component is a complete bipartite graph. The latter means that:
\begin{enumerate}
\item the set of vertices $V$ has a disjoint partitioning
$V = S \dot \cup R$ into the set of {\em source} vertices $S$ and
{\em range} vertices $R$,
\item $S$ and $R$ in turn have disjoint partitionings
$S = \dot \cup_{k=1}^{K} S_{k}$ and $R = \dot \cup_{k=1}^{K} R_{k}$ into
nonempty subsets $S_{1}, \dots, S_{K}$ and $R_{1}, \dots, R_{K}$ such
that, for each $s_{k} \in S_{k}$ and $r_{k} \in R_{k}$ (with the same
value of $k$) there is a unique edge $e = e_{s_{k},r_{k}}$
connecting $s_{k}$ to $r_{k}$ ($\bs(e) = s_{k}$, $\br(e) = r_{k}$), and
\item {\em every} edge of $\G$ is of this form.
\end{enumerate}
If $v$ is a vertex of $\G$ (so either $v \in S$ or $v \in R$) we denote
by $[v]$ the path-connected component $p$ (i.e., the complete bipartite
graph $p=\G_{k}$ with set of source vertices equal to $S_{k}$
and set of range vertices equal to $R_{k}$ for some $k = 1, \dots, K$)
containing $v$.  Thus, given two distinct vertices
$v_{1}, v_{2} \in S \cup R$, there is a
path of $\G$ connecting $v_{1}$ to $v_{2}$ if and only if $[v_{1}] =
[v_{2}]$ and this path has length 2 if both $v_{1}$ and $v_{2}$ are
either in $S$ or in $R$ and has length 1 otherwise. In case $s \in S$
and $r \in R$ are such that $[s] = [r]$, we shall use the notation
$e_{s,r}$ for the unique edge having $s$ as source vertex and $r$ as
range vertex:
$$
e_{s,r} \in E \text{ determined by }
\bs(e_{s,r}) = s, \, \br(e_{s,r}) = r.
$$
Note that $e_{s,r}$ is well defined only for  $s\in S$ and $ r \in R$
with  $[s] = [r]$.

\medskip

For an admissible graph $\G$, let $\cF_{E}$ be the free semigroup
generated by the edge set $E$ of $\G$.  An element of $\cF_{E}$ is then a
word $w$ of the form $w = e_{N} \cdots e_{1}$ where each $e_{k}$ is an
edge of $\G$ for $k = 1, \dots, N$.  We denote the empty word
(consisting of no letters) by $\emptyset$.  The semigroup operation is
concatenation: if $w = e_{N}\cdots e_{1}$ and $w' = e'_{N'} \cdots
e'_{1}$, then $ww'$ is defined to be
$$
ww' =e_{N} \cdots e_{1}e'_{N'} \cdots e'_{1}.
$$
Note that the empty word $\emptyset$ acts as the identity element for
this semigroup. On occasion we shall have use of the notation $w e^{-1}$
for a word $w \in \cF_{E}$ and an edge $e \in E$; by this notation we mean
\begin{equation}
\label{1.10}
     w e^{-1} = \begin{cases} w' &\text{if } w = w'e, \\
\text{undefined} & \text{otherwise.}\end{cases}
\end{equation}
with a similar convention for $e^{-1} w$. By $w^{\top}$ we mean
$e_{1}\cdots e_{N}$, the transpose of $w = e_{N}\cdots e_{1}$.

\medskip

For each $e \in E$, we define a matrix
$I_{\G, e} = [I_{\G, e; s, r}]_{s \in S, r \in R}$
(with rows indexed by $S$ and columns indexed by $R$) with matrix
entries given by
\begin{equation}
\label{1.11}
I_{\G, e; s, r} = \begin{cases}
1 & \text{ if } (s,r) = (\bs(e),\br(e)),\\ 0 & \text{ otherwise. }
\end{cases}
\end{equation}
We then define the \emph{structure matrix} $Z_{\G}(z)$
associated with each admissible graph $\G$ to be the linear form in the
noncommuting indeterminants $z = (z_{e}: e \in E)$ given by
\begin{equation}
Z_{\G}(z) = \sum_{e \in E} I_{\G, e} \, z_{e}.
\label{1.12}
\end{equation}
The latter function is the noncommutative analogue of ${\bf Q}(z)$ in
\eqref{1.5}. However, if we let the variables $(z_{e}: e \in E)$ in
\eqref{1.12} commute, we pick up only special examples of the
polynomial matrix functions ${\bf Q}(z)$, as will be clear from the
Examples below.

\begin{Ex}
\label{E:NCFMimp}
\textbf{Structure matrix for the noncommutative ball.}
{\rm
In this case, we take the admissible graph $\G^{FM}$ (where the label
``FM'' refers to {\em Fornasini-Marchesini} for system-theoretic
reasons explained in \cite{BGM1, BGM2}) to be a complete
bipartite graph having only one source vertex.  Thus we take
$S^{FM} = \{1\},$ and $R^{FM} = E^{FM} = \{1, \dots, d\}$ with
$\bs^{FM}(i)=1$, $\br^{FM}(i)=i$, i.e., $n = 1, m = d$.  Thus we have
$$
I_{\G^{FM}, i} = \begin{bmatrix} 0 & \cdots & 0 & 1 & 0 & \cdots &
0 \end{bmatrix},
$$
where 1 is located in the $i$-th slot. Thus, the structure
matrix for the noncommutative ball case is
given by
$$
Z_{\G^{FM}}(z) = \sum_{i=1}^{d} I_{\G^{FM}, i}z_{i} =
\begin{bmatrix} z_{1} & \cdots & z_{d} \end{bmatrix}.
$$
Note that when the variables $z_{1}, \dots, z_{d}$ commute, then the
associated domain $\{ z = (z_{1}, \dots, z_{d}) \colon \|
Z_{\G^{FM}}(z) \| < 1\}$ is the unit ball in ${\mathbb C}^{d}$.}
\end{Ex}
\begin{Ex}
\label{E:NCGRimp} \textbf{Structure matrix for the noncommutative polydisk.}
{\rm In this case, we take the
admissible graph $\G^{GR}$ (where the label ``GR'' refers to {\em
Givone-Roesser} for system-theoretic reasons explained in \cite{BGM1,
BGM2}) to have
$d$ path-connected components with each path-connected component
containing only one source and one range vertex.  Thus, we take
$S^{GR}= R^{GR} = E^{GR} = \{1,\dots,d\}$ with $\bs^{GR}(i)=i$,
$\br^{GR}(i)=i$ and thus $n = d = m$. Then $I_{\G^{GR}, i}$ is the
$d\times d$ matrix with $1$ located at the $(i, i)$-th entry
and with all other entries are zeros. Therefore, the
structure matrix for the noncommutative Givone-Roesser case has the
diagonal form
$$
Z_{\G^{GR}}(z) = \sum_{i=1}^{d} z_{i}I_{\G^{GR}, i} =
\begin{bmatrix} z_{1} & & \\ & \ddots & \\ & & z_{d} \end{bmatrix}.
$$
If the variables $z_{1}, \dots, z_{d}$ commute, then the associated
domain $\{ z = (z_{1}, \dots, z_{d}) \colon
\newline
\| Z_{\G^{GR}}(z) \| <
1\}$ is the unit polydisk in ${\mathbb C}^{d}$.}
\end{Ex}
\begin{Ex} \label{E:FullMatriximp} \textbf{Full matrix block
structure matrix.} {\rm In this case, we take $\G^{\text{full}}$ to be a
general finite, complete bipartite graph.  Thus we take $S = \{1,\dots,
n\}$, $R = \{1,\dots,m\}$, and $E = \{(i, j): i\in S, j\in R\}$ with
$\bs^{\text{full}}(i, j)=i$, $\br^{\text{full}}(i, j)=j$ where $d =
nm$. Then $I_{\G^{\text{full}}, (i, j)}$ is the $d\times d$ matrix with
$1$ located at the $(i, j)$-th entry and all other entries are zeros.
Thus the structure matrix for this case has the full-block structure
$$
Z_{\G^{\text{full}}}(z) =
\begin{bmatrix} z_{1,1} & \cdots & z_{1,m} \\ \vdots & & \vdots \\
z_{n,1} & \cdots & z_{n,m} \end{bmatrix}.
$$}
\end{Ex}
\begin{Ex} \label{E:GeneralSNMLSimp}
\textbf{The general structure matrix.}
{\rm Suppose that the admissible graph $\G$ has path connected components
$\G_{k}$
with  source vertices
$S_k = \{(k, 1),\dots , (k, n_k)\}$,  range vertices
$R_k = \{(k, 1),\dots , (k, m_k)\}$ and  edge sets
$E_k = \{(k, i, j): 1\leq i\leq n_k, 1\leq j\leq m_k\}$ for
$k = 1,\dots , K$. Define a graph $\G$ to have source vertex set
$$
S = \cup_{k=1}^K S_k =\{(k, i): 1\leq k\leq K, 1\leq i\leq n_k\},
$$
range vertex set
$$
R = \cup_{k=1}^K R_k =\{(k, j): 1\leq k\leq K, 1\leq j\leq m_k\}
$$
and edge set
$$
E = \cup_{k=1}^K E_k =\{(k, i, j): 1\leq k\leq K, 1\leq i\leq n_k,
1\leq j\leq m_k\}
$$
with $\bs(k, i, j)=(k, i)$, $\br(k, i, j)=(k, j)$ for $(k, i, j)\in
E$.
Then the associated structure matrix $Z_\G(z)$ is given by
$$
Z_{\G}(z) =\begin{bmatrix} Z_{\text{full}, 1}(z^{1}) & & \\ & \ddots & \\
& & Z_{\text{full}, K}(z^{K})
\end{bmatrix}
$$
where we let $z^{k}$ denote the $(n_{k} \cdot m_{k})$-tuple of
variables $z^{k} = (z_{k,i,j} \colon 1 \le i \le n_{k}; 1 \le j \le
m_{k})$ and where
$$
Z_{\text{full}, k}(z^{k}) =
\begin{bmatrix} z_{k, 1, 1} & \cdots & z_{k, 1, m_k} \\ \vdots & &
\vdots \\
z_{k, n_k, 1} & \cdots & z_{k, n_k , m_k} \end{bmatrix}
$$
is as in Example \ref{E:FullMatriximp} for $k = 1, \dots, K$.}
\end{Ex}

By the definition of an admissible graph as a graph with
path-connected components equal to complete bipartite graphs, we see that
Example \ref{E:GeneralSNMLSimp} amounts to the general case. Thus, the
case considered in the present framework corresponds
(in the commutative setting) not to
arbitrary polynomials \eqref{1.5}, but just to homogeneous linear
functions in which case, the corresponding domain ${\mathcal D}_{\bf Q}$
is the Cartesian product of finitely many Cartan domains of type $I$.
The proofs of realization and interpolation results in this particular
case are not substantially easier; however, most of needed constructions
can be expressed in terms of uniformly converging power series rather
than the Vasilescu's operator analogue of the Martinelli-Bochner kernel.
Thus, the transfer to the noncommutative setting via noncommutative
formal power series in this situation is much more clear.

\medskip

In what follows, $\cL(\cU,\cY)\langle\langle z\rangle\rangle$ will
stand for the space of formal power series $\T$ of the form
\begin{equation}
\T(z) = \sum_{v \in \cF_{E}} \T_{v} z^{v},\qquad \T_{v}\in\cL(\cU,\cY)
\label{1.13}
\end{equation}
in noncommutative variables $z=\{z_e: \; e\in E\}$ indexed by the edge
set $E$ of the admissible graph $\G$, with
coefficients $\T_{v}$ equal to bounded operators acting between Hilbert
spaces $\cU$ and $\cY$. Here $z^{\emptyset}=1$ and
$z^{w} = z_{e_{N}}z_{e_{N-1}}\dotsm z_{e_{1}}$ if
$w = {e_{N}}{e_{N-1}}\dotsm {e_{1}}$. Thus
$$
z^{w} \cdot z^{w'} = z^{ww'}, \qquad z^{w} \cdot z_{e} = z^{we} \;
\text{ for } \; w,w' \in \cF_{E} \; \text{ and } \; e \in E.
$$
On occasion we shall have need of multiplication on the right or
left by $z_{e}^{-1}$; we use the convention
\begin{equation}
z^{w}  z_{e}^{-1} = \begin{cases} z^{w e^{-1}} &\text{if } w
e^{-1} \in \cF_{E} \text{ is defined;} \\
0 &\text{if } w e^{-1} \text{ is undefined,}\end{cases}
\label{1.14}
\end{equation}
where we use the convention \eqref{1.10} for the meaning of $w e^{-1}$.
We use the obvious analogous convention to define $z_{e}^{-1} z^{w}$.

\medskip

Let $\del = (\del_{e}\colon e \in E)$ be a collection of bounded,
linear operators (not necessarily commuting) on some separable
infinite-dimensional Hilbert space $\cK$ (also indexed by the edge set $E$
of $\G$). We define an operator
$\T(\del) \colon \cU \otimes \cK \to  \cY \otimes \cK$ by
\begin{align}
& \T(\del) : =
\lim_{N \to \infty} \sum_{v \in \cF_{E} \colon |v| \le N}
\T_{v} \otimes \del^{v} \notag \\
& \qquad \text{ where } \del^{\emptyset} = I_{\cK} \text{ and }
\del^{v} = \del_{e_{N}} \cdots \del_{e_{1}} \text{ if }
v = e_{N} \cdots e_{1}
       \label{1.15}
      \end{align}
whenever the limit exists in the weak-operator topology. \footnote{In
\cite{BGM2} the limit is taken in the norm-operator topology; the
weak-operator topology is more convenient for our purposes here.}
In general there
is no reason for the limit in \eqref{1.15} to exist; on the other hand
if  $\T$ is a polynomial in $z$, its action on noncommutative tuples is
well defined. Alternatively, if $T$ is a nilpotent tuple (so that
$\T^{v} = 0$ once the length $|v|$ of $v$ is large enough),
then the expression \eqref{1.15} is well defined.
More generally, it is well defined if $\T(z)$ is a rational formal
power series and the tuple $\del$ is in the domain of $\T(z)$---see
\cite{AK1,AK2,HMV}.
Take the function $Z_{\G}$ as in \eqref{1.12}, define
(according to \eqref{1.15}) the operator
$$
     Z_{\G}(\del) := \sum_{e \in E} I_{\G; e} \otimes \del_{e} \in
\cL(\oplus_{r \in R} \cK, \oplus_{s \in S}\cK)
$$
and introduce the noncommutative structured ball
\begin{equation}
\label{1.16}
{\mathcal B}_{\G}\cL(\cK) = \{ \del =
(\del_{e})_{e \in E} \colon \del_{e} \in \cL(\cK) \text{ for
} e \in E \text{ and } \|Z_{\G}(\del)\| < 1 \}.
\end{equation}
Now we are in position to define the noncommutative Schur-Agler class.
\begin{Dn}
\label{D:SchurAgler}
Given an admissible graph $\G$, a formal power series \eqref{1.13} is
said to belong to the noncommutative Schur-Agler class
${\mathcal{SA}}_{\G}(\cU, \cY)$
if, for  each Hilbert space $\cK$ and each $\del= (\del_{e})_{e \in
E} \in{\mathcal B}_{\G}\cL(\cK)$, the limit
\begin{equation}  \label{funccal1}
\T(\del) = \lim_{N \to \infty} \sum_{v \in \cF_{E} \colon |v| \le N}
\T_{v} \otimes \del^{v}
\end{equation}
exists in the weak-operator topology and defines a contractive operator
$$
\T(\del)\colon\cU\otimes\cK\to\cY\otimes\cK,\quad\|\T(\del)\| \le 1.
$$
\end{Dn}
We remark that, for the noncommutative polydisk setting of Example
\ref{E:NCGRimp},
  Alpay and Kalyuzhny\u{\i}-Verbovetzki\u{\i} in \cite{AK2}
show that it suffices to check that the expression \eqref{funccal1} is
a contraction only for $\del \in {\mathcal B}_{\G}\cL(\cK)$ with $\cK$
a Hilbert space of arbitrarily large but finite dimension.
The noncommutative analogue of the unitary realization \eqref{1.6}
for the Schur-Agler class $\mathcal{ S A}_{\G}(\cU, \cY)$
was obtained in \cite{BGM2}. To formulate the result we
shall need some additional notation and terminology.

\medskip

First, given a collection $\cH=\left\{\cH_p\colon \; p\in P\right\}$ of
Hilbert spaces indexed by the  set $P$ of path-connected components of
$\G$, let
\begin{equation}
Z_{\G,\cH}(z) = \sum_{e \in E}I_{\G,\cH; e} z_{e}
\label{1.17}
\end{equation}
where $I_{\G, \cH; e} \colon \oplus_{r \in R} \cH_{[r]} \to
\oplus_{s \in S} \cH_{[s]}$ is given via matrix entries
$$
[I_{\G,\cH; e}]_{s,r} = \left\{\begin{array}{cc}
I_{\cH_{[\bs(e)]}} =I_{\cH_{[\br(e)]}}&\text{ if }s=\bs(e)\text{ and }
r = \br(e), \\ 0 & \text{otherwise.}\end{array}\right.
$$
Furthermore, let $z' = (z'_{e} \colon e \in E)$ be another
system of noncommuting indeterminants; while $z_{e} z_{e'} \ne z_{e'} z_{e}$
and $z'_{e} z'_{e'} \ne z'_{e'} z'_{e}$ unless $e = e'$,
we will use the convention that $z_{e}z'_{e'} = z'_{e'}z_{e}$ for all
$e,e' \in  E$.  We also shall need the convention
\eqref{1.14} to give meaning to expressions of the form
$$
z_{e}^{\prime -1}z^{v} z^{\prime v'} z_{e}^{-1} = (z^{v}z_{e}^{-1})
\cdot (z_{e}^{\prime -1}z^{\prime v'}) = z^{v e^{-1}} z^{\prime e^{-1} v'}.
$$
For $\T(z)$ of the form \eqref{1.13}, we will use the convention that
$$
\T(z)^{*} = \left( \sum_{v \in \cF_{E}} \T_{v} z^{v}\right)^{*} :=
\sum_{v \in \cF_{E}} \T_{v}^{*} z^{v^{\top}} = \sum_{v \in \cF_{E}}
\T_{v^{\top}}^{*} z^{v}.
$$
We also use the notation
$$
\operatorname{Row}_{x \in X} M_{x} = \begin{bmatrix} M_{x_{1}} &
\cdots & M_{x_{N}} \end{bmatrix}, \qquad
\operatorname{Col}_{x \in X} M_{x} = \begin{bmatrix} M_{x_{1}} \\
\vdots \\ M_{x_{N}} \end{bmatrix} \text{ if } X = \{x_{1}, \dots, x_{N} \}
$$
for block row and column matrices with rows or columns indexed by the
set $X$.

\begin{Tm}
Let $\T(z)$ be a formal power series in noncommuting indeterminants $z =
(z_{e} \colon e \in E)$ indexed by the $E$ of edges of the admissible
graph $\G$ with coefficients $\T_{v} \in \cL(\cU, \cY)$ for two Hilbert
spaces $\cU$ and $\cY$. The following are equivalent:
\begin{enumerate}
\item $\T$ belongs to the noncommutative Schur-Agler class
${\mathcal{SA}}_{\G}(\cU, \cY)$.
\item There exist a collection
$\cH=\left\{\cH_p\colon \;
p\in P\right\}$ of Hilbert spaces indexed by the
set $P$ of path-connected components  $\G$ and a unitary operator
$$
U = \begin{bmatrix} A & B \\ C & D \end{bmatrix} \colon
\begin{bmatrix}\oplus_{s \in S} \cH_{[s]} \\
\cU \end{bmatrix} \to  \begin{bmatrix} \oplus_{r \in R} \cH_{[r]}   \\
\cY\end{bmatrix}
$$
such that
\begin{equation}
\T(z) = D + C (I - Z_{\G,\cH}(z) A)^{-1} Z_{\G,\cH}(z) B
\label{1.18}
\end{equation}
where $Z_{\G,\cH}$ is defined in \eqref{1.17}.
\item There exist a collection of Hilbert spaces
$\cH=\{\cH_{p} \colon p \in P\}$ and a formal power series
\begin{equation}
H(z)=\operatorname{Row}_{s\in S} H_s(z)\in
\cL(\oplus_{s \in S} \cH_{[s]}, \cY)\langle\langle z \rangle\rangle
\label{1.19}
\end{equation}
     so that
\begin{equation}
\label{1.20}
I_\cY - \T(z) \T(z')^{*} = H(z) \left( I - Z_{\G,\cH}(z)
Z_{\G, \cH}(z')^{*} \right) H(z')^{*}.
\end{equation}
\item There exist a collection of Hilbert spaces
$\cH=\{\cH_{p} \colon p \in P\}$ and a formal power series
\begin{equation}
G(z) = \operatorname{Col}_{r\in R} G_r(z)\in
\cL(\cU,\oplus_{r \in R} \cH_{[r]})\langle\langle z \rangle\rangle
\label{1.21}
\end{equation}
     so that
\begin{equation}
\label{1.22}
I_\cU - \T(z)^* \T(z') =G(z)^* \left( I - Z_{\G,\cH}(z)^*
Z_{\G, \cH}(z') \right) G(z').
\end{equation}
\item There exist a collection of Hilbert spaces
$\cH=\{\cH_{p} \colon p \in P\}$ and formal power series
$H(z)$ and $G(z)$ as in \eqref{1.19}, \eqref{1.21} so that
relations \eqref{1.20}, \eqref{1.22} hold along with
\begin{equation}
\label{1.23}
\T(z)-\T(z')=H(z)\left( Z_{\G,\cH}(z)-Z_{\G,\cH}(z')\right)G(z').
\end{equation}
\end{enumerate}
\label{T:conservative-real}
\end{Tm}

A representation of the form \eqref{1.18} with $U = \sbm{ A & B \\
C & D}$ is called a {\em unitary realization} for $\T$, or, in more
detail in the terminology from \cite{BGM2},
a realization of $\T$ as the transfer function for the {\em
conservative Structured Noncommutative Multidimensional Linear System}
    $\boldsymbol{\Sigma} = \{ \G, \cH, \cU,
\cY, \bU\}$ (see  Section \ref{S:SNMLS} for further details).
Note that if $\T$ is of the form \eqref{1.18}, then representations
\eqref{1.20}, \eqref{1.22} and \eqref{1.23} are valid with
\begin{equation}
H(z)=C\left(I-Z_{\G,\cH}(z)A\right)^{-1}\quad\mbox{and}\quad
G(z)=\left(I- AZ_{\G,\cH}(z)\right)^{-1}B.
\label{1.24}
\end{equation}

Now we turn to the subject of the paper.
We shall consider bitangential interpolation problems
with the data set consisting of two Hilbert spaces $\cK_L$ and $\cK_R$,
two tuples $\del_L=\{\del_{L,e}\colon e\in E\}$ and
$\del_R=\{\del_{R,e}\colon e\in E\}$ of operators acting on $\cK_{L}$
and $\cK_{R}$ respectively,  and bounded operators
$$
X_L \colon \cY\to  \cK_L,\quad Y_L \colon \cU \to  \cK_L,\quad
X_R \colon \cK_R\to \cY,\quad Y_R \colon  \cK_R\to\cU.
$$
The pair $(T_{L},X_{L})$ will be said to be {\em left admissible}
(with respect to the Schur-Agler class $\mathcal{S A}_{\G}(\cU, \cY)$)
if the left-tangential evaluation map (with operator argument)
$H \mapsto (X_LH)^{\wedge L}(T_{L})$ given by
\begin{equation} \label{1.25'}
        (X_{L}H)^{\wedge L}(\del_{L})=\sum_{v\in\cF_E}\del_{L}^{v^\top}
        X_{L}H_v
\end{equation}
is well-defined (with convergence of the infinite series in the
weak-operator topology) whenever $H(z) = \sum_{v \in \cF_{E}}H_{v}
z^{v}$ is a formal power series of the form \eqref{1.19} appearing in
the representation \eqref{1.20} for a Schur-Agler class  formal power
series
$\T(z) = \sum_{v \in \cF_{E}} \T_{v} z^{v}
\in \mathcal{S A}_{\G}(\cU, \cY)$.  Whenever this is the case, from
the identity $\T(z) = D + H(z) Z_{\G, \cH}(z) B$ we read off that
then the left-tangential map  is also well-defined on the associated
Schur-Agler class formal power series $\T(z)$:
\begin{equation} \label{lefttan}
       (X_{L}\T)^{\wedge L}(\del_{L}) =
       \sum_{v \in \cF_{E}} \del_L^{v^{\top}} X_{L} \T_{v} =
       X_{L} D + \sum_{e \in E} \del_{L,e}
       [(X_{L} H_{\bs(e)})^{\wedge L}(\del_{L})]B_{\br(e)}.
\end{equation}
Similarly, we say that the
pair $(Y_{R},T_{R})$ is {\em right admissible} (with respect to the
Schur-Agler class $\mathcal{ S A}_{\G}(\cU, \cY)$) if the
right-tangential evaluation map (with operator argument)
$G\mapsto (G Y_{R})^{\wedge R}(\del_{R})$ given by
\begin{equation} \label{1.25''}
        (G Y_{R})^{\wedge R} (\del_{R})=\sum_{v\in\cF_E}
G_{v} Y_{R}\del_{R}^{v^\top}
\end{equation}
exists (with convergence of the infinite series in the weak-operator
topology) whenever $G(z) = \sum_{v \in \cF_{E}} G_{v} z^{v}$ is a
formal power series of the form \eqref{1.21} appearing in the
representation \eqref{1.22} for a Schur-Agler class formal power
series $\T(z) = \sum_{v \in \cF_{E}} \T_{v} z^{v} \in \mathcal{S
A}_{\cG}(\cU, \cY)$.  Using the identity $\T(z) = D + C Z_{\G, \cH}(z)
G(z)$ we then see that the right-tangential
evaluation map (with operator argument) is well-defined on the
associated Schur-Agler class formal power series $\T(z)$ as well:
\begin{equation}  \label{righttan}
       (\T Y_{R})^{\wedge R}(\del_{R}) =
       \sum_{v \in \cF_{E}} \T_{v} Y_{R} \del_{R}^{v^{\top}} =
       D Y_{R} + \sum_{e \in E} C_{\bs(e)}
       [(G_{\br(e)} Y_{R})^{\wedge R}(\del_{R}) ] T_{R,e}.
\end{equation}
The connections between left and right point evaluation with operator
argument given by \eqref{1.25'} and \eqref{1.25''} versus the
tensor-product functional calculus given by \eqref{funccal1} will be
discussed in Section \ref{S:admissible}.
We say that the data set
\begin{equation}
{\mathcal D}=\left\{\del_L, \; \del_R, \; X_L, \; Y_L, \; X_R, \; Y_R
\right\},
\label{1.26}
\end{equation}
is {\em admissible} (with respect to $\mathcal{ S A}_{\G}(\cU, \cY)$)
if $( \del_{L}, X_{L})$ is left admissible and $(Y_{R}, \del_{R})$ is
right admissible.
We shall give examples and further details on admissible interpolation
data sets in Section \ref{S:admissible} below.

Given an admissible interpolation data set \eqref{1.26},
the formal statement of the associated bitangential interpolation
problem is:
\begin{Pb}
\label{P:1.4}
Find necessary and sufficient conditions for existence of a power series
$\T\in\mathcal{ S A}_{\G}(\cU, \cY)$ such that
\begin{equation}
\left(X_L\T\right)^{\wedge L}(\del_L)=Y_L \quad\text{and}\quad
\left(\T Y_R\right)^{\wedge R}(\del_R)=X_R.
\label{1.27}
\end{equation}
\end{Pb}

To formulate the solution criterion we need some additional notation.
Let ${\boldsymbol \delta}_{s,s'}$ be the Kronecker delta function
$$
{\boldsymbol \delta}_{s,s'} = \begin{cases} 1 &\text{if } s=s', \\
0 & \text{otherwise}.\end{cases}
$$
For $s\in S$ and $r\in R$, define operators
\begin{eqnarray}
E_{L,s}&=&\operatorname{Col}_{s'\in S\colon
[s']=[s]}{\boldsymbol \delta}_{s,s'}I_{\cK_L}
\colon \; \cK_L\to \bigoplus_{s'\in S:[s']=[s]}\cK_L,\label{1.28}\\
E_{R,r}&=&\operatorname{Col}_{r'\in R\colon
[r']=[r]}{\boldsymbol \delta}_{r,r'}I_{\cK_R}
\colon \; \cK_R\to \bigoplus_{r'\in R:[r']=[r]}\cK_R,\label{1.28a}\\
\widetilde{N}_{r}(\del_L)&=&\operatorname{Col}_{s'\in S\colon [s']=[r]}
\del_{L,e_{s',r}}^*\colon \; \cK_L\to \bigoplus_{s'\in S:[s']=[r]}\cK_L,
\label{1.29}\\
\widetilde{M}_{s}(\del_R)&=& \operatorname{Col}_{r'\in R\colon [r']=[s]}
\del_{R,e_{s,r'}}\colon \; \cK_R\to \bigoplus_{r'\in
R:[r']=[s]}\cK_R.\label{1.29a}
\end{eqnarray}
Define also the operators
\begin{eqnarray}
M_s&=&M_s(\del_R)=\begin{bmatrix} E_{L,s} & 0  \\ 0 &
\widetilde{M}_{s}(\del_R) \end{bmatrix}\quad (s\in S),\label{1.30}\\
\nonumber\\
N_r& =&N_r(\del_L)=\begin{bmatrix}\widetilde{N}_{r}(\del_L) & 0 \\ 0 &
E_{R,r}\end{bmatrix}\quad (r\in R).\label{1.31}
\end{eqnarray}
\begin{Tm}
There is a power series $\T\in{\mathcal {SA}}_\G(\cU, \, \cY)$
satisfying interpolation conditions \eqref{1.27} if and only if there
exists a collection $\K=\{\K_p\colon \, p\in P\}$ of positive semidefinite
operators
$$
\K_p\in\cL((\oplus_{s\in S\colon [s]=p}\cK_L)\oplus(\oplus_{r\in R\colon
[r]=p}\cK_R))
$$
indexed by the  set $P$ of path-connected components of $\G$,
which satisfies the Stein identity
\begin{equation}
\sum_{s\in S}M_s^*\K_{[s]} M_s-\sum_{r\in R}N_r^*\K_{[r]}
N_r=X^*X-Y^*Y,
\label{1.32}
\end{equation}
where $M_s$ and $N_r$ are the operators defined via formulas
\eqref{1.30}, \eqref{1.31} and where
\begin{equation}
X=\begin{bmatrix} X_L^* & X_R\end{bmatrix} \quad\mbox{and}\quad
Y=\begin{bmatrix} Y_L^* & Y_R\end{bmatrix}.
\label{1.33}
\end{equation}
\label{T:1.9}
\end{Tm}
Let $\K=\{\K_p\colon \, p\in P\}$ be any collection of operators
satisfying the conditions in Theorem \ref{T:1.9}.
Let us represent these operators more explicitly as
\begin{equation}
\K_p=\begin{bmatrix}\K_{p,L} & \K_{p,LR} \\ \K_{p,LR}^* &
\K_{p,R}\end{bmatrix}
\label{1.34}
\end{equation}
where
\begin{equation}
\K_{p,L}=\left[\Psi_{s,s'}\right],\quad
\K_{p,R}=\left[\Phi_{r,r'}\right],\quad
\K_{p,LR}=\left[\Lambda_{s,r}\right]
\label{1.35}
\end{equation}
for $s,s'\in S$ and $r,r'\in R$ such that $[s]=[s']=[r]=[r']=p$ and with
\begin{equation}
\Psi_{s,s'}\in\cL(\cK_L), \quad \Phi_{r,r'}\in\cL(\cK_R),\quad
\Lambda_{s,r}\in\cL(\cK_R,\cK_L).
\label{1.36}
\end{equation}
It turns out that for every collection $\K=\{\K_p\colon \, p\in P\}$
of positive semidefinite operators satisfying
\eqref{1.32}, there is a solution $\T$ of the bitangential
interpolation Problem \ref{P:1.4} such that, for some choice of associated
functions $H(z)$ and $G(z)$
of the form \eqref{1.19} and \eqref{1.21} in representations \eqref{1.20},
\eqref{1.22}, \eqref{1.23}, it holds that
\begin{align}
(X_LH_s)^{\wedge L}(\del_L)\left[(X_LH_{s'})^{\wedge L}(\del_L)\right]^*
& = \Psi_{s,s'} \;  \text{ for } s, s'\in
S\colon \; [s]=[s'],\label{1.37}\\
       (X_LH_s)^{\wedge L}(\del_L)\left(G_r Y_R\right)^{\wedge R}(\del_R)
        & = \Lambda_{s,r} \; \text{ for } s\in S; \; r\in
R\colon \; [s]=[r],\label{1.38}\\
       \left[(G_rY_R)^{\wedge R}(\del_R)\right]^*(G_{r'} Y_R)^{\wedge
R}(\del_R) & = \Phi_{r,r'} \; \text{ for } r,r'\in R\colon \;
[r]=[r'].\label{1.39}
\end{align}
Furthermore, it turns out that conversely, for every solution $\T$
of Problem \ref{P:1.4} with representations \eqref{1.20},
\eqref{1.22}, \eqref{1.23} (existence of these representations is
guaranteed by Theorem \ref{T:conservative-real}), the operators $\K_p$
defined via  (\ref{1.34})--(\ref{1.36})
and (\ref{1.37})--(\ref{1.39}) satisfy conditions of Theorem
\ref{T:1.9}. These observations suggest the following modification of
Problem \ref{P:1.4} with the data set
\begin{equation}
{\mathcal D}=\left\{\del_L, \; \del_R, \; X_L, \; Y_L, \; X_R, \; Y_R, \;
\Psi_{s,s'}, \; \Phi_{r,r'}, \, \Lambda_{s,r} \right\}.
\label{1.40}
\end{equation}
\begin{Pb}
\label{P:1.10}
Given the data ${\mathcal D}$ as in \eqref{1.40}, find all power series
$\T \in {\mathcal {SA}}_\G(\cU, \, \cY)$ satisfying
interpolation conditions \eqref{1.27} and such that for  some choice of
associated functions $H_s$ and $G_r$ in the representations
\eqref{1.20}, \eqref{1.22}, \eqref{1.23}, the
equalities \eqref{1.37}--\eqref{1.39} hold.
\end{Pb}
In contrast to Problem \ref{P:1.4}, the solvability criterion for Problem
\ref{P:1.10} can be given explicitly in terms of the interpolation data.

\begin{Tm}
Problem \ref{P:1.10} has a solution if and only if the operators
$\K_p$ ($p\in P$) given  by \eqref{1.34}, \eqref{1.35} are positive
semidefinite  and satisfy the Stein identity \eqref{1.32}.

Moreover, there exist Hilbert spaces $\widetilde \Delta$ and
$\widetilde \Delta_{*}$, a collection of Hilbert spaces
$\hH=\{\hH_p\colon \; p\in P\}$ indexed by set of path-connected
components $P$ of $\G$, and a formal power series
$$
\Sigma(z) = \begin{bmatrix} \Sigma_{11}(z) &
\Sigma_{12}(z) \\ \Sigma_{21}(z) & \Sigma_{22}(z) \end{bmatrix} \colon
\begin{bmatrix} \cU \\ \widetilde \Delta_{*} \end{bmatrix} \to
          \begin{bmatrix} \cY \\ \widetilde \Delta \end{bmatrix}
$$
from the noncommutative Schur-Agler class
$\mathcal{ S A}_{{\G}}(\cU\oplus\widetilde{\Delta}_*,
\; \cY\oplus \widetilde{\Delta})$ of the form
\begin{equation}
\Sigma(z) = \begin{bmatrix} U_{22} & U_{23} \\ U_{32} & 0
\end{bmatrix} + \begin{bmatrix} U_{21} \\ U_{31} \end{bmatrix}
(I_{\widetilde \Delta_{*}} - Z_{\G,\hH}(z)
U_{11})^{-1} Z_{\G,\hH}(z) \begin{bmatrix}
U_{12} & U_{13} \end{bmatrix}
\label{1.41}
\end{equation}
with
$$
\bU_{0} = \begin{bmatrix} U_{11} & U_{12} & U_{13} \\
U_{21} & U_{22} & U_{23} \\ U_{31} & U_{32} & 0 \end{bmatrix} \colon
       \begin{bmatrix} \oplus_{s\in S}\hH_{[s]} \\ \cU \\ \widetilde \Delta_{*}
\end{bmatrix} \to \begin{bmatrix} \widehat \oplus_{r\in R}\hH_{[r]} \\
\cY \\ \widetilde \Delta \end{bmatrix}
$$
unitary and completely determined by the interpolation data set
${\mathcal D}$ so that $\T$ is a solution of Problem \ref{P:1.10} if and
only if $\T$ has the form
\begin{equation}
\T(z) = \Sigma_{11}(z) + \Sigma_{12}(z)
\left( I_{\widetilde \Delta_{*}} - {\mathcal T}(z)
\Sigma_{22}(z)\right)^{-1}{\mathcal T}(z) \Sigma_{21}(z)
\label{1.42}
\end{equation}
for a power series ${\mathcal T}(z) \in \mathcal{SA}_\G(\widetilde \Delta,
\widetilde \Delta_{*})$.
\label{T:1.11}
\end{Tm}

As a corollary we have the following less satisfactory parametrization
of the set of all solutions of Problem \ref{P:1.4}.

\begin{Cy}  \label{C:parametrization}
     Suppose that we are given a noncommutative interpolation data set
     ${\mathcal D}$ as in \eqref{1.26} and let ${\mathbf K}$
     be the set of all collections ${\mathbb K} = \{ {\mathbb K}_{p}
     \colon p \in P\}$ of positive semi-definite operators
     ${\mathbb K}_{p} \in \cL((\oplus_{s \in S \colon [s] = p} \cK_{L}) \oplus
     (\oplus_{r \in R \colon [r] =p}\cK_{R}) )$ which satisfy the Stein
     identity \eqref{1.32}.  For each ${\mathbb K}
     \in {\mathbf K}$, let
     $$
      \Sigma^{{\mathbb K}}(z) =
      \begin{bmatrix} \Sigma^{{\mathbb K}}_{11}(z) &
      \Sigma^{{\mathbb K}}_{12}(z) \\ \Sigma_{21}^{{\mathbb K}}(z) &
      \Sigma^{{\mathbb K}}_{22}(z) \end{bmatrix} \colon \begin{bmatrix} \cU \\
      \widetilde{\Delta}^{{\mathbb K}}_{*} \end{bmatrix} \to
      \begin{bmatrix} \cY \\ \widetilde{\Delta}^{{\mathbb K}} \end{bmatrix}
      $$
      be the characteristic function associated with ${\mathbb K}$ as in
      Theorem \ref{T:1.11}.  Then the formal power series $\T(z) =
      \sum_{v \in \cF_{E}} T_{v} z^{v}$ with coefficients $T_{v} \in
      \cL(\cU, \cY)$ is a solution of Problem \ref{P:1.4} if and only if
      there is a choice of ${\mathbb K} \in {\mathbf K}$ and a
      free-parameter formal power series ${\mathcal T}(z)$ in the
      Schur-Agler class $\mathcal{S A}_{\G}
      (\widetilde{\Delta}^{{\mathbb K}}, \widetilde{\Delta}^{{\mathbb
      K}}_{*})$ so that $\T(z)$ has the form
      $$ \T(z) = \Sigma^{{\mathbb K}}_{11}(z) + \Sigma^{{\mathbb
      K}}_{12}(z) (I_{\widetilde{\Delta}^{{\mathbb K}}_{*}} - {\mathcal
      T}(z) \Sigma^{{\mathbb K}}_{22}(z))^{-1} {\mathcal T}(z)
      \Sigma^{{\mathbb K}}_{21}(z).
      $$
\end{Cy}

There has been some work on noncommutative interpolation theory of the
sort discussed here, but to this point it is not nearly as well
developed as the commutative theory.  All the previous work of which
we are aware has been in the context of the noncommutative-ball case
(see Example \ref{E:NCFMimp} above).  In this case the Schur-Agler
class $\mathcal
{S A}_{\G^{FM}}(\cU, \cY)$ can be identified with the space of
contractive multipliers on a Fock space of formal power series in
noncommuting indeterminants with norm-square-summable vector
coefficients, a noncommutative analogue of the unit ball of analytic
Toeplitz operators acting on the classical Hardy space
(see \cite{BV} and the references there). In particular,
Constantinescu and Johnson \cite{CJ} formulated and obtained a
necessary and sufficient condition (in terms of positivity of an
associated Pick matrix) for the existence of solutions for
an interpolation problem of the form (when translated to our notation)
$F^{\wedge R}(Z_{i}) = W_{i}$ ($i = 1,
\dots, N$) for the class $\mathcal{S A}_{\G^{FM}}({\mathbb C},
{\mathbb C})$.  A number of authors (see \cite{ariaspopescu,davidsonpitts,
popescu1}) have analyzed noncommutative analogues of the Sarason
formulation of interpolation for the noncommutative-ball setting; one
approach for these problems is as an application of the Commutant
Lifting Theorem developed by Popescu for this setting (see
\cite{Popescu-CLT1, Popescu-multi}).  A direction for future work is to
understand the connections of our approach via evaluation with
operator argument with the Sarason formulation and commutant lifting theory.
We mention that a very general version of commutant lifting theory
(with applications
to new sorts of interpolation problems) has recently been worked out by
Muhly and Solel \cite{MS}.

The paper is organized as follows. After the present Introduction,
Section \ref{S:SNMLS}  derives some consequences of the energy
balance relations
encoded in the conservative SNMLSs beyond what was derived in
\cite{BGM2} which are needed in the sequel.  These
consequences are then used in Section \ref{S:admissible} to derive
some necessary conditions for a given pair of operators $(X_{L},
\del_{L})$ (or $(\del_{R},Y_{R})$) to induce a well-defined
left (or right) tangential point evaluation
with operator argument on a given noncommutative Schur-Agler class
$\mathcal {S A}_{\G}(\cU, \cY)$.  Section \ref{S:solvability} then
establishes the criterion for existence of solutions in Theorems
\ref{T:1.9} and \ref{T:1.11}.  Section \ref{S:solutions}
establishes a correspondence between
solutions of Problem \ref{P:1.10} and unitary extensions of a certain
partially-defined isometry constructed from the data of the problem,
while Section \ref{S:universal} then uses the idea of Arov-Grossman
(see \cite{arovgros1}) to obtain the linear-fractional
parametrization for the set of all solutions of Problem \ref{P:1.10}
as described in Theorem \ref{T:1.11}.  Sections \ref{S:solvability},
\ref{S:solutions} and \ref{S:universal} closely parallel the analysis
of \cite{bb5} worked out for the commutative case.  The final Section
\ref{S:examples} discusses various examples and special cases.

\section{Conservative structured noncommutative multidimensional
linear systems} \label{S:SNMLS}
\setcounter{equation}{0}

Following \cite{BGM1, BGM2}
we define a {\em structured noncommutative multidimensional linear
system} (SNMLS)
to be a collection
\begin{equation}
\boldsymbol{\Sigma}=\{\G, \, \cH, \; \cU, \; \cY, \; \bU\}
\label{3.1}
\end{equation}
where $\G$ is an admissible graph, $\cH = \{\cH_{p} \colon p \in P\}$ is a
collection of (separable) Hilbert spaces (called state spaces)
indexed by the path-connected components $p$ of the graph $\G$,
where $\cU$ and $\cY$ are additional (separable) Hilbert spaces (to be
interpreted as the {\em input space} and the {\em output space}
respectively) and where $\bU$ is a {\em connection matrix}
(sometimes also called {\em colligation}) of the form
\begin{equation}
\label{3.2}
\bU=\begin{bmatrix} A & B \\ C & D \end{bmatrix} =
\begin{bmatrix} [A_{r,s}] & [B_{r}] \\ [C_{s}] & D \end{bmatrix}
\colon \begin{bmatrix} \oplus_{s \in S} \cH_{[s]} \\ \cU
\end{bmatrix} \to
\begin{bmatrix} \oplus_{r \in R} \cH_{[r]} \\ \cY \end{bmatrix}
\end{equation}
In case the connection matrix $\bU$ is unitary, we shall say that
$\boldsymbol{\Sigma}$ is a {\em conservative} or {\em unitary}
SNMLS. Associated with any SNMLS $\Col$ as
in \eqref{3.1} is the collection of
system equations with evolution along the free semigroup $\cF_{E}$
\begin{align} \label{syseq}
\boldsymbol{\Sigma} \colon \left\{ \begin{array}{rcl} x_{\bs(e)}(ew) & = &
{\Sigma}_{s\in S} A_{\br(e), s} x_s(w) + B_{\br(e)} u(w) \\
x_{s^{\prime}}(ew) & = & 0 \; {\text{ if }} \; s^{\prime} \neq \bs(e) \\
y(w) & = & {\Sigma}_{s\in S} C_s x_s(w) + D u(w) \; \text{ for } \; w \in
\cF_{E}.
\end{array} \right.
\end{align}

\begin{Rk}  {\em Suppose that
\begin{equation}
\widetilde{\Col}=\{\G, \; \widetilde{\cH}, \; \cU, \; \cY, \;
\widetilde{\bU}\}
\label{3.4}
\end{equation}
is another SNMLS with the same
structure graph $\G$ and the
same input and output spaces as in \eqref{3.1} and with the connecting
matrix
\begin{equation}
\label{3.5}
\widetilde{\bU}=\begin{bmatrix} \widetilde{A} & \widetilde{B} \\
\widetilde{C} & \widetilde{D} \end{bmatrix} =
\begin{bmatrix} [\widetilde{A}_{r,s}] & [\widetilde{B}_{r}] \\
[\widetilde{C}_{s}] & \widetilde{D} \end{bmatrix}
\colon \begin{bmatrix}\oplus_{s \in S} \widetilde{\cH}_{[s]} \\ \cU
\end{bmatrix}\to\begin{bmatrix}\oplus_{r\in R} \widetilde{\cH}_{[r]} \\ \cY
\end{bmatrix}.
\end{equation}

The colligations $\Col$ and $\widetilde{\Col}$ are said to be
{\it unitarily equivalent} if there is a collection
$\Upsilon = \{ \Upsilon_{p} \colon p \in P \}$ of unitary operators
$\Upsilon_{p} \colon \cH_{p} \to  \widetilde{\cH}_{p}$ (for each path
connected component $p$ of $\G$) such that
\begin{equation}
\label{3.6}
\begin{bmatrix} \oplus_{r \in R} \Upsilon_{[r]} & 0 \\
0 & I_{\cY} \end{bmatrix}\begin{bmatrix} A & B \\ C & D \end{bmatrix}
=\begin{bmatrix} \widetilde{A} & \widetilde{B} \\ \widetilde{C} &
\widetilde{D}\end{bmatrix}\begin{bmatrix}\oplus_{s \in S} \Upsilon_{[s]}
& 0 \\ 0 & I_{\cU}
\end{bmatrix}.
\end{equation}
   It is an easy computation to see that unitarily equivalent colligations
have the same transfer functions. It is much less obvious that under certain
minimality conditions (structure controllability and observability), the
colligations having the same characteristic functions are unitarily equivalent
(see \cite[Theorem 7.2]{BGM1} for the proof).}
\label{R:3.1}
\end{Rk}

It will be convenient to have the notation $p \mapsto s_{p}$ for a
{\em source-vertex cross-section}, i.e., for each path-connected
component $p$ of $\G$, $s_{p}$ is the assignment of a one particular
source vertex in the path-connected component $p$.
   From the structure of the system equations \eqref{syseq} and under the
assumption that $\bU$ is unitary (or more generally, under the
assumption that $\bU$ is contractive), we read off the following
properties for system trajectories $w \mapsto (u(w), x(w), y(w))$
satisfying equations \eqref{syseq}:
\begin{align}
       & x_{s}(e_{s,r} w) \text{ is independent of } s \text{ for any given }
       r \in R \text{ and } w \in \cF_{E}, \label{readoff1} \\
      & \sum_{r \in R} \| x_{s_{[r]}}(e_{s_{[r]},r}w) \|^{2} - \|
      x(w)\|^{2} \le \| u(w)\|^{2} - \|y(w)\|^{2}, \label{readoff2} \\
      & x_{s'}(ew) = 0 \; \text{ if } \; s' \ne \bs(e). \label{readoff3}
\end{align}
We may then compute
\begin{align*}
       \sum_{e \in E} \| x(ew)\|^{2} & =
       \sum_{s \in S} \sum_{e \in E} \| x_{s}(e w)\|^{2} \\
       & = \sum_{e \in E} \| x_{\bs(e)}(e w) \|^{2} \text{ (by
       \eqref{readoff3})} \\
       & = \sum_{p \in P} \sum_{s \in S \text{ and } r \in R \colon
       [s] = [r]} \| x_{s}(e_{s,r}w)\|^{2} \\
       & = \sum_{ p \in P} \sum_{r \in R \colon [r] = p}
       n_{s_{p}} \| x_{s_{p}} (e_{s_{p},r}w)\|^{2} \text{ (by
       \eqref{readoff1})}
\end{align*}
where we have set $n_{s_{p}}$ equal to the number of source vertices
$s$ in the path-connected component $p$ of $\G$.  If we now set
$N_{S}$ equal to the maximum number of source vertices in any
path-connected component of $\G$
\begin{equation} \label{NS}
     N_{S} = \max\{n_{s_{p}} \colon p \in P \},
\end{equation}
then
\begin{align*}
\sum_{e \in E} \frac{1}{N_{S}} \| x(ew)\|^{2} & =
\sum_{p \in P} \sum_{r \colon [r] = p} \frac{n_{s_{p}}}{N_{S}} \|
x_{s_{p}}(e_{s_{p,r}}w)\|^{2} \\
& \le \sum_{p \in P} \sum_{r \colon [r] =p} \|
x_{s_{p}}(e_{s_{p},r}w)\|^{2} \\
& = \sum_{r \in R} \| x_{s_{p}}(e_{s_{p},r}w)\|^{2} \\
& \le \| x(w)\|^{2} + \| u(w)\|^{2} - \| y(w) \|^{2} \text{ (by
\eqref{readoff2}).}
\end{align*}
Summing over all words $w$ of a fixed length $n$ and then multiplying
by $N_{S}^{-n}$ then gives
\begin{align}
      &  \sum_{w \colon |w| = n+1} \frac{1}{N_{S}^{n+1}} \|
       x(w)\|^{2} - \sum_{w \colon |w| = n} \frac{1}{N_{S}^{n}}\|
       x(w)\|^{2} \notag \\
     & \qquad   \qquad  \le \sum_{w  \colon |w| = n}\frac{1}{N_{S}^{n}}
\| u(w) \|^{2}
     - \sum_{w \colon |w| = n}\frac{1}{N_{S}^{n}} \| y(w)\|^{2}.
       \label{EB}
\end{align}
If we now sum over $n=0,1, \dots, N$, the left-hand side of \eqref{EB}
telescopes and we arrive at
\begin{equation} \label{EB'}
       \sum_{w \colon |w| = N+1} \frac{1}{N_{S}^{|w|}} \| x(w) \|^{2}
       - \| x(\emptyset)\|^{2} \le \sum_{w \colon |w| \le N}
       \frac{1}{N_{S}^{|w|}} \| u(w)
       \|^{2} - \sum_{w \colon |w| \le N}  \frac{1}{N_{S}^{{|w|}}}\| 
y(w) \|^{2}.
\end{equation}
In particular, we get the estimate
$$ \sum_{w \colon |w| \le N} \frac{1}{N_{S}^{|w|}} \|y(w)\|^{2} \le \|
x(\emptyset)\|^{2} + \sum_{w \colon |w| \le N} \frac{1}{N_{S}^{|w|}} \|
u(w) \|^{2}.
$$
Letting $N \to \infty$ then gives
\begin{equation} \label{finalEB}
       \sum_{w \in \cF_{E}} \frac{1}{N_{S}^{|w|}} \| y(w) \|^{2} \le
       \| x(\emptyset)\|^{2} + \sum_{w \in \cF_{E}}
       \frac{1}{N_{S}^{|w|}} \| u(w) \|^{2}
\end{equation}
for all system trajectories $(u,x,y)$ of the SNMLS $\Col$ as long as
the connection matrix $\bU$ satisfies $\| \bU \| \le 1$.

If $\{u(w)\}_{w \in \cF_{E}}$ is a $\cU$-valued input
string and $x(\emptyset)$ the initial state fed into the system equations
to produce a $\cY$-valued output string
$\{y(w)\}_{w \in \cF_{E}}$ and if we introduce the formal $Z$-transform
of the $\{u(w)\}_{w \in \cF_{E}}$ and $\{y(w)\}_{w \in \cF_{E}}$ according to
$$ \widehat u(z) = \sum_{w \in \cF_{E}} u(w) z^{w}, \qquad
     \widehat y(z) = \sum_{w \in \cF_{E}} y(w) z^{w},
$$
then it follows that
\begin{equation} \label{freqdomIO}
    \widehat y(z) = C(I - Z_{\G,\cH}(z) A)^{-1} x(\emptyset) +
\T_{\Col}(z) \cdot \widehat u(z)
\end{equation}
where $\T_{\Col}(z)$ is the formal noncommutative power series with
coefficients in $\cL(\cU, \cY)$ given by
\begin{eqnarray} \label{transfunc}
\T_{\Col}(z)&=&D+C(I - Z_{\G,\cH}(z) A)^{-1} Z_{\G,\cH}(z) B \label{3.3}\\
&=& \T_{\emptyset}+\sum_{N=1}^{\infty} \sum_{e_1, \dots ,e_N \in E}
C_{\bs(e_{N})} A_{\br(e_{N}), \bs(e_{N-1})} \cdots A_{\br(e_{2}),
\bs(e_{1})} B_{\br(e_{1})}z_{e_{N}}z_{e_{N-1}} \cdots z_{e_{2}}z_{e_{1}}
\nonumber
\end{eqnarray}
with $Z_{\G,\cH}$ given by \eqref{1.17}. In particular, if we take the
initial state $x(\emptyset)$ equal to $0$, we obtain the
relation $\widehat y(z) = \T_{\Col}(z) \cdot \widehat u(z)$ between
the $Z$-transformed input signal $\widehat u(z)$ and
the $Z$-transformed output signal $\widehat y(z)$.
We shall call $\T_{\Col}(z)$
the {\em transfer function} of the SNMLS $\Col$
(see \cite{BGM1, BGM2}).
The assertion of Theorem
\ref{T:conservative-real} then is that a power series $\T$ belongs to the
noncommutative Schur-Agler class $\mathcal{ SA}_{\G}(\cU, \cY)$ if and
only if
it is the transfer function of a conservative SNMLS $\Col$
of the form \eqref{3.1}.

\begin{Rk}
{\rm For future reference, we note that
the action of $\T_{\Col}(z)$ on a vector $u\in\cU$, namely
$$
\T_{\Col}(z)=D+C\left(I-Z_{\G,\cH}(z)A\right)^{-1}Z_{\G,\cH}(z)B: \;
u\to y
$$
is the result of the feedback connection
$$
\begin{bmatrix}A & B \\ C& D\end{bmatrix}\begin{bmatrix}h\\ u\end{bmatrix}
=\begin{bmatrix}h^\prime\\ y\end{bmatrix},\qquad h=Z_{\G,\cH}(z)h^\prime
$$
where $h\in \oplus_{s \in S} \widetilde{\cH}_{[s]}$ and
$h'\in  \oplus_{r \in R} \widetilde{\cH}_{[r]}$}.
\label{R:3.2}
\end{Rk}

\begin{Rk}  \rm{For the special case where $n_{s_{p}} = 1$ for each
path-connected component and $\bU$ is isometric, it is easily verified
that one gets equality in \eqref{EB} and \eqref{EB'}.  Thus, in this
case $N_{S} = 1$ and \eqref{finalEB} holds with $N_{S} = 1$. All this
has already been noted in \cite{BGM2} (see Remark 5.14 there) where such
graphs $G$ are called {\em row-sum} graphs.  A particularly nice case
of a row sum graph is a Fornasini-Marchesini graph (a row-sum graph
with one path-connected component)---see Example
\ref{E:NCFMimp}.  Then the system and the associated noncommutative
function theory have a particularly nice structure---see
\cite{Popescu-multi, BV}.}
\end{Rk}

\section{Admissible interpolation data sets} \label{S:admissible}
\setcounter{equation}{0}

With these preliminaries out of the way, we now turn to the issue of
identifying large classes of examples of left-admissible
and right-admissible pairs $(\del_{L}, X_{L})$ and $(Y_{R}, \del_{R})$
for a general admissible graph $\G$.  In particular, we shall see that
the class of interpolation problems covered in Problem \ref{P:1.4}
and \ref{P:1.10} is nonempty.

We first note the following relations (stated here without proof)
between left and right
       evaluation with operator argument \eqref{1.25'} and \eqref{1.25''}
       and tensor-product functional calculus \eqref{funccal1}.

\begin{Pn}
\label{P:funccal}
Assuming that all the functional
evaluations below exist, we have the following relations among
tensor-product evaluation  with operator argument \eqref{funccal1},
left evaluation with operator argument \eqref{1.25'} and right
evaluation with operator argument \eqref{1.25''}.
        \begin{enumerate}
	\item
        Let $\T(z) = \sum_{v \in \cF_{E}} \T_{v} z^{v}$ be a formal power
        series with coefficients in  $\cL(\cU, \cY)$ and let $\del =
        (\del_{e} \colon e \in E)$ be a tuple of operators on the space
        $\cK$.  Define a new power series $\T^{\sim}(z)$ by
        $$
        \T^{\sim}(z) = \sum_{v \in \cF_{E}} \T_{v^{\top}} z^{v}\quad
        \text{if}\quad
          \T(z) = \sum_{v \in \cF_{E}} \T_{v} z^{v}.
        $$
        Denote by $(\T^{\sim} \otimes I_{\cK})(z)$ the power series
        $$
          (\T^{\sim} \otimes I_{\cK}) (z) = \sum_{ v \in \cF_{E}}
          (\T_{v^{\top}} \otimes I_{\cK}) z^{v}.
        $$
       Then
       \begin{equation}  \label{relation1}
      (\T^{\sim} \otimes I_{\cK})^{\wedge L}(I_{\cY} \otimes \del) =
       \T(\del) = (\T^{\sim}\otimes I_{\cK})^{\wedge R}(I_{\cU}\otimes
       \del).
       \end{equation}

       \item
       If $f(z) = \sum_{v \in \cF_{E}} f_{v} z^{v}$ is a formal power
       series with scalar coefficients (so $f_{v} \in {\mathbb C}$ for
       all $v \in \cF_{E}$), then
       \begin{equation}  \label{relation2}
         (f^{\sim}\otimes I_{\cK})^{\wedge L}(\del) = f(\del)  =
         (f^{\sim}\otimes I_{\cK})^{\wedge R}(\del).
       \end{equation}

       \item
       If $f(z) = \sum_{v \in \cF_{E}} f_{v} z^{v}$ is a formal power
       series with scalar coefficients as in \#2 above and if $x$ is a
       vector in $\cK$, then
       \begin{equation}   \label{relation3}
       f(T) x = (x \cdot f^{\sim})^{\wedge L}(\del).
       \end{equation}

       \item If $\T(z) = \sum_{v \in \cF_{E}} \T_{v} z^{v}$ is a formal
       power series with coefficients $\T_{v} \in \cL(\cU, \cY)$,
       $\lambda = (\lambda_{e})_{e} $ is a tuple of complex numbers
       considered as operators on ${\mathbb C}$, and $X_{L} \in \cL(\cY,
       \cK_{L})$ and $Y_{R} \in \cL(\cK_{R}, \cU)$, then
       \begin{align}
           (X_{L}\T)^{\wedge L}(\lambda \cdot I_{\cK_{L}}) & = X_{L}
           \T(\lambda), \label{relation4} \\
           (\T Y_{R})^{\wedge R}(\lambda \cdot I_{\cK_{R}}) & =
           \T(\lambda) Y_{R}.  \label{relation5}
       \end{align}

       \end{enumerate}
        \end{Pn}

        \begin{Rk}  {\rm The left-side of \eqref{relation3} is the type of
        point evaluation used by Rosenblum-Rovnyak to formulate the
        so-called Nudelman interpolation problem in \cite{RR}.  Relation
        \eqref{relation3} shows how this type of interpolation condition
        can be converted to the version of Nudelman interpolation for the
        classical case used in \cite{BGR}.  An alternative extension of
        the Rosenblum-Rovnyak Nudelman problem to the formal power series
        setting is given in \cite{popescu1}.  In the sequel we shall
        have use of only part (4) of Proposition \ref{P:funccal}. }
        \end{Rk}

By definition, a formal power series $\T(z) = \sum_{v \in \cF_{E}} \T_{v}
z^{v}$ is in the Schur-Agler class if and only if $\T(\del) $ (defined
via \eqref{funccal1}) is a contraction for all $\del \in {\mathcal
B}_{\G}(\cK)$.  Given operators $X_{L} \in \cL(\cY, \cK_{L})$ and
$Y_{R}\in \cL(\cK_{R}, \cU)$ and operator tuples $\del_{L} \in
\cL(\cK_{L})^{n_{E}}$ and $\del_{R} \in \cL(\cK_{R})^{n_{E}}$ (here
we use $n_{E}$ to denote the number of edges  $e \in E$ for the
admissible graph $\G$),  the hope would
be that $(\del_{L}, X_{L})$ would be left admissible as soon as
$\del_{L} \in {\mathcal B}_{\G}\cL(\cK_{L})$ and that $(Y_{R},
\del_{R})$ would be right admissible (with respect to
$\mathcal{ S A}_{\G}(\cU, \cY)$) as soon as  $\del_{R}$
is in ${\mathcal B}_{\G}\cL(\cK_{R})$.  As we shall see below,
this is indeed
correct in some special cases while we obtain only partial results in this
direction for the case of a general admissible graph $\G$.  We begin
with the situation of part (4) in Proposition \ref{P:funccal}.

\begin{Pn}
\label{P:scalareval}
        Suppose that $\T(z) = \sum_{v \in \cF_{E}} \T_{v} z^{v}$ is a
        formal power series in the class $\mathcal{S A}_{\Gamma}(\cU, \cY)$
        and suppose that $\lambda = (\lambda_{e})_{e \in E}$ is a tuple of
        complex numbers.
        Then:
     \begin{enumerate}
         \item
         Suppose that $X_{L} \in \cL(\cY, \cK_{L})$ and that we let
         $\del_{L}$ be the tuple of scalar operators $\del_{L}
         = \lambda \cdot I_{\cK_{L}} = (\lambda_{e} \cdot
         I_{\cK_{L}})_{e \in E}$. Then $(\lambda_{e} \cdot I_{\cK_{L}},
         X_{L})$ is left admissible whenever $\| Z_{\G}(\lambda) \| < 1$.

         \item Suppose that $Y_{R} \in \cL(\cK_{R}, \cU)$ and that we let
         $\del_{R}$ be the tuple of scalar operators $\del_{R}
         = \lambda \cdot I_{\cK_{L}} = (\lambda_{e} \cdot
	  I_{\cK_{L}})_{e \in E}$.  Then $(Y_{R}, \lambda \cdot I_{\cK_{R}})$
	  is right admissible whenever $\| Z_{\G}(\lambda)\| < 1$.

     \item Suppose that the hypotheses of parts (1) and (2) hold with
     $\cK_{L} = \cY$ and $X_{L} = I_{\cY}$ and with $\cK_{R} = \cU$ and
     $Y_{R} = I_{\cU}$.  Then
     \begin{equation} \label{scalarvariable}
       F^{\wedge L}(\lambda \cdot I_{\cY}) =
       F^{\wedge R}(\lambda \cdot I_{\cU})
     \end{equation}
      \end{enumerate}
        \end{Pn}

     \begin{proof}  This is an immediate consequence of relations
     \eqref{relation4} and \eqref{relation5} in Proposition \ref{P:funccal}
     and the definitions.
     \end{proof}

We next explore the function of the scalar-tuple variable
$\lambda = (\lambda_{1}, \dots, \lambda_{n_{E}})$ a little further.
To simplify notation, in the statement of the next result we label the
edges of the graph $G$ by the letters $1,2, \dots, d$ where $d =
n_{E}$ is the number of edges of $G$.  Then words in $\cF_{E}$ have
the form $w = i_{N} i_{N-1} \cdots i_{1}$ where each $i_{\ell} \in
\{1, \dots, d \}$.  If $F(z) = \sum_{v \in \cF_{E}} F_{v} z^{v}$ is a
formal power series with coefficients in $\cL(\cU, \cY)$, the
function $F^{{\mathbf a}}(\lambda)$ of the scalar $d$-tuple
$(\lambda_{1}, \dots, \lambda_{d})$ given by either the left-hand
side or the right-hand side of \eqref{scalarvariable} (under the
assumption that the series converges)
can be expressed as
\begin{align*}
    F^{{\mathbf a}}(\lambda) & = \sum_{v \in \cF_{E}} F_{v} (\lambda 
I_{\cU})^{v}
     = \sum_{{\mathbf n} \in {\mathbb Z}^{d}_{+}} \left[ \sum_{v \colon v \in
     {\mathbf a}^{-1}({\mathbf n})} F_{v}\right] \lambda^{{\mathbf n}}
     = : \sum_{{\mathbf n} \in {\mathbb Z}^{d}_{+}} F^{{\mathbf
     a}}_{{\mathbf n}} \lambda^{{\mathbf n}}
\end{align*}
where we have introduced the {\em abelianization map} ${\mathbf a}
\colon {\mathcal F}_{d} \to {\mathbb Z}^{d}_{+}$ given by
$$
{\mathbf a}(i_{N} \cdots i_{1}) = (n_{1}, \dots, n_{d}) \; \text{ if } \;
n_{j} = \# \{ \ell \colon i_{\ell} = j \} \;
\text{ for } \; j = 1, \dots, d,
$$
where $\lambda^{v} = \lambda_{i_{N}} \cdots \lambda_{i_{1}}$
if  $v = i_{N} \cdots i_{1}$  and where
$\lambda^{{\mathbf n}} = \lambda_{1}^{n_{1}} \cdots
\lambda_{d}^{n_{d}}$  if  ${\mathbf n} = (n_{1}, \dots, n_{d})$, and
where we have set
$$
   F^{{\mathbf a}}_{{\mathbf n}} = \sum_{v \colon v \in {\mathbf
   a}^{-1}({\mathbf n})} F_{v}.
$$
If $F \in \mathcal{S A}_{\G}(\cU, \cY)$ then necessarily $F^{{\mathbf
a}}$ is analytic on ${\mathcal D}_{Z_{\G}^{{\mathbf a}}}$ where
$Z_{\G}^{{\mathbf a}}(\lambda)$ is just the abelianization of the
structure matrix $Z_{\G}(z)$ for $\G$.
For a general  matrix-valued polynomial ${\mathbf Q}(\lambda)$ in the
commuting variables $\lambda = (\lambda_{1}, \dots, \lambda_{d})$,
the associated {\em commutative} Schur-Agler class $\mathcal{S
A}_{{\mathbf Q}}(\cU, \cY)$ was defined in \cite{bb5} to consist of
holomorphic functions $\lambda \mapsto F(\lambda)$ defined on the domain
${\mathcal D}_{{\mathbf Q}}: = \{ \lambda \in {\mathbb C}^{d} \colon
\| {\mathbf Q}(\lambda) \| < 1 \}$ such that $\|F(\del)\| \le 1$ for any
commuting $d$-tuple of operators $(\del_{1}, \dots, \del_{d})$ on
${\mathcal K}$ such that ${\mathbf Q}(\del_{1} \dots, \del_{d}) \| <
1$.  For the special case where ${\mathbf Q}$ is taken to be the
abelianized structure matrix ${\mathbf Q}(\lambda) = Z^{{\mathbf
a}}_{\G}(\lambda)$, then we see that the the set of commuting
$d$-tuples $\del$ with
$\| Z_{\G}^{{\mathbf a}}(\del) \| < 1$ is just the intersection of
${\mathcal B}_{\G}\cL(\cK)$ with commutative operator tuples. A consequence of
Lemma 1 from \cite{at} is that a commuting $d$-tuple $\del  =
(\del_{1}, \dots, \del_{d})$ has its Taylor spectrum in the domain
${\mathcal D}_{Z_{\G}^{{\mathbf a}}}$ whenever
$\|Z_{\G}^{{\mathbf a}}(\del) \| < 1$.
Moreover, as $Z_{\G}^{{\mathbf a}}$ is a {\em linear} polynomial, the
associated domain ${\mathcal
D}_{Z_{\G}^{{\mathbf a}}}$ is a logarithmically convex Rinehardt domain,
and the functional calculus with operator argument defined via the Taylor
functional calculus can equivalently be carried out by using
power series centered at the origin (see \cite[Remark 2.2]{bb5}).
Hence, if $\del_{L} = (\del_{L,j})_{j=1,\dots,d}$ is a {\em commuting}
$d$-tuple of operators on $\cK$ and $X_{L} \in \cL(\cY, \cK_{L})$, then
\begin{equation} \label{abelianization}
      (X_{L}F)^{\wedge L}(\del_{L})  =
      \sum_{v \in \cF_{E}} \del_{L}^{v} X_{L} F_{v}
      = \sum_{{\mathbf n} \in {\mathbb Z}^{d}_{+}} \del_{L}^{{\mathbf n}}
      X_{L} F^{{\mathbf a}}_{{\mathbf n}}
      = (X_{L} F^{{\mathbf a}})^{\wedge L}(\del_{L})
\end{equation}
where $(X_{L} F^{{\mathbf a}})^{\wedge L}(\del_{L})$ is the
functional calculus with commuting operator argument used in \cite{bb5}.
We conclude that: {\em if
the formal power series $F(z) \in \cL(\cU, \cY)\langle \langle z
\rangle \rangle$ is in the noncommutative Schur-Agler class $\mathcal
{S A}_{\G}(\cU, \cY)$, then its abelianization $F^{{\mathbf
a}}(\lambda)$ is in the commutative Schur-Agler class
$\mathcal{S A}_{Z^{\mathbf a}_{\G}}$ associated with
${\mathbf Q}(\lambda) : =
Z_{\G}^{{\mathbf a}}(\lambda)$ as defined in \cite{bb5}.}
Moreover, we see that the pair $(X_{L}, \del_{L})$ is left
admissible whenever $\del_{L} = (\del_{L,1}, \dots, \del_{L,d})$ is a
commutative operator-tuple in ${\mathcal B}_{\G}\cL(\cK)$, and then,
from the identity \eqref{abelianization}, we see in addition that
$$
    (X_{L}F)^{\wedge L}(\del_{L}) = (X_{L}F^{a})^{ \wedge L}(\del_{L}).
$$
More generally, if $\del_{L} = (\del_{L,1}, \dots, \del_{L,d})$ is a
commuting operator-tuple with Taylor spectrum contained in ${\mathcal
D}_{Z^{{\mathbf a}}_{\G}}$, one can use Theorem 2.1 from \cite{CH} to
see that then $\del_{L}$ is similar to a commuting operator-tuple
$\del_{L}'$ satisfying $\| Z^{{\mathbf a}}_{\G}(\del_{L}')\| < 1$, and
hence $(X_{L}, \del_{L})$ is admissible in this case as well.
We have arrived at the following result.

\begin{Pn} Suppose that $F(z) = \sum_{v \in \cF_{E}} F_{v} z^{v}$ is a formal
power series in the class $\mathcal{ S A}_{\G}(\cU, \cY)$.
\begin{enumerate}
      \item Suppose that $X_{L} \in {\mathcal L}(\cY, \cK_{L})$ and
      that $\del_{L} = (\del_{L,1}, \dots, \del_{L,d})$ is a
      commutative tuple of operators on $\cK_{L}$ with Taylor joint spectrum
      $\sigma_{\text{Taylor}}(\del_{L})$ contained in ${\mathcal
      D}_{Z^{{\mathbf a}}_{\G}}$.  Then the pair $(\del_{L}, X_{L})$ is
      left-admissible. In particular, $(\del_{L},X_{L})$ is in
      left-admissible whenever $\del_{L}$ is a commutative
      tuple in ${\mathcal B}_{\G}\cL(\cK_{L})$.

      \item Suppose that $Y_{R} \in {\mathcal L}(\cK_{R}, \cU)$ and
      that $\del_{R} = (\del_{R,1}, \dots, \del_{R,d})$ a commutative
      tuple of operators on $\cK_{R}$ with Taylor joint spectrum
      contained in ${\mathcal D}_{Z^{{\mathbf a}}_{\G}}$. Then the
      pair $(Y_{R}, \del_{R})$ is right-admissible.  In particular,
      $(Y_{R}, \del_{R})$ is right-admissible whenever $\del_{R}$ is a
      commutative tuple in ${\mathcal B}_{\G}\cL(\cK_{R})$.
   \end{enumerate}
   \end{Pn}

   \begin{proof}
       Statement (1) follows from the discussion immediately preceding
       the statement of the Proposition.  A completely parallel argument
       proves statement (2).
       \end{proof}

We now give a  sufficient condition for
left-admissibility for the general case.

\begin{Pn}  \label{P:Gadmissible}
        Suppose that $\G$ is an admissible graph, $\del = \{ \del_{e}
        \colon e \in E\}$ is a tuple of operators on the Hilbert space
        $\cK_{L}$ and that $X_{L} \in \cL(\cY, \cK_{L})$. Set
        $\rho_{_{\G,L}} = 1/N_{S}$ with $N_{S}$ defined as in \eqref{NS}.
        Then a sufficient condition for the pair $(\del_{L}, X_{L})$ to be
        left admissible with respect to the Schur-Agler class
        $\mathcal{SA}_{\G}(\cU, \cY)$ is that
      \begin{equation}
\label{Lsummable''}
          \sum_{v \in \cF_{E}} \rho_{_{\G,L}}^{-|v|} \| X_{L}^{*}
         \del_{L}^{*v} k \|_{\cY}^{2} < \infty\quad
         \text{for all } \; k \in \cK_{L}.
      \end{equation}
        \end{Pn}

        \begin{proof}
     Suppose that $H(z) = \sum_{v \in \cF_{E}} H_{v} z^{v}$ is of the
     form \eqref{1.19} in a representation \eqref{1.20} for a
     Schur-Agler class formal power series $\T(z)$ as in \eqref{1.21}.
     Then
     $$
     H(z) = C (I - Z_{\G, \cH}(z) A)^{-1}
     $$
     where $\T(z) = D + C
     (I - Z_{\G, \cH}(z) A)^{-1} Z_{\G, \cH}(z) B$
     is a unitary
     realization for $\T(z)$.  Let $x \in \oplus_{s \in S} \cH_{[s]}$.
     Then from \eqref{freqdomIO} which now takes the form
     $$\widehat y(z)  = H(z) x(\emptyset) + F(z) \cdot \widehat u(z)
     $$
     we see that the coefficients
     $H_{v} x$ of $H(z) x$ amount to the output
     string $y(v) = H_{v}x$ associated with running the SNMLS $\Col
     = (\G, \cH, \cU, \cY, \bU = \left[ \begin{smallmatrix} A & B \\ C & D
     \end{smallmatrix} \right])$ with zero input string $u(v) = 0$ for
     all $v \in \cF_{E}$ and with initial state $x(\emptyset) = x_{0}$.
     Hence, from \eqref{finalEB} we see that
     $$
       \sum_{v \in \cF_{E}} \rho_{\G, L}^{|v|}\| H_{v} x\|^{2}_{\cY}
        \le \| x_{0}\|^{2} < \infty.
     $$
     Hence
\begin{align*}
&  \sum_{n=0}^{N}
\sum_{v \in \cF_{E} \colon |v| = n} \left| \left\langle
\del_{L}^{v^{\top}} X_{L} H_{v}
x_{0}, \; k \right\rangle_{\cK_{L}} \right|  =
\sum_{v \in \cF_{E}} \left|
\left\langle \rho_{_{\G,L}}^{|v|/2} H_{v} x_{0}, \;
\rho_{_{\G,L}}^{-|v|/2}
X_{L}^{*}\del_{L}^{*v} k \right\rangle_{\cY} \right| \\
    & \qquad \le  \left( \sum_{v \in \cF_{E}} \rho_{\G,L}^{|v|} \| H_{v}
    x_{0}\|_{\cY}^{2} \right)^{1/2} \cdot \left( \sum_{v \in \cF_{E}}
    \rho_{\G, L}^{-|v|} \|X_{L}^{*} \del_{L}^{*v} k\|_{\cY}^{2}
    \right)^{1/2} < \infty.
\end{align*}
and it follows that $(\del_{L}, X_{L})$ is left-admissible as wanted.
\end{proof}

Given an admissible graph $\G$, we can always associate a new graph
     $\G^{FM}$ of Fornasini-Marchesini type (as in Example \ref{E:NCFMimp})
     by letting $\G^{FM}$ be the admissible graph of Fornasini-Marchesini
     type having the same edge set $E$ as $\G$.  This notation appears
     in the next corollary.

     \begin{Cy} \label{C:strictrowcon}
         Let $\G$ be an admissible graph with associated $\rho_{_{\G,L}} =
1/N_{S}$
         given by  \eqref{NS}, let $\del_{L} = (\del_{L,e})_{e \in E}$ be a
         tuple of operators in $\cL(\cK_{L})$ and let $X_{L} \in
         \cL(\cY, \cK_{L})$.  Then a sufficient condition for $(\del_{L},
         X_{L})$ to be admissible is that $\| Z_{\G^{FM}}(\del_{L}) \| <
         \sqrt{\rho_{\G, L}}$, i.e., that
      $$
      \| \operatorname{row}_{e \in E}  \del_{L,e}  \| < \sqrt{\rho_{\G,L}}.
      $$
         In particular, if $\G$ is a Fornasini-Marchesini
         graph (see Example \ref{E:NCFMimp}), then $(\del_{L},X_{L})$ is
         left admissible with respect to $\mathcal{ S A}_{\G}(\cU, \cY)$
         whenever $\del_{L} \in {\mathcal B}_{\G}\cL(\cK_{L})$.
         \end{Cy}

     \begin{proof}
     Set $r := \| Z_{G^{FM}}(\del_{L}) \| = \left\|
     \begin{bmatrix} T_{L,e_{1}} &
     T_{L,e_{2}} & \cdots & T_{L,e_{d}} \end{bmatrix}\right\|$.  Then
     the operator
     $$ Z_{G^{FM}}(\del)^{*} = \begin{bmatrix} T_{L,e_{1}}^{*} \\ \vdots \\
     T_{L,e_{d}}^{*} \end{bmatrix} \colon \cK \to \oplus_{e \in E} \cK
     $$
     also has norm $r$.  Hence, for each $k \in \cK_{L}$ we have
     $$ \sum_{e \in E} \|  T_{L,e}^{*} k\|^{2} \le r^{2}
      \|k\|^{2}
     $$
     and, more generally,
     $$ \sum_{v \colon |v| = N+1} \| \del_{L}^{*v}k \|^{2} \le
       r^{2}  \sum_{v \in \cF_{E} \colon |v| = N} \|
       \del_{L}^{*v} k \|^{2}.
     $$
     An easy induction argument then gives
     $$
        \sum_{v \in \cF_{E} \colon |v| = N} \| \del_{L}^{*v} k
        \|^{2}  \le r^{2N}  \|k\|^{2}
       $$
    and hence also
     $$
      \sum_{v \in \cF_{E} \colon |v| = N} \| X_{L}^{*} \del_{L}^{N} k
      \|^{2}  \le r^{2N} \| X_{L}^{*} \|^{2} \|k\|^{2}.
     $$
     Hence
     $$
     \sum_{N=0}^{\infty} \sum_{v \in \cF_{E} \colon |v| = N}
     \rho_{_{\G,L}}^{-|v|} \| X_{L}^{*} \del_{L}^{*v}k\|^{2} \le
     \| X_{L}^{*}\|^{2} \|k\|^{2} \cdot  \sum_{N=0}^{\infty} \left(
     \frac{r}{\sqrt{\rho_{_{\G,L}}}}\right)^{2N} < \infty
     $$
if $r < \sqrt{\rho_{_{\G,L}}}$. An application of the criterion
\eqref{Lsummable''} from Proposition \ref{P:Gadmissible} now completes the
proof of  Corollary  \ref{C:strictrowcon}.
\end{proof}

     Given an admissible graph $\G$ together with a tuple of
     operators $\del_{R} = (\del_{R,e})_{e \in E}$ of operators
     on a Hilbert space $\cK_{R}$ and an operator
     $Y_{R} \in \cL(\cK_{R}, \cU)$,
     there is a sufficient condition for right admissibility of the
     $(Y_{R}, \del_{R})$  in the sense of \eqref{1.25''} dual to
     condition \eqref{Lsummable''} which can be obtained as follows.
     Note that weak convergence of the series
     $ \sum_{v \in \cF_{E}} \T_{v} Y_{R} \del_{R}^{v^{\top}}$ is equivalent
     to weak convergence of the adjoint series
     $$
     \sum_{v \in \cF_{e}} \del_{R}^{* v} Y_{R}^{*} \T_{v}^{*} =
     \sum_{v \in \cF_{e}} \del_{R}^{*v^{\top}} Y_{R}^{*}
        \T_{v^{\top}}^{*}
     $$
     which has the same form as \eqref{1.25'} with $\del_{R}^{*}$ in place of
     $\del_{L}$, $Y_{R}^{*}$ in place of $X_{L}$ and $\T_{v^{\top}}^{*}$ in
     place of $\T_{v}$.  To apply the results on left-admissibility to get
     results on right admissibility, we wish to consider $(\del_{R}^{*},
     Y_{R}^{*})$ as a left pair acting on the formal power series
     $$
      \T(z)^{*} = \sum_{v \in \cF_{E}} \T_{v^{\top}}^{*} z^{v}
     $$
     in place of $\T(z) = \sum_{v \in \cF_{E}} \T_{v} z^{v}$.
     We know from Theorem \ref{T:conservative-real} that the formal power
     series
$\T(z) = \sum_{v \in \cF_{E}} \T_{v} z^{v}$ is in the Schur-Agler
     class $\mathcal{ S A}_{\G}(\cU, \cY)$ if and only if $\T(z)$ has a
     representation \eqref{1.18} with $U = \sbm{A & B \\ C & D}$ unitary.
     If $\T(z)$ has the form \eqref{1.18}, then we compute
     \begin{align*}
     \T(z)^{*} & = \sum_{v \in \cF_{E}} \T_{v}^{*} z^{v^{\top}} \\
      & = D^{*} + B^{*} Z_{\G, \cH}(z)^{*}
      (I - A^{*} Z_{\G, \cH}(z)^{*})^{-1} C^{*} \\
      & = D^{*} + B^{*} (I - Z_{\G, \cH}(z)^{*} A^{*})^{-1} Z_{\G,
      \cH}(z)^{*} C^{*}.
     \end{align*}
     This suggests that, given a SNMLS
     $\boldsymbol{\Sigma} = (\G,\cH, \cU, \cY, \bU)$ as defined in
     \eqref{3.1},
     we define a dual SNMLS $\boldsymbol{\Sigma}^{\prime} =
     (\G', \cH, \cY, \cU, \bU')$ where
     \begin{enumerate}
         \item the admissible graph $\G'$ for $\boldsymbol{\Sigma}'$
         is the same graph as
         the admissible graph $\G$, but with the source vertices
         for $\G$ taken to be the range vertices for $\G'$ and with the
         range vertices for $\G$ taken to be the source vertices for $\G'$;
         thus the set of path-components remains unchanged: $P' = P$,
          and
         \item the connection matrix $\bU'$ for $\boldsymbol{\Sigma}'$ is
         simply the adjoint
         $$ \bU' = \bU^{*} = \begin{bmatrix} A^{*} & C^{*} \\ B^{*} & D^{*}
         \end{bmatrix} \colon \begin{bmatrix} \oplus_{r \in R} \cH_{[r]}
         \\ \cY \end{bmatrix} \to \begin{bmatrix} \oplus_{s \in S}
         \cH_{[s]} \\ \cU \end{bmatrix}
         $$
         of the connection matrix $U$ for $\boldsymbol{\Sigma}$.
     \end{enumerate}
     Then it is easily checked: {\em if $\T(z)$ is the transfer
     function of the SNMLS $\boldsymbol{\Sigma}$, then $\T(z)^{*}$ is
     the transfer function of the SNMLS $\boldsymbol{\Sigma}'$.}
     Moreover $\boldsymbol{\Sigma}$ is conservative (i.e., $\bU$ is unitary)
     if and only
     if $\boldsymbol{\Sigma}'$ is conservative (i.e., $\bU' = \bU^{*}$
is unitary).
     By the equivalence (1) $\Longleftrightarrow$ (2) in Theorem
     \ref{T:conservative-real}, we conclude that: {\em the formal power
     series $\T(z) = \sum_{v \in \cF_{E}} \T_{v} z^{v}$ is in the
     Schur-Agler class $\mathcal{ S A}_{\G}(\cU, \cY)$ if and only if its
     adjoint $\T(z)^{*} = \sum_{v \in \cF_{E}} \T_{v^{\top}}^{*} z^{v}$ is
     in the Schur-Agler class $\mathcal{ S A}_{\G'}(\cY, \cU)$, where
     $\G'$ is the reflection of $\G$ induced by interchanging source
     vertices with range vertices.}

     A consequence of this analysis is that we have the following
     analogues of Proposition \ref{P:Gadmissible} and Corollary
     \ref{C:strictrowcon}.  We leave the details
     of the proof to the reader.  In the statement of the theorem we use
     the notation
     \begin{equation}  \label{NR}
         n_{R} = \max\{ n_{r_{p}} \colon p \in P\}
      \end{equation}
      where $n_{r_{p}}$ is the number of range vertices in component $P$
      of the graph $G$.

     \begin{Pn}  \label{P:Gadmissible'}
         Let $\G$ be an admissible graph with associated constant
         $\rho_{\G,R}:= 1/n_{R}$ with $n_{R}$ as in \eqref{NR}, let
         $\del_{R} = (\del_{R,e})_{e \in E}$ be a tuple of operators
         acting on a Hilbert space $\cK_{R}$ and let $Y_{R} \in
         \cL(\cK_{R}, \cU)$.  Then a sufficient condition for the pair
         $(Y_{R}, \del_{R})$ to be right admissible in the sense of
         \eqref{1.25''} is that
      \begin{equation}
      \label{Rsummable''}
         \sum_{v \in \cF_{E}} \rho_{\G,R}^{-|v|}
         \| Y_{R} \del_{R}^{v} k \|^{2} < \infty \quad \text{for all } \;
         k \in \cK_{R}.
      \end{equation}
      \end{Pn}

       For the statement of the following corollary, we use the notation
       $\G^{FM'}$ to denote the dual of the Fornasini-Marchesini graph
       $\G^{FM}$ associated with $\G$; thus $\G^{FM^\prime}$ has a single
       range vertex $\{r_{0}\}$, the same edge set $E$ as does $\G$ and the
       source-vertex set taken also equal to $E$ and with each edge
        $e$ considered to have source itself $e$ and range $r_{0}$.
       The associated structure matrix $Z_{\G}(z)$ is then a column
       $$
       Z_{\G^{FM^\prime}}(z) = \begin{bmatrix} z_{e_{1}} \\ \vdots \\
       z_{e_{d}} \end{bmatrix}
       $$
       where $d = n_{E}$ is the number of edges.

       \begin{Cy} \label{C:strictrowcon'}
       Let $\G$ be an admissible graph with associated $\rho_{\G,R} = 1/N_{R}$
      given by  \eqref{NR}, let $\del_{R} = (\del_{R,e})_{e \in E}$ be a
      tuple of operators in $\cL(\cK_{R})$ and let $X_{R} \in
      \cL(\cK_{R}, \cU)$.  Then a sufficient condition for $(Y_{R},
      \del_{R})$ to be right-admissible is that $\| Z_{G^{FM
      \prime}}(\del_{R}) \| <\sqrt{\rho_{\G, R}}$.
      In particular, if $\G = \G^{FM \prime}$ is itself the reflection of a
      Fornasini-Marchesini graph, then $(Y_{R}, \del_{R})$ is right
      admissible whenever $\del_{R} \in {\mathcal B}_{\G}\cL(\cK)$.
       \end{Cy}

\section{The solvability criterion} \label{S:solvability}
\setcounter{equation}{0}

In this section we prove the necessity part of Theorem \ref{T:1.9}.
First we need to note the following elementary properties of evaluations
\eqref{1.25'} and \eqref{1.25''}.
\begin{La}
\label{L:2.1}
Let $\del=\{\del_e\colon e\in E\}$ and
$\del^\prime=\{\del^\prime_e\colon e\in E\}$ be tuples of bounded
linear operators acting on Hilbert spaces $\cK$ and $\cK^\prime$,
respectively.
\begin{enumerate}

\item For every constant function $W(z)\equiv W\in\cL(\cK^\prime,\cK)$,
\begin{equation}
\left(W\right)^{\wedge L}(\del)=\left(W\right)^{\wedge
R}(\del^\prime)=W.
\label{2.1}
\end{equation}

\item For every $\T\in\cL(\cU,\cK)\langle\langle z\rangle\rangle$,
and $\widetilde{\T} \in\cL(\cK^\prime,\cY)\langle\langle
z\rangle\rangle$, $W\in \cL(\cU',\cU)$ and $\widetilde{W}\in\cL(\cY,\cY')$.
\begin{equation}
\left(\T\cdot W\right)^{\wedge L}(\del)=\T^{\wedge L}(\del)\cdot W
\quad\mbox{and}\quad
\left(\widetilde{W}\cdot\widetilde{\T}\right)^{\wedge R}(\del^\prime)
=\widetilde{W} \cdot \widetilde{\T}^{\wedge R}(\del^\prime)
\label{2.2}
\end{equation}
whenever $\T^{\wedge L}(\del)$ and $\widetilde{\T}^{\wedge
R}(\del^{\prime})$
are defined.

\item For every $\T$ and $\widetilde{\T}$ as in part (2) and every
$e\in E$,
\begin{equation}
\left(\T(z)z_e\right)^{\wedge L}(\del)=\del_e\cdot \T^{\wedge
L}(\del)\quad\mbox{and}\quad
\left(z_e\widetilde{\T}(z)\right)^{\wedge R}(\del^\prime)
=\widetilde{\T}^{\wedge R}(\del^\prime)\cdot \del^\prime_e
\label{2.3}
\end{equation}
whenever $\T^{\wedge L}(\del)$ and $\widetilde{\T}^{\wedge
R}(\del^{\prime})$ are defined.

\item For every choice of $\T\in\cL(\cU,\cK)\langle\langle
z\rangle\rangle$ and of
$\widetilde{\T}\in\cL(\cU^\prime,\cU)\langle\langle z\rangle\rangle$
\begin{equation}
\left(\T\cdot \widetilde{\T} \right)^{\wedge L}(\del)=(\T^{\wedge
L}(\del)\cdot \widetilde{\T})^{\wedge L}(\del)
\label{2.5}
\end{equation}
whenever $\T^{\wedge L}(\del)$ and $(\T^{\wedge
L}(\del)\cdot \widetilde{\T})^{\wedge L}(\del)$ are defined.

\item For every choice of
$\T\in\cL(\cY^\prime,\cY)\langle\langle z\rangle\rangle$ and of
$\widetilde{\T}\in\cL(\cK^\prime,\cY)\langle\langle z\rangle\rangle$,
\begin{equation}
\left(\T\cdot \widetilde{\T}\right)^{\wedge R}(\del^\prime)=(\T\cdot
\widetilde{\T}^{\wedge R}(\del^\prime))^{\wedge R}(\del^\prime)
\label{2.6}
\end{equation}
whenever $\widetilde{\T}^{\wedge R}(\del^\prime)$ and
$(\T\cdot
\widetilde{\T}^{\wedge R}(\del^\prime))^{\wedge R}(\del^\prime)$ are
defined.
\end{enumerate}
\label{lm2.1}
\end{La}

{\bf Proof:} The two first statements follow immediately from definitions
\eqref{1.25'} and \eqref{1.25''}. To prove \eqref{2.5}, take $\T$ and
$\widetilde{\T}$ in the
form
$$
\T(z) = \sum_{v\in\cF_{E}} \T_{v} z^{v},\quad
\widetilde{\T}(z)=\sum_{v \in \cF_{E}} \widetilde{\T}_{v} z^{v}.
$$
Then
$$
\T(z)\cdot\widetilde{\T}(z)=\sum_{v\in\cF_{E}} \left(\sum_{uw=v}
\T_{u}\widetilde{\T}_w\right)z^{v}
$$
and therefore, according to \eqref{1.25'},
\begin{equation}
\left(\T\cdot \widetilde{\T} \right)^{\wedge L}(\del)
=\sum_{v\in\cF_E}\del^{v^\top}\left(\sum_{uw=v}
\T_{u}\widetilde{\T}_w\right).
\label{2.6a}
\end{equation}
On the other hand, again by \eqref{1.25'},
\begin{eqnarray*}
(\T^{\wedge L}(\del)\cdot \widetilde{\T})^{\wedge L}(\del)
&=&\sum_{w\in\cF_E}\del^{v^\top}
\T^{\wedge L}(\del) \widetilde{\T}_w\\
&=&\sum_{w\in\cF_E}\del^{w^\top}\left(
\sum_{u\in\cF_E}\del^{u^\top}\T_{u}\right)\widetilde{\T}_w\\
&=&\sum_{w,u\in\cF_E}\del^{(uw)^\top}\T_{u}\widetilde{\T}_w\\
&=&\sum_{v\in\cF_E}\del^{v^\top}\left(\sum_{uw=v}
\T_{u}\widetilde{\T}_w\right).
\end{eqnarray*}
Comparison of the last equality with \eqref{2.6a} gives \eqref{2.5}.
Equality \eqref{2.6} is obtained in much the same way. The first
equality in \eqref{2.3} follows from \eqref{2.5} for the special case of
$\widetilde{\T}(z)=z_e I_\cU$. The second equality in \eqref{2.3} follows
from \eqref{2.6} for the special case of $\T(z)=z_e I_\cY$.\qed

\medskip

{\bf Proof of the necessity part in Theorem \ref{T:1.9} and \ref{T:1.11}:}
Let $\T$ belong to ${\mathcal {SA}}_\G(\cU, \, \cY)$ and suppose
that $\T$ is a solution of Problem \ref{P:1.4}.  Choose formal power
series $H$ and $G$ of the form \eqref{1.19} and \eqref{1.21}
so that the representations \eqref{1.20}, \eqref{1.22}, \eqref{1.23}
hold.  Use \eqref{1.37}--\eqref{1.39} to define operators
$\Psi_{s,s'}$, $\Lambda_{s,r}$ and $\Phi_{r,r'}$ for $s,s' \in S$ and
$r,r' \in R$.  Then use equations  \eqref{1.34}--\eqref{1.36} to
define the block operator matrix $\K_{p}$.  If $\T$ is assumed to be
a solution of Problem \ref{P:1.10} then we are given $\K_{p}$ via
\eqref{1.34}--\eqref{1.36} where $\Psi_{s,s'}$, $\Lambda_{s,r}$ and
$\phi_{r,r'}$ are part of the interpolation data and
\eqref{1.37}--\eqref{1.39} hold as part of the interpolation
conditions for some choice of $H$ and $G$ associated with the
representations \eqref{1.20}, \eqref{1.22}, \eqref{1.23} for $\T$.
In any case, the conditions
\eqref{1.37}--\eqref{1.39} hold and imply that $\K_p$ can be represented as
\begin{equation}
\K_p=\begin{bmatrix}{\mathbb T}_{p,L}^* \\ {\mathbb
T}_{p,R}^*\end{bmatrix}
\begin{bmatrix}{\mathbb T}_{p,L} & {\mathbb T}_{p,R} \end{bmatrix}
\label{2.7}
\end{equation}
where the operators ${\mathbb T}_{p,L}$ and ${\mathbb T}_{p,R}$
are given by
\begin{eqnarray}
{\mathbb T}_{p,L}&=&\operatorname{Row}_{s\in S\colon
[s]=p}\left[(X_LH_s)^{\wedge L}(\del_L)\right]^*\colon \;
\bigoplus_{s\in S\colon [s]=p}\cK_L\to \cH_p,\label{2.8} \\
{\mathbb T}_{p,R}&=& \operatorname{Row}_{r\in R\colon
[r]=p}(G_rY_R)^{\wedge
R}(\del_R)\colon \; \bigoplus_{r\in R\colon [r]=p}\cK_R\to \cH_p.
\label{2.9}
\end{eqnarray}
Comparing \eqref{2.7} with \eqref{1.34} we see that
\begin{equation}
\K_{p,L}={\mathbb T}_{p,L}^*{\mathbb T}_{p,L},\quad \K_{p,R}=
{\mathbb T}_{p,R}^*{\mathbb T}_{p,R}, \quad \K_{p,LR}
={\mathbb T}_{p,L}^*{\mathbb T}_{p,R}.
\label{2.9a}
\end{equation}
It follows from \eqref{2.7} that $\K_p\ge 0$ and thus, it remains to show
that these operators satisfy the Stein identity \eqref{1.32}. To
this end, note that by \eqref{1.17} and \eqref{1.19},
$$
H(z)Z_{\G,\cH}(z)=\operatorname{Row}_{r\in R}\sum_{s\in S\colon [s]=[r]}
H_s(z)z_{e_{s,r}}
$$
and therefore, by the first equality in \eqref{2.3},
$$
\left(X_LHZ_{\G,\cH}\right)^{\wedge L}(\del_L)=
\operatorname{Row}_{r\in R}\sum_{s\in S\colon [s]=[r]}
\del_{L,e_{s,r}}\left(X_LH_s\right)^{\wedge L}(\del_L)
$$
which can be written in terms of \eqref{1.29} and \eqref{2.8} as
\begin{equation}
\left(X_LHZ_{\G,\cH}\right)^{\wedge L}(\del_L)=
\operatorname{Row}_{r\in R}\widetilde{N}_r(\del_L)^*{\mathbb T}_{[r],L}^*.
\label{2.10}
\end{equation}
Note also that according to decompositions \eqref{1.19} and \eqref{1.28},
\begin{equation}
\left(X_LH\right)^{\wedge L}(\del_L)=\operatorname{Row}_{s\in S}
E_{L,s}^*{\mathbb T}_{[s],L}^*.
\label{2.10a}
\end{equation}
Similarly, by \eqref{1.17} and \eqref{1.21},
$$
Z_{\G,\cH}(z)G(z)=\operatorname{Col}_{s\in S}\sum_{r\in R\colon [r]=[s]}
z_{e_{s,r}}G_r(z)
$$
and therefore, by the second equality in \eqref{2.3},
$$
\left(Z_{\G,\cH}GY_R\right)^{\wedge R}(\del_R)=
\operatorname{Col}_{s\in S}\sum_{r\in R\colon [r]=[s]}
\left(G_rY_R\right)^{\wedge R}(\del_R)\del_{R,e_{s,r}},
$$
which can be written in terms of \eqref{1.29a} and
\eqref{2.9} as
\begin{equation}
\left(Z_{\G,\cH}GY_R\right)^{\wedge R}(\del_R)=
\operatorname{Col}_{s\in S}{\mathbb T}_{[s],R}\widetilde{M}_s(\del_R).
\label{2.11}
\end{equation}
Finally, by decompositions \eqref{1.21} and \eqref{1.28a},
\begin{equation}
\left(GY_R\right)^{\wedge R}(\del_R)=\operatorname{Col}_{r\in R}{\mathbb
T}_{[r],R}E_{R,r}.
\label{2.11a}
\end{equation}
Substituting the partitionings \eqref{1.30}, \eqref{1.31}, \eqref{1.33}
and \eqref{1.34}  into \eqref{1.32} we conclude that \eqref{1.32} is
equivalent to the following three equalities:
\begin{eqnarray}
\sum_{s\in S} E_{L,s}^*\K_{[s],L}E_{L,s}-\sum_{r\in R}
\widetilde{N}_r(\del_L)^*\K_{[r],L}\widetilde{N}_r(\del_L)
&=&X_LX_L^*-Y_LY_L^*,\label{2.12}\\
\sum_{s\in S} E_{L,s}^*\K_{[s],LR}\widetilde{M}_s(\del_R)
-\sum_{r\in R} \widetilde{N}_r(\del_L)^*\K_{[r],LR}E_{R,r}
&=&X_LX_R-Y_LY_R,\label{2.13}\\
\sum_{s\in S} \widetilde{M}_s(\del_R)^*\K_{[s],R}\widetilde{M}_s(\del_R)
-\sum_{r\in R} E_{R,r}^*\K_{[r],R}E_{R,r}
&=&X_R^*X_R-Y_R^*Y_R.\label{2.14}
\end{eqnarray}
To check (\ref{2.12}) we consider the equality
\begin{equation}  \label{equality}
X_LX_L^*-X_L\T(z)\T(z^\prime)^*X_L^*=X_LH(z)\left(I-
Z_{\G,\cH}(z)Z_{\G, \cH}(z')^{*} \right) H(z')^{*}X_L^*
\end{equation}
which is an immediate corollary of \eqref{1.20}. We may consider each
side of \eqref{equality} as a formal power series in $z'$ with coefficients
equal to formal power series in $z$, i.e., we have a natural
identification
$$ \cL({\mathcal K}_{L})\langle \langle z, z' \rangle  \rangle
       \cong  \left(\cL({\mathcal K}_{L})\langle \langle z \rangle \rangle
\right) \langle
       \langle z' \rangle \rangle.
$$
We then apply the left evaluation map (applied to formal power series
in the variable $z$) to each coefficient of the resulting formal power
series in the variable $z'$.  The result amounts to applying left
evaluation to both sides of \eqref{equality} in the variable $z$ with
the formal variable $z'$ considered as fixed.
Making use of
properties \eqref{2.1}, \eqref{2.2} and of relation (\ref{2.10}) and
taking into account the first interpolation condition in \eqref{1.27}, we
get
\begin{eqnarray*}
X_LX_L^*-Y_L  F(z')^*X_L^*&=&\left(X_LH\right)^{\wedge L}(\del_L)\cdot
H(z')^*X_L^*\\
&&-\left(X_LHZ_{\G, \cH}\right)^{\wedge L}(\del_L)\cdot
Z_{\G, \cH}(z')^{*}H(z')^*X_L^*.
\end{eqnarray*}
This equality holds as an identity in $\cL(\cK_{L})\langle \langle
z'\rangle \rangle$.
Taking adjoints  and replacing $z'$ by $z$, we get
\begin{eqnarray*}
X_LX_L^*-X_L\T(z)Y_L^*&=&X_LH(z)\left(\left(X_LH\right)^{\wedge
L}(\del_L)\right)^*\\ && -X_LH(z)Z_{\G,\cH}(z)
\left(\left(X_LHZ_{\G,\cH}\right)^{\wedge L}(\del_L)\right)^*.
\end{eqnarray*}
Applying again the left evaluation to the latter equality we get
\begin{eqnarray*}
X_LX_L^*-Y_LY_L^*&=&\left(X_LH\right)^{\wedge
L}(\del_L)\left(\left(X_LH\right)^{\wedge L}(\del_L)\right)^*\\
&&-\left(X_LH Z_{\G,\cH}\right)^{\wedge
L}(\del_L)\left(\left(X_LHZ_{\G,\cH}\right)^{\wedge L}(\del_L)\right)^*.
\end{eqnarray*}
Substituting \eqref{2.10} and \eqref{2.10a} into the right hand side
expression we come to
$$
X_LX_L^*-Y_LY_L^*=\sum_{s\in S} E_{L,s}^*{\mathbb
T}_{[s],L}^*{\mathbb T}_{[s],L}E_{L,s}-\sum_{r\in
R}\widetilde{N}_r(\del_L)^*
{\mathbb T}_{[s],L}^*{\mathbb T}_{[s],L}\widetilde{N}_r(\del_L)
$$
which is equivalent to \eqref{2.12}, since
$$
{\mathbb T}_{[s],L}^*{\mathbb T}_{[s],L}=\K_{[s],L}\quad\mbox{and}\quad
{\mathbb T}_{[r],L}^*{\mathbb T}_{[r],L}=\K_{[r],L}.
$$

\medskip

To prove \eqref{2.13} we start with equality
$$
X_L\T(z)Y_R-X_L\T(z')Y_R=X_LH(z)\left(Z_{\G,\cH}(z)-Z_{\G,\cH}(z')\right)G(z')
Y_R
$$
which is a consequence of \eqref{1.23}.
We apply the left evaluation in the $z$ variable:
by the first interpolation condition in
\eqref{1.27} we have
$$
Y_LY_R-X_L\T(z')Y_R=\left(X_LHZ_{\G,\cH}\right)^{\wedge
L}(\del_L)G(z')Y_R-\left(X_LH\right)^{\wedge L}(\del_L)
Z_{\G,\cH}(z')G(z')Y_R.
$$
The last identity equality holds true as an identity between formal
power series in the variable $z'$; we then apply the right
evaluation \eqref{1.25''} to both sides. In view of the second interpolation
condition in \eqref{1.27} and of properties \eqref{2.1}, \eqref{2.2}, we
obtain
\begin{eqnarray*}
Y_LY_R-X_LX_R&=&\left(X_LHZ_{\G,\cH}\right)^{\wedge
L}(\del_L)\left(GY_R\right)^{\wedge R}(\del_R)\\
&&-\left(X_LH\right)^{\wedge
L}(\del_L)\left(Z_{\G,\cH}GY_R\right)^{\wedge R}(\del_R).
\end{eqnarray*}
Substituting equalities \eqref{2.10}, \eqref{2.10a}, \eqref{2.11} and
\eqref{2.11a} into the right-hand side expression in the last equality we
come to
$$
X_{L}X_{R} - Y_{L}Y_{R} =
      \sum_{s\in S} E_{L,s}^*{\mathbb
T}_{[s],L}^*{\mathbb T}_{[s],R}\widetilde{M}_s(\del_R)
-\sum_{r\in R} \widetilde{N}_r(\del_L)^*{\mathbb T}_{[r],L}^*{\mathbb
T}_{[r],R}E_{R,r}
$$
which is equivalent to \eqref{2.13}, since
$$
{\mathbb T}_{[s],L}^*{\mathbb T}_{[s],R}=\K_{[s],LR}\quad\mbox{and}\quad
{\mathbb T}_{[r],L}^*{\mathbb T}_{[r],R}=\K_{[r],LR},
$$
by \eqref{2.9a}. The proof of \eqref{2.14} is quite similar:
we start with the equality
$$
Y_R^*Y_R-Y_R^*\T(z)^*\T(z')Y_R=Y_R^*G(z)^*\left(I-
Z_{\G,\cH}(z)^*Z_{\G,\cH}(z')\right)G(z')Y_R
$$
(which follows from \eqref{1.22}) and apply the right evaluation in
the $z'$ variable.
     Then we take adjoints in the resulting formal power series identity
(in the variable $z$) and apply again the right evaluation
map. The obtained equality together with relations (\ref{2.11}) and
(\ref{2.11a}) leads to (\ref{2.14}).
This completes the proof of necessity in both Theorem \ref{T:1.9} and
Theorem \ref{T:1.11}.

\medskip

\section{Solutions to the interpolation problem and unitary
extensions} \label{S:solutions}
\setcounter{equation}{0}

In this Section we shall show that there is a correspondence between
solutions to Problem \ref{P:1.10} and unitary extensions of a
partially defined isometry determined by the problem data set
${\mathcal D}$.

    From now on we assume that we are given an interpolation data set
${\mathcal D}$ as in \eqref{1.40} and that the necessary conditions
for Problem \ref{P:1.10} to have a solution are in force:  the
operators $\K_p$ defined in \eqref{1.34}, \eqref{1.35} are each positive
semidefinite on the space
\begin{equation}
\cH^\circ_p=
\left(\bigoplus_{s\in S\colon [s]=p}\cK_L\right)
\oplus \left(\bigoplus_{r\in R\colon [r]=p}\cK_R \right)
\label{3.6a}
\end{equation}
and satisfy the Stein identity \eqref{1.32} which we write now as
\begin{equation}
\sum_{s\in S}M_s^*\K_{[s]}M_s+Y^*Y= \sum_{r\in R}N_r^*\K_{[r]}N_r+X^*X.
\label{3.7}
\end{equation}
For every $p\in P$, we introduce the equivalence $\buildrel{p}\over \sim$
on $\cH^\circ_p$ by
$$
h_1 \buildrel{p}\over \sim h_2 \; \; \mbox{if and only if} \; \; \langle
\K_p(h_1-h_2), \; y\rangle_{\cH^\circ_p}=0\; \; \mbox{for all} \; \;
y\in\cH^\circ_p,
$$
denote $[h]_p$ the equivalence class of $h$ with respect to the above
equivalence and endow the linear space of equivalence
classes  with the inner product
\begin{equation}
\langle [h]_p, \; [y]_p\rangle=\langle \K_p h, \; y\rangle_{\cH^\circ_p}.
\label{3.8}
\end{equation}
We get  a  prehilbert space whose completion is $\hH_p$.
It is readily seen from definitions \eqref{1.30}, \eqref{1.31} of
operators $M_s$ and $N_r$ that $M_sf$ and $N_rf$ belong to
$\cH^\circ_{[s]}$ and $\cH^\circ_{[r]}$, respectively, for every
choice of $f\in\cK_L\oplus\cK_R$. Furthermore, identity \eqref{3.7} can
be written as
$$
\sum_{s\in S}\left\|[M_sf]_{[s]}\right\|_{\hH_{[s]}}^2+\left\|
Yf\right\|_{\cU}^2=
\sum_{r\in R}\left\| [{N_r}f]_{[r]}\right\|_{\hH_{[r]}}^2+
\left\| Xf\right\|_{\cY}^2,
$$
holding for every choice of $f\in\cK_L\oplus\cK_R$.
Therefore the linear map defined by the rule
\begin{equation}
{\bf V}: \; \left[\begin{array}{c} \operatorname{Col}_{s\in S}[M_sf]_{[s]}\\
Yf\end{array}\right]\rightarrow \begin{bmatrix}
\operatorname{Col}_{r\in R}[N_rf]_{[r]} \\ Xf\end{bmatrix}
\label{3.9}
\end{equation}
extends by linearity to define an isometry from
\begin{equation}
{\mathcal {D}}_{\bf V}={\rm Clos}
\left\{\left[\begin{array}{c}\operatorname{Col}_{s\in S}[M_sf]_{[s]}
\\ Yf\end{array}\right], \; f\in\cK_L\oplus\cK_R\right\}
\subset\left[\begin{array}{c}\oplus_{s\in S}\hH_{[s]}\\ \cU\end{array}\right]
\label{3.10}
\end{equation}
onto
\begin{equation}
{\mathcal {R}}_{\bf V}={\rm Clos}
\left\{\left[\begin{array}{c}\operatorname{Col}_{r\in R}[N_rf]_{[r]}
\\ Xf\end{array}\right], \; f\in\cK_L\oplus\cK_R\right\}
\subset\left[\begin{array}{c}\oplus_{r\in R}\hH_{[r]}\\ \cY\end{array}\right].
\label{3.11}
\end{equation}
The next two lemmas establish a correspondence between solutions $\T$
to Problem \ref{P:1.10} and unitary extensions of the partially defined
isometry ${\bf V}$ given in (\ref{3.9}).
\begin{La}
Any solution $\T$ to Problem \ref{P:1.10} is a characteristic function of a
unitary colligation
\begin{equation}
\widetilde{\Col}=\{\G, \, \tH, \, \cU, \, \cY, \widetilde{\bU}\}
\label{3.11a}
\end{equation}
with the state space
$$
\tH=\hH\oplus \cH':=\{\tH_p=\hH_p\oplus \cH'_p\colon \, p\in P\}
$$
and the connecting operator
\begin{equation}
\widetilde{\bU}=\begin{bmatrix} \widetilde{A} & \widetilde{B} \\
\widetilde{C} & \widetilde{D} \end{bmatrix}\colon
\begin{bmatrix}\oplus_{s\in S}(\hH_{[s]}\oplus\tH_{[s]}) \\
\cU  \end{bmatrix}\to\begin{bmatrix}\oplus_{r\in R}(\hH_{[r]}\oplus\tH_{[r]})
\\ \cY  \end{bmatrix}
\label{3.12}
\end{equation}
being an extension of the isometry ${\bf V}$ given in \eqref{3.9}.
\label{L:3.3}
\end{La}
{\bf Proof:} Let $\T$ be a solution to Problem \ref{P:1.10}. In
particular, $\T$ belongs to the noncommutative Schur-Agler class
${\mathcal {SA}}_\G(\cU, \, \cY)$ and, by Theorem
\ref{T:conservative-real}, it is the characteristic function of some
unitary colligation $\Col$ of the form (\ref{3.1}). In other words, $\T$
admits a unitary realization (\ref{1.18}) with the state space
${\cH}=\{\cH_p\colon p\in P\}$ and representations (\ref{1.20}),
(\ref{1.22}), (\ref{1.23}) hold for power series $H$ and $G$ defined via
(\ref{1.24}) and decomposed as in (\ref{1.19}) and (\ref{1.21}). These
series lead to the following two representations
\begin{equation}
\T(z)=D+H(z)Z_{\G,\cH}(z)B=D+CZ_{\G,\cH}(z)G(z),
\label{3.13}
\end{equation}
of $\T$, each of which is equivalent to (\ref{1.18}).

\smallskip

The interpolation conditions (\ref{1.27}) and (\ref{1.37})--(\ref{1.39})
which hold for $\T$ by assumption force certain restrictions on
the connecting operator ${\bf U}=\begin{bmatrix} A & B \\ C & D \end{bmatrix}$.
Substituting (\ref{3.13}) into (\ref{1.27}) we get equalities
$$
\left(X_LD+X_LHZ_{\G,\cH} B\right)^{\wedge L}(\del_L)=Y_L
$$
and
$$
\left(DY_R+CZ_{\G,\cH} GY_R\right)^{\wedge R}(\del_R)=X_R
$$
which are equivalent, due to properties \eqref{2.1}, \eqref{2.2}, to
\begin{equation}
X_LD+\left(X_LHZ_{\G,\cH} \right)^{\wedge L}(\del_L)B=Y_L
\label{3.14}
\end{equation}
and
\begin{equation}
DY_R+C\left(Z_{\G,\cH} GY_R\right)^{\wedge R}(\del_R)=X_R,
\label{3.15}
\end{equation}
respectively. It also follows from (\ref{1.24}) that
$$
C+H(z)Z_{\G,\cH}(z)A=H(z), \qquad B+AZ_{\G,\cH}(z)G(z)=G(z)
$$
and therefore, that
\begin{equation}
X_LC+\left(X_LHZ_{\G,\cH} \right)^{\wedge L}(\del_L)A=
\left(X_LH \right)^{\wedge L}(\del_L)
\label{3.16}
\end{equation}
and
\begin{equation}
BY_R+A\left(Z_{\G,\cH} GY_R\right)^{\wedge R}(\del_R)=\left(
GY_R\right)^{\wedge R}(\del_R).
\label{3.17}
\end{equation}
The equalities (\ref{3.14}) and (\ref{3.16}) can be written in matrix
form as
\begin{equation}
\begin{bmatrix}\left(X_LHZ_{\G,\cH} \right)^{\wedge L}(\del_L) &
X_L\end{bmatrix}
\left[\begin{array}{cc} A & B \\ C & D\end{array}\right]=
\begin{bmatrix}\left(X_LH \right)^{\wedge L}(\del_L)& Y_L\end{bmatrix},
\label{3.18}
\end{equation}
whereas the equalities (\ref{3.15}) and (\ref{3.17}) are equivalent to
\begin{equation}
\left[\begin{array}{cc} A & B \\ C & D\end{array}\right]
\begin{bmatrix}\left(Z_{\G,\cH} GY_R\right)^{\wedge R}(\del_R) \\
Y_R\end{bmatrix}=\begin{bmatrix}\left(
GY_R\right)^{\wedge R}(\del_R)\\ X_R\end{bmatrix}.
\label{3.19}
\end{equation}
Since the operator $\begin{bmatrix} A & B \\ C & D \end{bmatrix}$ is
unitary, we conclude from (\ref{3.18}) that
\begin{equation}
\left[\begin{array}{cc} A & B \\ C & D\end{array}\right]
\begin{bmatrix} \left[\left(X_LH \right)^{\wedge L}(\del_L)\right]^* \\
Y_L^*\end{bmatrix}
=\begin{bmatrix}\left[\left(X_LHZ_{\G,\cH} \right)^{\wedge
L}(\del_L)\right]^* \\
X_L^*\end{bmatrix}.
\label{3.20}
\end{equation}
Combining \eqref{3.19} and \eqref{3.20} we conclude that for every choice
of $f\in\cK_L\oplus\cK_R$,
\begin{eqnarray}
&&\left[\begin{array}{cc}A & B \\ C&
D\end{array}\right]\left[\begin{array}{cc}
\left[\left(X_LH \right)^{\wedge L}(\del_L)\right]^* & \left(Z_{\G,\cH}
GY_R\right)^{\wedge
R}(\del_R) \\ Y_L^* & Y_R\end{array}\right]f\nonumber\\
&&=\left[\begin{array}{cc}
\left[\left(X_LHZ_{\G,\cH} \right)^{\wedge L}(\del_L)\right]^* & \left(
GY_R\right)^{\wedge R}(\del_R) \\ X_L^* & X_R\end{array}\right]f.
\label{3.23}
\end{eqnarray}
Let ${\mathbb T}_{p,L}$ and ${\mathbb T}_{p,R}$ be the operators given by
\eqref{2.8} and \eqref{2.9}, respectively, and let
\begin{equation}
{\mathbb T}_p:=\begin{bmatrix} {\mathbb T}_{p,L} & {\mathbb
T}_{p,R}\end{bmatrix}: \; \cH^\circ_p \to  \cH_p.
\label{3.21}
\end{equation}
Now we use the interpolation conditions (\ref{1.37})--(\ref{1.39}),
which provide the factorization (\ref{2.7}) of the operator $\K_p$. Thus,
$$
\K_p={\mathbb T}_p^*{\mathbb T}_p
$$
and
$$
\langle [h]_p, \; [y]_p\rangle_{\hH_p}=\langle \K_p h, \;
y\rangle_{\cH^\circ_p}
=\langle {\mathbb T}_ph, \; {\mathbb T}_py\rangle_{\cH_{p}}
$$
for every $h, \, y\in\cH^\circ_p$. Therefore, the linear transformation
$U_p$
defined by the rule
\begin{equation}
U_p \, : \; \; {\mathbb T}_p h\rightarrow [h]_p\qquad (h\in\cH^\circ_p)
\label{3.24}
\end{equation}
can be extended to the unitary map (which still is denoted by $U_p$) from
$\overline{{\rm Ran} \, {\mathbb T}_p}$ onto $\hH_p$. Noticing that
$\overline{{\rm Ran} \, {\mathbb T}_p}$ is a subspace of $\cH_p$ and setting
$$
{\cN}_p:={\cH}_p\ominus\overline{{\rm Ran} \, {\mathbb T}_p}\quad {\rm
and}\quad\tH_p:=\hH_p\oplus{\cN}_p,
$$
we define the unitary map $\tU_p \, : \; \; {\cH}_p\rightarrow\tH_p$ by
the rule
\begin{equation}
\tU_p \, g=\left\{\begin{array}{cll} U_p g & {\rm for} &
g\in\overline{{\rm
Ran} \,  {\mathbb T}_p}, \\ g & {\rm for} & g\in{\cN}_p.\end{array}\right.
\label{3.25}
\end{equation}
Introducing the operators
\begin{equation}
\widetilde{A}=\left[\widetilde{U}_{[r]}A_{r,s}\widetilde{U}_{[s]}^*
\right]_{r\in R, s\in S},\quad \widetilde{B}=\operatorname{Col}_{r\in
R}\widetilde{U}_{[r]}B_r,\quad
\widetilde{C}=\operatorname{Row}_{s\in S}C_s\widetilde{U}_{[s]}^*,\quad
\widetilde{D}=D,
\label{3.24a}
\end{equation}
we construct the colligation $\widetilde{\Col}$ via (\ref{3.4}) and
\eqref{3.12}. By definition, $\widetilde{\Col}$ is unitarily equivalent to
the initial colligation ${\Col}$ defined in (\ref{3.1}). By Remark
\ref{R:3.1}, $\widetilde{\Col}$ has the same characteristic function as
${\Col}$, that is, $\T(z)$. It remains to check that the connecting
operator of $\widetilde{\Col}$ is an extension of ${\bf V}$, that is
\begin{equation}
\left[\begin{array}{cc}\widetilde{A} & \widetilde{B} \\
\widetilde{C}&\widetilde{D}\end{array}\right]
\left[\begin{array}{c} \operatorname{Col}_{s\in S}[M_sf]_{[s]}\\
Yf\end{array}\right]=\begin{bmatrix}
\operatorname{Col}_{r\in R}[N_rf]_{[r]} \\ Xf\end{bmatrix} \;
\text{ for every } f \in \cK_L\oplus\cK_R.
\label{3.26}
\end{equation}
To this end, note that by (\ref{3.24}), (\ref{3.25}) and block
partitionings (\ref{1.30}) and (\ref{3.21}) of $M_s$ and ${\mathbb T}$,
it holds that
$$
\tU_{[s]}^*[M_sf]_{[s]}={\mathbb T}_{[s]}(M_sf)=\begin{bmatrix}
{\mathbb T}_{[s,],L}E_{L,s} & {\mathbb
T}_{[s],R}\widetilde{M}_{s}(\del_R)\end{bmatrix}f
$$
for every $f\in\cK_L\oplus\cK_R$ and for every $s\in S$. Therefore,
\begin{equation}
\operatorname{Col}_{s\in S}\tU_{[s]}^*[M_sf]_{[s]}=
\operatorname{Col}_{s\in S} \begin{bmatrix}
{\mathbb T}_{[s],L}E_{L,s} &
{\mathbb T}_{[s],R}\widetilde{M}_{s}(\del_R)\end{bmatrix}f
\label{3.26a}
\end{equation}
which, on account of (\ref{2.10a}) and (\ref{2.11}) can be written as
\begin{equation}
\operatorname{Col}_{s\in S}\tU_{[s]}^*[M_sf]_{[s]}=
\begin{bmatrix}\left[\left(X_LH \right)^{\wedge L}(\del_L)\right]^* &
\left(Z_{\G,\cH} GY_R\right)^{\wedge R}(\del_R)\end{bmatrix}f.
\label{3.27}
\end{equation}
Similarly,  by (\ref{3.24}), (\ref{3.25}) and block
partitionings (\ref{1.31}) and (\ref{3.21}) of $N_r$ and ${\mathbb T}$,
it holds that
$$
[N_rf]_{[r]}=\tU_{[r]}{\mathbb T}_{[r]}(N_rf)=
\tU_{[r]}\begin{bmatrix}{\mathbb T}_{[r],L}\widetilde{N}_{r}(\del_L) &
{\mathbb T}_{[r],R}E_{R,r}\end{bmatrix}f\qquad (r\in R).
$$
Therefore,
\begin{equation}
\operatorname{Col}_{r\in R}\tU_{[r]}^*[N_rf]_{[r]}=
\operatorname{Col}_{r\in R}
\begin{bmatrix}{\mathbb T}_{[r],L}\widetilde{N}_{r}(\del_L)
& {\mathbb T}_{[r],R}E_{R,r}\end{bmatrix}f
\label{3.27a}
\end{equation}
which, on account of (\ref{2.10}) and (\ref{2.11a}) can be written as
\begin{equation}
\operatorname{Col}_{r\in R}\tU_{[r]}^*[N_rf]_{[r]}=
\begin{bmatrix} \left[\left(X_LH Z_{\G,\cH}\right)^{\wedge L}
(\del_L)\right]^*& \left(GY_R\right)^{\wedge R}(\del_R)\end{bmatrix}f.
\label{3.28}
\end{equation}
Thus, by (\ref{3.23}) and in view of \eqref{1.33}, (\ref{3.27}) and
(\ref{3.28}),
{\allowdisplaybreaks
\begin{eqnarray}
&&\left[\begin{array}{cc}\widetilde{A} & \widetilde{B} \\
\widetilde{C}&\widetilde{D}\end{array}\right]
\left[\begin{array}{c} \operatorname{Col}_{s\in S}[M_sf]_{[s]}
\\ Yf\end{array}\right]\label{3.29}\\
&&=\left[\begin{array}{cc}\oplus_{r\in R}\tU_{[r]} & 0 \\ 0 &
I\end{array}\right]\left[\begin{array}{cc}A & B \\ C& D\end{array}\right]
\left[\begin{array}{c} \operatorname{Col}_{s\in S}\tU_{[s]}^*[M_sf]_{[s]}
\\ Yf\end{array}\right]\nonumber\\
&&=\left[\begin{array}{cc}\oplus_{r\in R}\tU_{[r]} & 0 \\ 0 &
I\end{array}\right]\left[\begin{array}{cc}A & B \\ C& D\end{array}\right]
\begin{bmatrix}\left[\left(X_LH \right)^{\wedge L}(\del_L)\right]^* &
\left(Z_{\G,\cH} GY_R\right)^{\wedge R}(\del_R)\\ Y_L^*
& Y_R\end{bmatrix}f\nonumber\\
&&=\left[\begin{array}{cc}\oplus_{r\in R}\tU_{[r]} & 0 \\ 0 &
I\end{array}\right]
\begin{bmatrix}
\left[\left(X_LHZ_{\G,\cH} \right)^{\wedge L}(\del_L)\right]^* & \left(
GY_R\right)^{\wedge R}(\del_R) \\ X_L^* & X_R\end{bmatrix}=
\begin{bmatrix}\operatorname{Col}_{r\in R}[N_rf]_{[r]}\\
Xf\end{bmatrix},\nonumber
\end{eqnarray}
which proves (\ref{3.26}) and completes the proof of the lemma.\qed
\begin{La}
Let $\widetilde{\bU}$ of the form \eqref{3.12} be a unitary
extension of the isometry ${\bf V}$ given in \eqref{3.9}. Then the
characteristic function $\T$ of the unitary colligation of the form
\eqref{3.11a},
$$
\T(z)=\widetilde{D}+\widetilde{C}\left(
I-Z_{\G,\tH}(z)\widetilde{A}\right)^{-1}
Z_{\G,\tH}(z)\widetilde{B},
$$
is a solution to Problem \ref{P:1.10}.
\label{L:3.4}
\end{La}
{\bf Proof:} We use the arguments from the proof of the previous
lemma in the reverse order. We start with positive semidefinite operators
$\K_p\in\cL(\cH^\circ_p)$ (the spaces $\cH^\circ_p$ are given in
\eqref{3.6a}) and fix their factorizations
\begin{equation}
\K_p={\mathbb T}_p^*{\mathbb T}_p\quad\mbox{with}\quad
{\mathbb T}_p=\begin{bmatrix}{\mathbb T}_{p,L} & {\mathbb T}_{p,R}
\end{bmatrix}\colon\cH^\circ_p \to \cH_p
\label{1}
\end{equation}
where $\cH=\{\cH_p\colon \, p\in P\}$ is a collection of auxiliary
Hilbert spaces. Comparing \eqref{1} with \eqref{1.34} we get
factorizations
$$
\K_{p,L}={\mathbb T}_{p,L}^*{\mathbb T}_{p,L},\quad \K_{p,R}=
{\mathbb T}_{p,R}^*{\mathbb T}_{p,R}, \quad \K_{p,LR}
={\mathbb T}_{p,L}^*{\mathbb T}_{p,R}.
$$
for the block entries in $\K_p$ and more detailed decompositions
\eqref{1.35} lead us to equalities
\begin{eqnarray}
E_{L,s}^*{\mathbb T}_{[s],L}^*{\mathbb T}_{[s'],L}E_{L,s'}&=&
E_{L,s}^*\K_{[s],L}E_{L,s'}=\left[\K_{[s],L}\right]_{s,s'}=
\Psi_{s,s'},\label{3.52}\\
E_{L,s}^*{\mathbb T}_{[s],L}^*{\mathbb T}_{[r],R}E_{L,r}&=&
E_{L,s}^*\K_{[s],LR}E_{R,r}=\left[\K_{[s],LR}\right]_{s,r}
=\Lambda_{s,r},\label{3.53}\\
E_{R,r}^*{\mathbb T}_{[r],R}^*{\mathbb T}_{[r'],R}E_{R,r'}
&=&E_{R,r}^*\K_{[r],R}E_{R,r'}=\left[\K_{[r],L}\right]_{r,r'}=
\Phi_{r,r'}\label{3.54}
\end{eqnarray}
(where $E_{L,s}$ and $E_{R,r}$ are given by \eqref{1.28}, \eqref{1.28a})
holding for every choice of $s,s'\in S$ and $r,r'\in R$ so that
$[s]=[s']=[r]=[r']$. The latter equalities suggest the introduction
     of the operators
\begin{eqnarray}
{\mathbb F}_{L}&=&\operatorname{Col}_{s\in S}
{\mathbb T}_{[s],L}E_{L,s}\colon \; \cK_L\to \bigoplus_{s\in S} \cH_{[s]},
\label{3.57}\\
{\mathbb F}_{R}&=&\operatorname{Col}_{r\in R}
{\mathbb T}_{[r],R}E_{R,r}\colon \; \cK_R\to \bigoplus_{r\in R} \cH_{[r]}.
\label{3.57a}
\end{eqnarray}
We note the following two formulas
\begin{eqnarray}
\operatorname{Col}_{r\in R}{\mathbb T}_{[r],L}\widetilde{N}_r(\del_L)
&=&\left[\left({\mathbb F}_L^*\cdot Z_{\G,\cH}\right)^{\wedge
L}(\del_L)\right]^*,\label{3.55}\\
\operatorname{Col}_{s\in S}{\mathbb T}_{[s],R}\widetilde{M}_s(\del_R)
&=&\left(Z_{\G,\cH} \cdot {\mathbb F}_R\right)^{\wedge R}(\del_R),
\label{3.56}
\end{eqnarray}
which are similar to formulas \eqref{2.10} and \eqref{2.11} and are
verified in much the same way.

\smallskip

Let $\tU=\{\tU_p\colon p\in P\}$ be the collection of
unitary maps indexed by the  set of path-connected
components $P$ of $\G$ and defined via formulas  \eqref{3.24},
\eqref{3.25}. Then relations (\ref{3.26a}) and (\ref{3.27a})  hold by
construction; in view of \eqref{3.57}--\eqref{3.56} these relations can
be written as
\begin{eqnarray}
\operatorname{Col}_{s\in S}\tU_{[s]}^*[M_sf]_{[s]}
&=&\begin{bmatrix}{\mathbb F}_L & \left(Z_{\G,\cH} \cdot
{\mathbb F}_R\right)^{\wedge R}(\del_R)\end{bmatrix}f,\label{3.26b}\\
\operatorname{Col}_{r\in R}\tU_{[r]}^*[N_rf]_{[r]}
&=&\begin{bmatrix}\left[\left({\mathbb F}_L^*\cdot
Z_{\G,\cH}\right)^{\wedge L}(\del_L)\right]^* & {\mathbb
F}_R\end{bmatrix}f.\label{3.27b}
\end{eqnarray}
Now we define the operator
$$
\bU=\begin{bmatrix} A & B \\ C & D \end{bmatrix} =
\begin{bmatrix} [A_{r,s}] & [B_{r}] \\ [C_{s}] & D \end{bmatrix}
\colon \begin{bmatrix} \oplus_{s \in S} \cH_{[s]} \\ \cU
\end{bmatrix} \to
\begin{bmatrix} \oplus_{r \in R} \cH_{[r]} \\ \cY \end{bmatrix}
$$
in accordance to \eqref{3.24a} by
$$
\bU=\left[\begin{array}{cc}A & B \\ C& D\end{array}\right]=
\left[\begin{array}{cc}\oplus_{r\in R}\tU_{[r]}^* & 0 \\ 0 &
I\end{array}\right]
\left[\begin{array}{cc}\widetilde{A} & \widetilde{B} \\
\widetilde{C}&\widetilde{D}\end{array}\right]
\left[\begin{array}{cc}\oplus_{s\in S}\tU_{[s]} & 0 \\ 0 &
I\end{array}\right].
$$
By the assumption of the lemma, $\widetilde{\bU}$ extends ${\bf V}$:
$$
\left[\begin{array}{cc}\widetilde{A} & \widetilde{B} \\
\widetilde{C}&\widetilde{D}\end{array}\right]
\left[\begin{array}{c} \operatorname{Col}_{s\in S}[M_sf]_{[s]}
\\ Yf\end{array}\right]=\begin{bmatrix}
\operatorname{Col}_{r\in R}[N_rf]_{[r]}\\ Xf\end{bmatrix} \quad\mbox{for
every} \; \; f\in\cK_L\oplus\cK_R,
$$
which can be written in terms of $\bU$ as
$$
\left[\begin{array}{cc} A & B \\ C & D\end{array}\right]
\begin{bmatrix}\operatorname{Col}_{s\in S}\tU_{[s]}^*[M_sf]_{[s]}\\
     Yf\end{bmatrix}
=\begin{bmatrix}\operatorname{Col}_{r\in R}\tU_{[r]}^*[N_rf]_{[r]}\\
Xf\end{bmatrix}\quad (f\in\cK_L\oplus\cK_R).
$$
Upon substituting equalities \eqref{3.26b} and \eqref{3.27b}
and block decompositions \eqref{1.33} for $X$ and $Y$ in the latter
equality we get
\begin{equation}
\left[\begin{array}{cc} A & B \\ C & D\end{array}\right]
\begin{bmatrix}{\mathbb F}_L &
\left(Z_{\G,\cH} \cdot {\mathbb F}_R\right)^{\wedge R}(\del_R) \\
Y_L^* & Y_R\end{bmatrix}=\begin{bmatrix}\left[\left({\mathbb F}_L^*
\cdot Z_{\G,\cH}\right)^{\wedge L}(\del_L)\right]^* & {\mathbb F}_R\\
X_L^* & X_R\end{bmatrix}.
\label{3.50}
\end{equation}
By Remark \ref{R:3.1}, the colligations $\Col$ and $\widetilde{\Col}$
defined in \eqref{3.1} and \eqref{3.11a} have the same
characteristic functions and thus, $\T$ can be taken in the form
\eqref{1.18}. Let $H(z)$ and $G(z)$ be defined as in (\ref{1.24}) and
decomposed as in (\ref{1.19}) and (\ref{1.21}). We shall use the
representations  (\ref{3.13}) of $\T(z)$ which  are equivalent  to
(\ref{1.18}).

\smallskip

Since $\bU$ is unitary, it follows from (\ref{3.50}) that
\begin{eqnarray}
A^*\left[\left({\mathbb F}_{L}^*Z_{\G,\cH}\right)^{\wedge
L}(\del_L)\right]^*
+C^*X_L^*&=&{\mathbb F}_{L},\label{3.58}\\
B^*\left[\left({\mathbb F}_{L}^*Z_{\G,\cH}\right)^{\wedge
L}(\del_L)\right]^*+D^*X_L^*&=&Y_L^*,\label{3.59}\\
A\left(Z_{\G,\cH} {\mathbb F}_{R}\right)^{\wedge R}(\del_R)
+BY_R&=&{\mathbb F}_{R},\label{3.60}\\
C\left(Z_{\G,\cH} {\mathbb F}_{R}\right)^{\wedge R}(\del_R)
+DY_R&=&X_R.\label{3.61}
\end{eqnarray}
Taking adjoints in \eqref{3.58} we get
$$
X_L C={\mathbb F}_{L}^*-\left(
{\mathbb F}_{L}^*Z_{\G,\cH}\right)^{\wedge L}(\del_L)A
$$
which can be written, by properties \eqref{2.1} and \eqref{2.2} of the
left evaluation map, as
$$
X_LC=\left({\mathbb F}_{L}^*(I-Z_{\G,\cH} A)\right)^{\wedge L}(\del_L).
$$
Multiplying both sides in the last equality by $(I-Z_{\G,\cH}(z)A)^{-1}$
on
the
right and applying the left evaluation map to the resulting identity
$$
X_LH(z)=\left({\mathbb F}_{L}^*(I-Z_{\G,\cH} A)\right)^{\wedge
L}(\del_L)\cdot
(I-Z_{\G,\cH}(z)A)^{-1},
$$
we get
\begin{eqnarray}
\left(X_LH\right)^{\wedge L}(\del_L)&=&\left(\left({\mathbb F}_{L}^*
(I-Z_{\G,\cH} A)\right)^{\wedge L}(\del_L)(I-Z_{\G,\cH}
A)^{-1}\right)^{\wedge L}(\del_L) \nonumber\\
&=&\left({\mathbb F}_{L}^*(I-Z_{\G,\cH} A)(I-Z_{\G,\cH}
A)^{-1}\right)^{\wedge
L}(\del_L)\nonumber\\
&=&\left({\mathbb F}_{L}^*\right)^{\wedge L}(\del_L)={\mathbb
F}_{L}^*.\label{3.62}
\end{eqnarray}
Note that the second equality in the last chain has been obtained upon
applying \eqref{2.5} to
$$
T(z)={\mathbb F}_{L}^*(I-Z_{\G,\cH}(z)A)\quad\mbox{and}\quad
\widetilde{T}(z)=(I-Z_{\G,\cH}(z)A)^{-1},
$$
whereas the third equality follows by the property \eqref{2.1}.

Next  we take adjoints in
\eqref{3.59} to get
\begin{equation}
Y_L=\left({\mathbb F}_{L}^*Z_{\G,\cH}\right)^{\wedge L}(\del_L)B+X_LD=
\left({\mathbb F}_{L}^*Z_{\G,\cH} B\right)^{\wedge L}(\del_L)+X_LD.
\label{mn}
\end{equation}
By \eqref{3.62},
$$
\left({\mathbb F}_{L}^*Z_{\G,\cH} B\right)^{\wedge L}(\del_L)=
\left(\left(X_LH\right)^{\wedge L}(\del_L)\cdot Z_{\G,\cH}
B\right)^{\wedge L}(\del_L)
$$
and applying \eqref{2.5} to
$$
T(z)=X_LH(z)\quad\mbox{and}\quad
\widetilde{T}(z)=Z_{\G,\cH}(z)B
$$
leads us to
$$
\left({\mathbb F}_{L}^*Z_{\G,\cH} B\right)^{\wedge
L}(\del_L)=\left((X_LHZ_{\G,\cH} B\right)^{\wedge L}(\del_L).
$$
Substituting the latter equality into the left hand side expression in
\eqref{mn} and making use of the first representation of $S$ in
(\ref{3.13}), we get
\begin{eqnarray*}
Y_L&=& \left(X_LHZ_{\G,\cH} B\right)^{\wedge L}(\del_L)+X_LD\\
&=&\left(X_LHZ_{\G,\cH} B+X_LD\right)^{\wedge
L}(\del_L)=\left(X_LS\right)^{\wedge L}(\del_L),
\end{eqnarray*}
which proves  the first interpolation condition in (\ref{1.27}).

\medskip

To get the second interpolation condition in (\ref{1.27}) write
\eqref{3.60} in the form
$$
BY_R=\left(I-AZ_{\G,\cH}){\mathbb F}_{R}\right)^{\wedge R}(\del_R),
$$
multiply the latter equality by $(I-AZ_{\G,\cH}(z))^{-1}$ on the left and
apply
the right evaluation map to the resulting identity
$$
G(z)Y_R=(I-AZ_{\G,\cH}(z))^{-1}\left(I-AZ_{\G,\cH}){\mathbb
F}_{R}\right)^{\wedge R}(\del_R).
$$
We have
\begin{eqnarray}
\left(GY_R\right)^{\wedge R}(\del_R)
&=&\left((I-AZ_{\G,\cH})^{-1}\left((I-AZ_{\G,\cH}){\mathbb
F}_{R}\right)^{\wedge
R}(\del_R)\right)^{\wedge R}(\del_R)\nonumber\\
&=&\left((I-AZ_{\G,\cH})^{-1}(I-AZ_{\G,\cH}){\mathbb F}_{R}\right)^{\wedge
R}(\del_R)\nonumber\\
&=&\left({\mathbb F}_{R}\right)^{\wedge R}(\del_R)={\mathbb
F}_{R}.\label{3.63}
\end{eqnarray}
Note that the third equality in the last chain has been obtained upon
applying \eqref{2.6} to
$$
T(z)=(I-AZ_{\G,\cH}(z))^{-1}\quad\mbox{and}\quad
\widetilde{T}(z)=(I-AZ_{\G,\cH}(z)){\mathbb F}_{R}.
$$
Substituting \eqref{3.63} into \eqref{3.61} and applying
\eqref{2.6} to
$$
T(z)=CZ_{\G,\cH}(z)\quad\mbox{and}\quad
\widetilde{T}(z)=G(z)Y_R,
$$
we get
\begin{eqnarray*}
X_R&=&\left(CZ_{\G,\cH}\left(GY_R\right)^{\wedge
R}(\del_R)\right)^{\wedge R}(\del_R)+DY_R\\
&=& \left(CZ_{\G,\cH} GY_R\right)^{\wedge R}(\del_R)+DY_R\\
&=& \left(CZ_{\G,\cH} GY_R+DY_R\right)^{\wedge R}(\del_R)
\end{eqnarray*}
which coincides with the second equality in (\ref{1.27}), due to
the second representation in (\ref{3.13}).

\smallskip

Thus, $\T$ belongs to ${\mathcal {SA}}_\G(\cU, \, \cY)$ as the
characteristic function of a unitary colligation \eqref{3.1} and
satisfies interpolation conditions (\ref{1.27}).
It remains to show that it satisfies also conditions
(\ref{1.37})--(\ref{1.39}). But it follows from
\eqref{3.62}, \eqref{3.63} and \eqref{3.57} that
$$
\left(X_LH_s\right)^{\wedge
L}(\del_L)=E_{L,s}^*{\mathbb T}_{[s],L}^*\quad\mbox{and}\quad
\left(G_r Y_R\right)^{\wedge L}(\del_R)={\mathbb T}_{[r],R}E_{R,r}
$$
for $s\in S$ and  $r\in R$. Now we pick any $s,s'\in S$ and $r,r'\in R$ so
that $[s]=[s']=[r]=[r']$ and combine the two latter equalities with
\eqref{3.52}--\eqref{3.54} to get \eqref{1.37}--\eqref{1.39}:
\begin{eqnarray*}
(X_LH_s)^{\wedge L}(\del_L)\left[(X_LH_{s'})^{\wedge L}(\del_L)\right]^*
&=&E_{L,s}^*{\mathbb T}_{[s],L}^*{\mathbb
T}_{[s'],L}E_{L,s'}=\Psi_{s,s'},\\
(X_LH_s)^{\wedge L}(\del_L)\left(G_r Y_R\right)^{\wedge R}(\del_R)
&=&E_{L,s}^*{\mathbb T}_{[s],L}^*{\mathbb
T}_{[r],R}E_{R,r}=\Lambda_{s,r},\\
\left[(G_rY_R)^{\wedge R}(\del_R)\right]^*(G_{r'} Y_R)^{\wedge R}(\del_R)
&=&E_{R,r}^*{\mathbb T}_{[r],R}^*{\mathbb T}_{[r'],R}E_{R,r'}=\Phi_{r,r'},
\end{eqnarray*}
and complete the proof.\qed

\medskip

\section{The universal unitary colligation associated with the
interpolation problem}  \label{S:universal}
\setcounter{equation}{0}

A general result of Arov and Grossman (see \cite{arovgros1},
\cite{arovgros2}) describes how to parametrize the set of all unitary
extensions of a given partially defined isometry ${\bf V}$.
Their result has been extended to the multivariable (commutative) case in
\cite{bt, BTV, bb5} and will be extended in this section to
the setting of noncommutative power series.

\smallskip

Let ${\bf V} : \; {\mathcal D}_{\bf V}\rightarrow {\mathcal R}_{\bf
V}$ be the isometry given in (\ref{3.9}) with ${\mathcal D}_{\bf V}$
and ${\mathcal R}_{\bf V}$ given in (\ref{3.10}) and (\ref{3.11}).
Introduce the defect spaces
$$
\Delta=\left[\begin{array}{c}\oplus_{s\in S}\hH_{[s]}\\
\cU\end{array}\right]
\ominus{\mathcal {D}}_{\bf V} \quad{\rm and}\quad
\Delta_*=\left[\begin{array}{c}\oplus_{r\in R}\hH_{[r]}\\
\cY\end{array}\right]\ominus{\mathcal {R}}_{\bf V}
$$
and let $\widetilde{\Delta}$ to be another copy of $\Delta$ and
$\widetilde{\Delta}_*$ to be another copy of $\Delta_*$ with unitary
identification maps
$$
i: \; \Delta\rightarrow \widetilde{\Delta}\quad\mbox{and}\quad i_*:
\; \Delta_*\rightarrow \widetilde{\Delta}_*.
$$
Define a unitary operator ${\bU}_0$ from ${\mathcal {D}}_{\bf
V}\oplus\Delta\oplus\widetilde{\Delta}_*$ onto
${\mathcal {R}}_{\bf V}\oplus\Delta_*\oplus\widetilde{\Delta}$
by the rule
\begin{equation}
{\bU}_0x=\left\{ \begin{array}{ll} {\bf V}x, & \mbox{if} \; \;
x\in{\mathcal {D}}_{\bf V},\\ i(x) & \mbox{if} \; \; x\in\Delta,\\
i_*^{-1}(x)& \mbox{if} \; \; x\in \widetilde{\Delta}_*.
\end{array}\right.
\label{4.2}
\end{equation}
Identifying $\begin{bmatrix}{\mathcal {D}}_{\bf V} \\ \Delta\end{bmatrix}$
with $\begin{bmatrix}\oplus_{s\in S}\hH_{[s]}\\ \cU\end{bmatrix}$
and $\begin{bmatrix}{\mathcal {R}}_{\bf V} \\ \Delta_*\end{bmatrix}$ with
$\begin{bmatrix}\oplus_{r\in R}\hH_{[r]} \\ \cY\end{bmatrix}$,
we decompose ${\bU}_0$ defined by (\ref{4.2}) according to
\begin{equation}
{\bf U}_0=\left[\begin{array}{ccc}U_{11} &U_{12}&U_{13} \\
U_{21}&U_{22}& U_{23}\\ U_{31}&U_{32}&0\end{array}\right]:
\quad\begin{bmatrix}\oplus_{s\in S}\hH_{[s]} \\ \cU \\
\widetilde{\Delta}_*\end{bmatrix}\rightarrow
\begin{bmatrix}\oplus_{r\in R}\hH_{[r]} \\
\cY \\ \widetilde{\Delta}\end{bmatrix}.
\label{4.3}
\end{equation}
The $(3,3)$ block in this decomposition is zero, since (by definition
(\ref{4.2})), for every $x\in\widetilde{\Delta}_*$, the vector
${\bU}_0x$ belongs to $\Delta$, which is a subspace of
$\left[\begin{array}{c}\oplus_{r\in R}\hH_{[r]} \\ \cY\end{array}\right]$
and therefore, is orthogonal to $\widetilde{\Delta}$ (in other words
${\bf P}_{\widetilde{\Delta}}\bU_0\vert_{\widetilde{\Delta}_*}=0$,
where ${\bf P}_{\widetilde{\Delta}}$ stands for the orthogonal
projection of ${\mathcal {R}}_{\bf V}\oplus\Delta_*\oplus\widetilde{\Delta}$
onto $\widetilde{\Delta}$).

\medskip

The unitary operator ${\bf U}_0$ is the
connecting operator of the unitary colligation
\begin{equation}
\Col_0=\left\{\G, \; \hH, \; \left[\begin{array}{c} \cU \\
\widetilde{\Delta}_*\end{array}\right], \; \left[\begin{array}{c} \cY
\\ \widetilde{\Delta}\end{array}\right], \; {\bU}_0\right\},
\label{4.4}
\end{equation}
which is called {\em the universal unitary colligation} associated
with the interpolation Problem \ref{P:1.10}.

\medskip

Let $\widetilde{\Col}$ be any colligation of the form
\begin{equation}
\widetilde{\Col}=\left\{\G, \; \tH, \; \widetilde{\Delta}, \;
\widetilde{\Delta}_*, \; \widetilde{\bU}\right\}.
\label{4.5}
\end{equation}
We define another colligation ${\cF}_{\Col_0}[\widetilde{\Col}]$,
called the {\em coupling} of $\Col_0$ and $\widetilde{\Col}$, to be
the colligation of the form
$$
{\cF}_{\Col_0}[\widetilde{\Col}]=\left\{\G, \; \hH\oplus\tH, \;
\cU, \; \cY, \; \cF_{\bU_0}[\widetilde{\bU}]\right\}
$$
with the connecting operator $\cF_{\bU_0}[\widetilde{\bU}]$
defined as follows:
\begin{equation}
\cF_{\bU_0}[\widetilde{\bU}]: \; \begin{bmatrix}c \\ h \\ u\end{bmatrix}
\to \begin{bmatrix}c^\prime \\ h^\prime \\ y\end{bmatrix}
\label{4.6}
\end{equation}
if the system of equations
\begin{equation}
\bU_0: \; \begin{bmatrix}c \\ u \\ \td_*\end{bmatrix}
\to \begin{bmatrix}c^\prime \\ y \\ \td\end{bmatrix}\quad\mbox{and}\quad
\widetilde{\bU}: \; \begin{bmatrix}h \\ \td\end{bmatrix}
\to \begin{bmatrix}h^\prime \\ \td_*\end{bmatrix}
\label{4.7}
\end{equation}
is satisfied for some choice of $\td\in\widetilde{\Delta}$ and
$\td_*\in\widetilde{\Delta}_*$. To show that the operator
$\cF_{\bU_0}[\widetilde{\bU}]$ is well defined, i.e., that for every
triple $(c, \, h, \, u)$, there exist $\td$ and $\td_*$ for which
the system (\ref{4.7}) is consistent and the resulting triple
$(c^\prime, \, h^\prime, \, y)$ does not depend on the choice of
$\td$ and $\td_*$, we note first that, on account of (\ref{4.2}) and
(\ref{4.3}), the the bottom component of the first equation in  (\ref{4.7})
determines $\td$ uniquely by
$$
\td={\bf P}_{\widetilde{\Delta}}\left({\bf V}{\bf P}_{{\mathcal D}_{\bf
V}}+i{\bf P}_\Delta\right)\begin{bmatrix} c\\ u\end{bmatrix}=
i{\bf P}_\Delta\begin{bmatrix} c\\ u\end{bmatrix}.
$$
With this $\td$, the the bottom component of the second equation
in (\ref{4.7}) determines uniquely
$\td_*$ and $h^\prime$. Using $\td_*$ one can recover now $c^\prime$ and
$y$ from the first and second components of the first equation in (\ref{4.7}).

\smallskip

Since operators $\bU_0$ and $\widetilde{\bU}$ are unitary, it follows
from (\ref{4.7}) that
\begin{eqnarray*}
\|c\|^2+\|u\|^2+\|\td_*\|^2&=&\|c^\prime\|^2+\|y\|^2+\|\td\|^2,\\
\|h\|^2+\|\td\|^2&=&\|h^\prime\|^2+\|\td_*\|^2
\end{eqnarray*}
and therefore, that
$$
\|c\|^2+\|u\|^2+\|h\|^2=\|c^\prime\|^2+\|y\|^2+\|h^\prime\|^2,
$$
which means that the coupling operator $\cF_{\bU_0}[\widetilde{\bU}]$
is isometric. A similar argument can be made with the adjoints of
${\mathbf U}_{0}$, $\widetilde{\mathbf U}$ and ${\mathcal
F}_{{\mathbf U}_{0}}[\widetilde{\mathbf U}]$, and hence ${\mathcal
F}_{{\mathbf U}_{0}}[\widetilde{\mathbf U}]$ is unitary.
Furthermore, by (\ref{4.6}) and (\ref{4.7}),
$$
\cF_{\bU_0}[\widetilde{\bU}]\vert_{(\oplus_{s\in S}\hH_{[s]})\oplus\cU}
={\bU_0}\vert_{(\oplus_{s\in S}\hH_{[s]})\oplus\cU}
$$
and since ${\mathcal D}_{\bf V}\subset (\oplus_{s\in
S}\hH_{[s]})\oplus\cU$, it follows that
\begin{equation}
\cF_{\bU_0}[\widetilde{\bU}]\vert_{{\mathcal D}_{\bf V}}=
{\bU_0}\vert_{{\mathcal D}_{\bf V}}={\bf V}.
\label{4.8}
\end{equation}
Thus, the coupling of the connecting operator $\bU_0$ of the universal
unitary colligation associated with Problem \ref{P:1.10} and any other
unitary operator is a unitary extension of the isometry ${\bf V}$ defined
in (\ref{3.9}). Conversely for every unitary colligation
$\Col=\{\G, \; \hH\oplus{\tH}, \; \cU, \; \cY, \; \bU\}$ with the
connecting operator being a unitary extension of ${\bf V}$, there
exists a unitary colligation $\widetilde{\Col}$ of the form
(\ref{4.5}) such that $\Col={\cF}_{\Col_0}[\widetilde{\Col}]$
(the proof is the same as in \cite[Theorem 6.2]{bt}). Thus,
all unitary extensions $\bU$ of the isometry ${\bf V}$ defined
in (\ref{3.9}) are parametrized by the formula
\begin{equation}
\bU=\cF_{\bU_0}[\widetilde{\bU}],\qquad\widetilde{\bU}: \;
(\oplus_{s\in S}\hH_{[s]})\oplus\widetilde{\Delta}\to
(\oplus_{r\in R}\hH_{[r]})\oplus\widetilde{\Delta}_*
\label{4.9}
\end{equation}
and $\tH = \{ \tH \colon p \in P\}$ is a collection of auxiliary
Hilbert spaces indexed by the path-connected components $p \in P =
P(\G)$ of the admissible graph $\G$.

\medskip

According to \eqref{3.3}, the characteristic function of the
colligation $\Col_0$ defined in (\ref{4.4}) with the connecting
operator ${\bf U}_0$ partitioned as in \eqref{4.3}, is given by
\begin{eqnarray}
\Sigma(z)&=&\left[\begin{array}{cc}\Sigma_{11}(z) & \Sigma_{12}(z) \\
\Sigma_{21}(z) & \Sigma_{22}(z)\end{array}\right]\nonumber\\
&=&\left[\begin{array}{cc}U_{22} &U_{23}\\ U_{32}&0\end{array}\right]
+\left[\begin{array}{c}U_{21}\\U_{31}\end{array}\right] \left(I-
Z_{\G,\hH}(z)U_{11}\right)^{-1}Z_{\G,\hH}(z)\begin{bmatrix}U_{12}&
U_{13}\end{bmatrix}
\label{4.10}
\end{eqnarray}
and belongs to the class
$\mathcal{ S A}_{{\G}}(\cU\oplus\widetilde{\Delta}_*,
\; \cY\oplus \widetilde{\Delta})$ by Theorem \ref{T:conservative-real}.

\begin{Tm}
Let ${\bf V}$ be the isometry defined in \eqref{3.9}, let $\Sigma$
be constructed as above and let $\T$ be an element in
$\cL(\cU,\cY)\langle\langle z\rangle\rangle$. Then the following are
equivalent:
\begin{enumerate}
\item $\T$ is a solution of Problem \ref{P:1.10}.
\item $\T$ is a characteristic function of a colligation
$\Col=\{\G, \; \hH\oplus{\tH}, \; \cU, \; \cY, \; \bU\}$ with the
connecting operator $\bU$ being a unitary extension of ${\bf V}$.
\item $\T$ is of the form
\begin{equation}
\T(z)=\Sigma_{11}(z)+\Sigma_{12}(z)\left(I_{\widetilde{\Delta}_*}-
{\CE}(z)\Sigma_{22}(z)\right)^{-1}{\CE}(z)\Sigma_{21}(z)
\label{4.11}
\end{equation}
where $\CE$ is a power series from the noncommutative Schur-Agler
class
\newline
$\mathcal{ S A}_{\G}(\widetilde{\Delta}, \;
\widetilde{\Delta}_*)$.
\end{enumerate}
\label{T:4.1}
\end{Tm}
{\bf Proof:} The equivalence ${\bf 1\Longleftrightarrow 2}$ follows by
Lemmas \ref{L:3.3} and \ref{L:3.4}.

\smallskip

${\bf 2\Longrightarrow 3}$. By the preceding analysis, the colligation
$\Col$ is the coupling of the universal colligation $\Col_0$ defined in
(\ref{4.4}) and some unitary colligation $\widetilde{\Col}$ of the
form (\ref{4.5}). The connecting operators $\bU$, $\bU_0$ and
$\widetilde{\bU}$ of these colligations are related as in (\ref{4.9}).
Let $\T$, $\Sigma$ and $\CE$ be characteristic functions of
$\Col$, $\Col_0$ and $\widetilde{\Col}$, respectively.
Applying Remark \ref{R:3.2} to (\ref{4.6}) and (\ref{4.7}), we get
\begin{equation}
\T(z)e=e_*,\quad \Sigma(z)\begin{bmatrix}u\\ \td_*\end{bmatrix}=
\begin{bmatrix}y\\ \td\end{bmatrix},\quad \CE(z)\td=\td_*.
\label{4.12}
\end{equation}
Substituting the third relation in (\ref{4.12}) into the second we get
$$
\Sigma(z)\begin{bmatrix}u\\ \CE(z)\td\end{bmatrix}=\begin{bmatrix}y\\
\td\end{bmatrix},
$$
which in view of the block decomposition (\ref{4.10}) of $\Sigma$ splits
into
$$
\Sigma_{11}(z)u+\Sigma_{12}(z)\CE(z)\td=y\quad\mbox{and}\quad
\Sigma_{21}(z)u+\Sigma_{22}(z)\CE(z)\td=\td.
$$
The second from the two last equalities gives
$$
\td=\left(I-\Sigma_{22}(z)\CE(z)\right)^{-1}\Sigma_{21}(z)u
$$
which, being substituted into the first equality, implies
$$
\left(\Sigma_{11}(z)+\Sigma_{12}(z)\CE(z)
\left(I-\Sigma_{22}(z)\CE(z)\right)^{-1}\Sigma_{21}(z)\right)u=y.
$$
The latter is equivalent to
$$
\left(\Sigma_{11}(z)+\Sigma_{12}(z)
\left(I-\CE(z)\Sigma_{22}(z)\right)^{-1}\CE(z)\Sigma_{21}(z)\right)u=y
$$
and the comparison of the last equality with the first relation in
(\ref{4.12}) leads to representation (\ref{4.11}) of $\T$, since a
vector $u\in\cU$ is arbitrary.

\smallskip

${\bf 3\Longrightarrow 2}$. Let $\T$ be of the form (\ref{4.11}) for some
$\CE\in \mathcal{ S A}_{\G}
(\widetilde{\Delta}, \;  \widetilde{\Delta}_*)$.
By Theorem \ref{T:conservative-real}, $\CE$ is the characteristic function of
a unitary colligation $\widetilde{\Col}$ of the form
(\ref{4.5}). Let $\Col$ be the unitary colligation defined by
$\Col={\cF}_{\Col_0}[\widetilde{\Col}]$. By the preceding
``${\bf 2 \Longrightarrow 3}$'' part,  $\T$ of the form (\ref{4.11}) is
the characteristic function of $\Col$. It remains to note that the
colligation $\Col$ is of required the form: its input and output spaces
coincide with $\cU$ and $\cY$, respectively (by the definition
of coupling) and its connecting operator is an extension of ${\bf V}$,
by \eqref{4.8}.\qed

\medskip

As a corollary we obtain the sufficiency part in both Theorem \ref{T:1.9}
and Theorem \ref{T:1.11},
including the parametrization of the set of all solutions of Problem
\ref{P:1.10} in Theorem \ref{T:1.11} and the parametrization of the
set of all solutions of Problem \ref{P:1.4} in Corollary
\ref{C:parametrization}.

\section{Examples and special cases} \label{S:examples}
\setcounter{equation}{0}

For certain special cases of Problems \ref{P:1.4} and \ref{P:1.10},
the general interpolation results stated in Theorems \ref{T:1.9} and
\ref{T:1.11} become much more transparent. Moreover, some of
these particular cases are quite important for applications and
are interesting in their own right; it seems reasonable therefore
to display them in more detail.

\subsection{Left sided interpolation problems}

The left sided problem
can be considered as the special case of Problem \ref{P:1.4} when $T_R$
is a tuple of operators acting on the space of dimension zero.
\begin{Pb}
Given an admissible data set ${\mathcal D}=\{T_L, \, X_L, \, Y_L\}$,
find necessary
and sufficient conditions for existence of a power series
$\T\in\mathcal{SA}_{\G}(\cU, \cY)$ such that
\begin{equation}
\left(X_L\T\right)^{\wedge L}(\del_L)=Y_L.
\label{ex1}
\end{equation}
\label{P:ex1}
\end{Pb}
The answer follows immediately from Theorem \ref{T:1.9}.
\begin{Tm}
There is a power series $\T\in{\mathcal {SA}}_\G(\cU, \, \cY)$
satisfying interpolation condition \eqref{ex1} if and only if there
exists a collection $\K_{L}=\{\K_{p,L}\colon \, p\in P\}$ of positive
semidefinite operators on the space $\oplus_{s\in S\colon [s]=p}\cK_L$
indexed by the  set of path-connected components $P$ of $\G$,
which satisfies the Stein identity
\begin{equation}
\sum_{s\in S}E_{L,s}^*\K_{[s],L}E_{L,s}-\sum_{r\in R}
\widetilde{N}_{r}(\del_L)^*\K_{[r],L}
\widetilde{N}_{r}(\del_L)=X_LX_L^*-Y_LY_L^*,
\label{ex2}
\end{equation}
where $E_{L,s}$ and $\widetilde{N}_{r}$ are the operators defined via
formulas \eqref{1.28} and \eqref{1.29}, respectively.
\label{T:ex2}
\end{Tm}
Furthermore, it follows by Theorem \ref{T:1.11} that for every choice
of a tuple $\K_{L}$ satisfying the conditions of Theorem \ref{T:ex2},
there exists a power series $\T\in{\mathcal {SA}}_\G (\cU, \, \cY)$
satisfying (besides the left interpolation condition \eqref{ex1})
supplementary interpolation conditions
\begin{equation}
(X_LH_s)^{\wedge L}(\del_L)\left[(X_LH_{s'})^{\wedge L}(\del_L)\right]^*
= \Psi_{s,s'} \;  \text{ for } s, s'\in
S\colon \; [s]=[s'],\label{ex3}
\end{equation}
for some choice of associated
function $H(z)$ in representation \eqref{1.20} of $\T$. Furthermore,
all such $\T$ can be parametrized by a linear fractional transformation.
We leave to the reader to formulate the right sided interpolation problem
and to derive the right sided version of Theorem \ref{T:ex2} from
Theorem \ref{T:1.9}.

Parallel results hold for right-sided interpolation problems; we
leave the formulation of explicit statements to the reader.

\subsection{The case of the noncommutative ball}  \label{S:noncomball}

Now we consider the Fornasini-Marchesini case (see Example
\ref{E:NCFMimp} above) where $S=\{1\}$ and
$R=E=\{1,\ldots, d\}$. In this case, from Corollaries
\ref{C:strictrowcon} and \ref{C:strictrowcon'} we see that a sufficient
condition for
$T_L=(T_{L,1},\ldots,T_{L,d})$ to be left-admissible is that $T_{L}$
be a strict row contraction and that a sufficient condition for
$T_R=(T_{R,1},\ldots,T_{R,d})$ to be right admissible is that
$T_{R}$ be a strict column contraction:
$$
\sum_{j=1}^d T_{L,j}T_{L,j}^{*}< I_{\cK_L}\quad\mbox{and}\quad
\sum_{j=1}^d T_{R,j}^{*}T_{R,j}< I_{\cK_R}.
$$
The left sided problem is of special interest.
\begin{Pb}
Given an admissible data set ${\mathcal D}=\{T_L, \, X_L, \, Y_L\}$,
find necessary and sufficient conditions for existence of a power series
$\T\in\mathcal{SA}_{\G^{\rm FM}}(\cU, \cY)$ satisfying the left sided
interpolation condition \eqref{ex1}.
\label{P:ex3}
\end{Pb}
In this particular case
\begin{equation}
E_{L}=I_{\cK_L},\qquad \widetilde{N}_{j}(\del_L)=\del_{L,j}^*
\label{ex4}
\end{equation}
and we conclude by Theorem \ref{T:ex2} that there exists a power series
$\T\in{\mathcal {SA}}_{\G^{\rm FM}}(\cU, \, \cY)$ satisfying \eqref{ex1}
if and only if there exists a positive semidefinite operator
$\K_{L}$ subject to the Stein identity
$$
\K_{L}-\sum_{j=1}^d \del_{L,j}\K_{L}\del_{L,j}^*=X_LX_L^*-Y_LY_L^*.
$$
Since the $d$-tuple $T_L$ is a strict row contraction, the latter Stein
equation has a unique solution given in terms of convergent series by
\begin{equation}
\K_{L}=\sum_{v\in\cF_E}T_L^v\left(X_LX_L^*-Y_LY_L^*\right)(T_L^*)^v
\label{ex5}
\end{equation}
and we come to the following.

\begin{Tm}
      Assume that $\del_{L} = (\del_{L,1}, \dots, \del_{L,d})$ is a
      strict row contraction.  Then
there is a power series $\T\in{\mathcal {SA}}_{\G^{\rm FM}}(\cU, \, \cY)$
satisfying interpolation condition \eqref{ex1} if and only if
the operator $\K_{L}$ defined in \eqref{ex5} is positive semidefinite.
\label{T:ex5}
\end{Tm}

A remarkable part about the left-sided interpolation for the
Fornasini-Marchesini case is that no supplementary conditions are needed
to get a parametrization of the solution set: since the operator
$\K_{L}$ is uniquely determined by the interpolation date, it follows by
Theorem \ref{T:1.11} that for every $\T\in{\mathcal {SA}}_{\G^{\rm
FM}}(\cU, \, \cY)$ satisfying \eqref{ex1}, the function $H(z)$ associated
with $\T$ via representation \eqref{1.20}, satisfies
$$
(X_LH)^{\wedge L}(\del_L)\left[(X_LH)^{\wedge L}(\del_L)\right]^*=\K_L.
$$
Furthermore, in this case the power series $\Sigma$ defined in
\eqref{1.41} depends on the data $\{T_L, \, X_L, \, Y_L\}$ only and the
linear fractional formula \eqref{1.42} parametrizes the solution set to
Problem \ref{P:ex3}.

\medskip

The two sided problem in the Fornasini-Marchesini case is less remarkable.
\begin{Pb}
Given an admissible interpolation data set \eqref{1.26}, find necessary
and sufficient conditions for existence of a power series
$\T\in\mathcal{SA}_{\G^{\rm FM}}(\cU, \cY)$ such that
\begin{equation}
\left(X_L\T\right)^{\wedge L}(\del_L)=Y_L \quad\text{and}\quad
\left(\T Y_R\right)^{\wedge R}(\del_R)=X_R.
\label{11.1}
\end{equation}
\label{P:11.1}
\end{Pb}
The formulas \eqref{1.30} and \eqref{1.31} read
\begin{equation}
M=\begin{bmatrix} I_{\cK_L} & 0  \\ 0 & T_{R,1} \\
\vdots & \vdots \\ 0 & T_{R,d}\end{bmatrix}\quad \mbox{and}\quad
N_j=\begin{bmatrix}T_{L,j}^* & 0 \\ 0 & E_j\end{bmatrix}\quad
(j=1,\ldots,d)
\label{11.2}
\end{equation}
where
$$
E_1=\begin{bmatrix}I_{\cK_R} \\ 0 \\ \vdots \\ 0\end{bmatrix},\quad
E_2=\begin{bmatrix}0 \\ I_{\cK_R} \\ \vdots \\ 0\end{bmatrix}, \; \ldots
\; , \; E_d=\begin{bmatrix}0 \\ \vdots \\ 0 \\ I_{\cK_R}\end{bmatrix}.
$$
Now Theorem \ref{T:1.9} leads us to the following conclusion:
\begin{Tm}
There is a power series
$\T\in{\mathcal {SA}}_{\G^{\rm FM}}(\cU, \, \cY)$
satisfying interpolation conditions \eqref{11.1} if and only if
there exists a positive semidefinite operator
\begin{equation}
\K=\begin{bmatrix}\K_{L} & \K_{LR} \\ \K_{LR}^* &
\K_{R}\end{bmatrix}\in\cL(\cK_L\oplus\cK_R^d)
\label{11.3}
\end{equation}
subject to the Stein identity
\begin{equation}
M^*\K M-\sum_{j=1}^d N_j^*\K N_j=X^*X-Y^*Y,
\label{11.4}
\end{equation}
where $M$, $N_j$, $X$ and $Y$ are defined in \eqref{11.2} and
\eqref{1.33}.
\label{T:ex7}
\end{Tm}
Since the block $\K_L$ in \eqref{11.3} is uniquely determined from the
left interpolation data via the Stein identity \eqref{11.4}, the latter
result can be displayed more explicitly in terms of a structured positive
completion problem.
\begin{Tm}
There is a power series
$\T\in{\mathcal {SA}}_{\G^{\rm FM}}(\cU, \, \cY)$
satisfying interpolation conditions \eqref{11.1} if and only if
there exist operators $\Lambda_j\in\cL(\cK_R,\cK_L)$ and
$\Phi_{ij}\in\cL(\cK_R)$ for $i,j=1,\ldots,d$ subject to Stein identities
\begin{eqnarray}
\sum_{j=1}^d(T_{L,j}\Lambda_j-\Lambda_jT_{R,j})&=&Y_LY_R-X_LX_R,\label{11.5}\\
\sum_{j=1}^d\Phi_{jj}-\sum_{i,j=1}^d
T_{R,i}^*\Phi_{ij}T_{R,j}&=&Y_R^*Y_R-X_R^*X_R\label{11.6}
\end{eqnarray}
and such that the operator
$$
\begin{bmatrix}\K_{L} & \Lambda_1 & \ldots & \Lambda_d\\
\Lambda_1^* & \Phi_{11} & \ldots & \Phi_{1d}\\
\vdots & \vdots & & \vdots \\
\Lambda_d^* & \Phi_{d1} & \ldots & \Phi_{dd}\end{bmatrix}
$$
is positive semidefinite, where $\K_L$ is defined in \eqref{ex5}.
\label{T:ex8}
\end{Tm}
To get Theorem \ref{T:ex8} from Theorem \ref{T:ex7}, it suffices to
let $\K_{LR}=\begin{bmatrix}\Lambda_1 & \ldots & \Lambda_d\end{bmatrix}$
and $\K_R=[\Phi_{ij}]_{i,j=1}^d$ and to make use of block decompositions
\eqref{11.2} and \eqref{11.3}.

\subsection{The case of the noncommutative polydisk}

Here we consider the Givone-Roesser case (see Example \ref{E:NCGRimp}
above) where $S=R=E=\{1,\ldots,d\}$
and the tuples $T_L$ and $T_R$ are just $d$-tuples
$T_L=(T_{L,1},\ldots,T_{L,d})$ and $T_R=(T_{R,1},\ldots,T_{R,d})$
of contractive operators acting on $\cK_L$ and $\cK_R$, respectively.
\begin{Pb}
Given an admissible interpolation data set \eqref{1.26}, find necessary
and sufficient conditions for existence of a power series
$\T\in\mathcal{SA}_{\G^{\rm GR}}(\cU, \cY)$ such that
\begin{equation}
\left(X_L\T\right)^{\wedge L}(\del_L)=Y_L \quad\text{and}\quad
\left(\T Y_R\right)^{\wedge R}(\del_R)=X_R.
\label{10.1}
\end{equation}
\label{P:10.1}
\end{Pb}
The formulas \eqref{1.28}--\eqref{1.29a} read
$$
E_{L,j}=I_{\cK_L},\quad E_{R,j}=I_{\cK_R},\quad
\widetilde{N}_{j}(\del_L)=\del_{L,j}^*,\quad
\widetilde{M}_{j}(\del_R)=\del_{R,j}\quad (j=1,\ldots,d)
$$
and therefore, formulas \eqref{1.30} and \eqref{1.31} take the form
\begin{equation}
M_j=\begin{bmatrix} I_{\cK_L} & 0  \\ 0 & \del_{R,j}\end{bmatrix},
\quad N_j=\begin{bmatrix} \del_{L,j}^* & 0 \\ 0 &
I_{\cK_R}\end{bmatrix}\quad (j=1,\ldots,d).
\label{10.2}
\end{equation}
Theorem \ref{T:1.9} now reduces to
\begin{Tm}
There is a power series $\T\in{\mathcal {SA}}_\G^{\rm GR}(\cU, \, \cY)$
satisfying interpolation conditions \eqref{10.1} if and only if there
exist positive semidefinite operators
\begin{equation}
\K_j=\begin{bmatrix}\K_{j,L} & \K_{j,LR} \\ \K_{j,LR}^* &
\K_{j,R}\end{bmatrix}\in \cL(\cK_L\oplus\cK_R) \quad\mbox{for} \; \;
j=1,\ldots,d,
\label{10.3}
\end{equation}
that satisfy the Stein identity
\begin{equation}
\sum_{j=1}^d \left(M_j^* \K_j M_j-N_j^*\K_jN_j\right)=X^*X-Y^*Y,
\label{10.4}
\end{equation}
where $M_j$ and $N_j$ are the operators defined via formulas \eqref{10.2}
and $X$ and $Y$ are the same as in \eqref{1.33}.
\label{T:10.2}
\end{Tm}
Furthermore, it follows by Theorem \ref{T:1.11} that for every
choice of positive semidefinite operators $\K_1,\ldots,\K_d$ of the form
\eqref{10.3}, satisfying the Stein identity \eqref{10.4}, there exists
   a power series $\T\in{\mathcal {SA}}_\G^{\rm GR}(\cU, \, \cY)$ satisfying
(besides \eqref{10.1}) supplementary interpolation conditions
\begin{align}
(X_LH_j)^{\wedge L}(\del_L)\left[(X_LH_{j})^{\wedge L}(\del_L)\right]^*
& = \K_{j,L},\label{10.7}\\
      (X_LH_j)^{\wedge L}(\del_L)\left(G_j Y_R\right)^{\wedge R}(\del_R)
       & = \K_{j,LR},\nonumber\\
      \left[(G_jY_R)^{\wedge R}(\del_R)\right]^*(G_{j} Y_R)^{\wedge
R}(\del_R) & = \K_{j,R}\nonumber
\end{align}
for $j=1,\ldots,d$ and for some choice of associated
functions $H(z)$ and $G(z)$ in representations \eqref{1.20},
\eqref{1.22}, \eqref{1.23} of $\T$. Furthermore, all such $\T$ can be
parametrized by a linear fractional transformation.

\begin{Cy}
There is a power series $\T\in{\mathcal {SA}}_\G^{\rm GR}(\cU, \, \cY)$
satisfying the left interpolation condition
\begin{equation}
\left(X_L\T\right)^{\wedge L}(\del_L)=Y_L
\label{10.5}
\end{equation}
if and only if there exist positive semidefinite operators
$\K_{1,L},\ldots, \K_{d,L}\in\cL(\cK_L)$ that satisfy the Stein identity
\begin{equation}
\sum_{j=1}^d \left(\K_{j,L}-N_j^*\K_{j,L}N_j\right)=X_L^*X_L-Y_L^*Y_L.
\label{10.6}
\end{equation}
\label{C:10.3}
\end{Cy}
Again, for every choice of operators $\K_{1,L},\ldots, \K_{d,L}$
meeting conditions of Corollary \ref{C:10.3}, there exists $\T\in{\mathcal
{SA}}_\G^{\rm GR}(\cU, \, \cY)$ satisfying (besides the left condition
\eqref{10.5}) conditions \eqref{10.7} for $j=1,\ldots,d$ and for some
choice of associated function $H(z)$ in representation \eqref{1.20}  of
$\T$.

\subsection{The Schur interpolation problem}  \label{S:Schurint}

The classical Schur problem \cite{schur} is concerned with necessary and
sufficient conditions for existence of a (scalar valued) Schur function
$S$ with the preassigned first $n+1$ Taylor coefficients at the origin
(sometimes, especially if the Taylor coefficients at a point of $\D$
different from the origin, this problem is called the
Carath\'eodory-Fej\'er problem). The operator-valued analogue of this
problem is the following {\bf SP}: {\em given a collection of operators
$S_0,\ldots,S_n\in\cL(\cU,\cY)$, find necessary and sufficient conditions
for existence of a Schur function $S\in\cS(\cU,\cY)$ of the form
$$
S(z)=S_0+zS_1+\ldots+z^{n-1}S_{n-1}+\ldots.
$$}
The answer is given in terms of the Toeplitz matrix
$$
{\bf S}=\left[\begin{array}{cccc}S_0&0& \cdots &0\\ S_1
&S_0&\ddots&\vdots\\ \vdots&\ddots&\ddots &0\\ S_{n}&
\cdots&S_1 & S_0\end{array}\right]
$$
with operator entries: the {\bf SP} has a solution if and only if the
operator ${\bf S}: \; \cU^{n+1}\to \cY^{n+1}$ is contractive.
An interpolation problem with the data string $S_0,\ldots,S_N$ containing
gaps (that is, with unspecified $S_k$ for some $k<n$) also makes sense.
In fact, this is a completion question: is it possible to complete
a partially defined operator ${\bf S}$ as above to a contractive operator?
This question (even in the scalar valued case) is beyond our current
interests and will not be discussed here.

\medskip

Let $\G$ be an admissible graph and let $\cF_{E}$ be the free
semigroup generated by the edge set $E$ of $\G$. A subset
$\cF\subset\cF_E$ will be called {\em lower inclusive} if whenever
$v\in\cF$ and $v=uw$ for some $u,\, w\in\cF_E$, then it is the case that
also $u\in\cF$. A natural noncommutative analogue of the Schur problem is
the following:

\medskip

{\bf NSP:} Let $\G$ be an admissible graph, let $\cF_{E}$ be the free
semigroup generated by the edge set $E$ of $\G$ and let $\cF$ be a finite
lower inclusive subset of $\cF_{E}$. Given a collection of operators
$\{S_v\in\cL(\cU,\cY)\colon \; v\in\cF\}$, find necessary and sufficient
conditions for a noncommutative Schur-Agler function
$$
\T(z)=\sum_{v\in\cF_E}\T_v z^v\in{\mathcal S}{\mathcal A}_{\G}(\cU, \,
\cY)
$$
to exist such that
\begin{equation}
\T_v=S_v\quad\mbox{for every} \; \; v\in\cF.
\label{9.2}
\end{equation}

\medskip
We will show that conditions \eqref{9.2} can be written in the form
\begin{equation}
\left(X_L\T\right)^{\wedge L}(\del)=Y_L
\label{9.3}
\end{equation}
for an appropriate choice of $X_L$, $Y_L$ and $T=\{T_e: \; e\in E\}$;
in other words we will show that the {\bf NSP} is a particular
left-sided case of Problem \ref{P:1.4}.
The construction does not depend on the structure of the graph $\G$
and proceeds as follows.

We are given a lower inclusive subset $\cF$ of the free
semigroup $\cF_{E}$ together with and operator $F_{v} \in {\mathcal
L}(\cU, \cY)$ for each $v \in \cF$.  We let $\ell^{2}(\cF)$ be the
Hilbert space with orthonormal basis $\{ \delta_{v} \colon v \in
\cF\}$ indexed by $\cF$ and set $\cK_{L} = \ell^{2}(\cF) \otimes \cY$.
Note that elements of $\cK_{L}$ can also be viewed as functions $v
\mapsto f(v)$ on $\cF$ with values in $\cY$ subject to $\sum_{v \in
\cF} \|f(v)\|_{\cY}^{2} < \infty$. Note that the empty word
$\emptyset$ is in $\cF$ since $\cF$ is lower-inclusive.  Define an
operator $X_{L} \in {\mathcal L}(\cY, \cK_{L})$ by
$$
    X_{L} \colon y \mapsto \delta_{\emptyset} \otimes y.
$$
For each $e \in E$, we define an operator $\del_{L,e}$ on $\cK_{L}$ in
terms of matrix entries $\del_{L,e} = [\del_{L,e;v,w}]_{v,w \in \cF}$
(where each $\del_{L,e;v,w} \in {\mathcal L}(\cY)$)
by
$$
    \del_{L,e; v,w} = \begin{cases} I_{\cY}  &\text{if } v = we, \\
           0 &\text{otherwise,}
    \end{cases}
$$
or via the equivalent functional form
$$
     (T_{e}f)(v) = f(v \cdot e^{-1}) \text{ for } f \in \cK_{L}
$$
where we use the convention \eqref{1.10} and we declare
   $f(\text{undefined})=0$.
   Then it is easily checked that, given a formal power series $F(z) = \sum_{v
   \in \cF_{E}} F_{v} z^{v}$, the left left evaluation with operator
   argument $(X_{L}F)^{\wedge L}(\del_{L})$ works out to be given by
   $$
    \left( (X_{L}F)^{\wedge L}(\del_{L})u \right)(v) = F_{v} u
    \text{ for } v \in {\mathcal F} \text{ and all } u \in \cU.
   $$
   Hence, if we define $Y_{L} \colon \cU \to \cK_{L}$ by
   $$
    (Y_{L}u)(v) = S_{v}u,
   $$
   then the left tangential interpolation problem with operator argument
   associated with the data set ${\mathcal D} = (\del_{L}, X_{L},Y_{L})$
   is exactly equivalent to ${\bf NSP}$, and hence necessary and
   sufficient conditions for the ${\bf NSP}$ to have a solution
   can be derived from Theorem \ref{T:ex2}.

   \subsection{Interpolation with commutative data}

   For this example we consider the general Problems \ref{P:1.4} and
\ref{P:1.10} when the tuples $T_L$ and $T_R$ are commutative.
As explained in Section \ref{S:admissible}, the interpolation
conditions \eqref{1.3} imposed on a formal power series $S \in
\mathcal{S A}_{\G}(\cU, \cY)$ associated with Problem \ref{pb0} can
be expressed as interpolation conditions on the abelianized function
$S^{{\mathbf a}}$ of commuting variables $\lambda_{e_{1}}, \dots,
\lambda_{e_{d}}$:
\begin{equation}  \label{A1}
   (X_{L} S^{{\mathbf a}})^{\wedge L}(\del_{L}) = Y_{L}, \qquad
   (S^{{\mathbf a}} Y_{R})^{\wedge R}(\del_{R}) = X_{R}.
\end{equation}
Similarly, the additional interpolation  conditions
\eqref{1.37}--\eqref{1.39} imposed on $S \in \mathcal{SA}_{\cG}(\cU, \cY)$
by Problem \ref{P:1.10}  can be expressed as
interpolation conditions on the abelianized function $S^{{\mathbf
a}}$:
\begin{align}
(X_L H_{s}^{{\mathbf a}})^{\wedge L}(\del_L)
\left[(X_LH_{s'}^{{\mathbf a}})^{\wedge L}(\del_L)\right]^*
& = \Psi_{s,s'} \;  \text{ for } s, s'\in
S\colon \; [s]=[s'],  \notag \\
     (X_LH^{{\mathbf a}}_s)^{\wedge L}(\del_L)
     \left(G^{{\mathbf a}}_r Y_R\right)^{\wedge R}(\del_R)
        & = \Lambda_{s,r} \; \text{ for } s\in S; \; r\in
R\colon \; [s]=[r], \notag \\
       \left[(G^{{\mathbf a}}_rY_R)^{\wedge R}
       (\del_R)\right]^*(G^{{\mathbf a}}_{r'} Y_R)^{\wedge
R}(\del_R) & = \Phi_{r,r'} \; \text{ for } r,r'\in R\colon \;
[r]=[r'].
\label{A2}
\end{align}
  From the characterization of the class $\mathcal {S A}_{\cG}(\cU,
\cY)$ as transfer functions of conservative SNMLSs with structure
graph $\G$ and the counterpart of this result for the commutative
Schur-Agler class $\mathcal{S A}_{Z^{{\mathbf a}}_{\G}}(\cU, \cY)$
found in \cite{bb4}, it is clear that the abelianization $S^{{\mathbf
a}}$ of any element $S \in \mathcal{S A}_{\G}(\cU, \cY)$ is an
element of $\mathcal{S A}_{Z^{{\mathbf a}}}(\cU, \cY)$ as studied in
\cite{bb4, bb5}, and, conversely, any element $F$ of $\mathcal{S
A}_{Z^{{\mathbf a}}_{\G}}(\cU, \cY)$ lifts to an element $S \in
{\mathcal S A}_{\G}(\cU, \cY)$ (so $F = S^{{\mathbf a}}$).  The
results of \cite{bb5} can be applied to the abelianized problems
involving interpolation conditions \eqref{A1} (and possibly
also \eqref{A2}) for a function $S^{{\mathbf a}}$ in the
commutative Schur-Agler class $\mathcal{ S A}_{Z^{{\mathbf
a}}_{\cG}}(\cU, \cY)$.  When this is done, the Stein equation
\eqref{1.32} is the same as the Stein equation in \cite{bb5} where it
was shown to be the necessary and sufficient condition for the
abelianized interpolation problem to have a solution in the
commutative Schur-Agler class $\mathcal{ S A}_{Z^{{\mathbf
a}}_{\G}}(\cU, \cY)$.  In this way, we see that
interpolation problems for formal power series in noncommuting
indeterminants involving
commutative data reduces to the more standard interpolation
problems for analytic functions in commuting variables.

As an example, let us consider the case with commutative data for
the noncommutative-ball
setting discussed in Section \ref{S:noncomball}.
Let, in particular, $\cK_L=\C^n$, let $T_L$ is the
$d$-tuple of diagonal matrices constructed from $n$ points
$\lambda_i=\left(\lambda_i^{(1)},\ldots,\lambda_i^{(d)}\right)\in\B^d$
($i=1,\ldots,n$) by
$$
T_{L,j}={\rm diag} \, (\lambda_1^{(j)},\ldots,\lambda_n^{(j)})
\quad\mbox{for} \; \; j=1,\ldots,d,
$$
and let $X_L$ and $Y_L$ be conformally decomposed as
$$
X_L=\begin{bmatrix} b_1 \\ \vdots \\ b_n\end{bmatrix}\quad\mbox{and}\quad
Y_L=\begin{bmatrix} c_1 \\ \vdots \\ c_n\end{bmatrix}.
$$
Then the pair $(T_L, \, X_L)$ is left admissible and it is easily seen
that
$$
(X_L\T)^{\wedge L}(T_L)=\operatorname{Col}_{1\le i\le n}b_i\T(\lambda_i)
$$
(where $\T(\lambda_i)$ is defined via \eqref{scalarvariable},
so that condition
\eqref{ex1} collapses to $n$ left sided conditions
\begin{equation}
b_i\T(\lambda_i)=c_i\quad (i=1,\ldots,n).
\label{ex50}
\end{equation}
Furthermore, the matrix $\K_{L}$ in \eqref{ex5} admits a more explicit
representation
\begin{equation}
\K_{L}=\left[\frac{b_ib_j^*-c_ic_j^*}{1-\langle
\lambda_i, \, \lambda_j\rangle}\right]_{i,j=1}^n
\label{ex51}
\end{equation}
where $\langle \lambda_i, \, \lambda_j\rangle$ stands for the standard
inner product in $\C^d$. Thus, Theorem \ref{T:ex5} gives \cite[Theorem
4.1]{popescu1}: {\em there exists
a power series $\T\in{\mathcal {SA}}_{\G^{\rm FM}}(\cU, \, \cY)$
satisfying interpolation conditions \eqref{ex50} if and only if
the matrix $\K_{L}$ defined in \eqref{ex51} is positive semidefinite}.

We note that the commutative (several-variable) analogue of the
Schur problem discussed above in Section \ref{S:Schurint} is one of
the examples for the commutative theory discussed in \cite{bb5}.

\bibliographystyle{amsplain}

\begin{thebibliography}{10}


\bibitem{agler2}
J.~Agler, {\it Some interpolation theorems of Nevanlinna--Pick type},
Preprint, 1988.

\bibitem{agler-hellinger}
J.~Agler, {\it On the representation of certain holomorphic functions
defined on a polydisk}, in {\it Topics in Operator Theory: Ernst D.
Hellinger memorial Volume} (L.~de~Branges, I.~Gohberg and J.~Rovnyak,
eds.), pp. 47--66,  \textbf{OT 48},
Birkh\"auser Verlag, Basel, 1990.

\bibitem{aglmccar-poly} J.~Agler and J.E.~McCarthy,
{\it Nevanlinna-Pick interpolation on the bidisk}, J.~Reine
Angew.~Math. {\bf 506} (1999), 191--204.

\bibitem{aglmccar} J.~Agler and J.~E.~McCarthy,  \emph{Complete
Nevanlinna-Pick kernels}, J. Funct. Anal., {\bf  175} (2000), 111--124.

\bibitem{AK1} D.~Alpay and D.S.~Kalyuzhny\u{\i}-Verbovetzki\u{\i},
\emph{On the intersection of null spaces for matrix substitutions in a
non-commutative rational formal power series}, Comptes rendus
Mathematiques Acad. Sci. Paris I {\bf 339} (2004) 533--538.

\bibitem{AK2} D.~Alpay and D.S.~Kalyuzhny\u{\i}-Verbovetzki\u{\i},
\emph{Matrix-$J$-unitary non-commutative rational formal power
series}, in {\it Linear Operators and Systems} (Ed. D.~Alpay and
I.~Gohberg), OT volume, Birkh\"auser-Verlag, Basel-Boston-Berlin, to
appear.


\bibitem{AE}  C.-G.~Ambrozie and J.~Eschmeier, A commutant
lifting theorem on analytic polyhedra, {\em Proceedings of Operator
Theory Conference Dedicated to Prof.  Wieslaw Zelazko}, Banach Center
publ., Warszawa, to appear.

\bibitem{at}
C.-G.~Ambrozie and D.~Timotin.
\newblock {\em A von Neumann type inequality for certain domains in
$\C^n$}, Proc.~Amer.~Math.~Soc., \textbf{131} (2003), 859--869.


\bibitem{ariaspopescu}
A.~Arias and G.~Popescu, {\em Noncommutative interpolation and Poisson
transforms}, Israel J. Math. {\bf 115} (2000), 205--234.

\bibitem{arovgros1}
D.Z.~Arov and L.~Z.~Grossman,
{\it Scattering matrices in the theory of unitary
extensions of isometric operators},
Soviet Math. Dokl. {\bf 270} (1983), 17--20, MR0705184 (85c:47008),
Zbl 0543.47010.

\bibitem{arovgros2}
D.Z.~Arov and L.~Z.~Grossman,
{\it Scattering matrices in the theory of unitary
extensions of isometric operators},
Math. Nachr. {\bf 157} (1992), 105--123.


\bibitem{arv}
W.~Arveson, \emph{ Subalgebras of $C^*$-algebras. III. Multivariable
operator theory}, Acta Math. {\textbf 181} (1998), no. 2,
159--228.

\bibitem{bb3}
J.A.~Ball and V.~Bolotnikov, {\em A tangential interpolation problem on
the distinguished boundary of the polydisk for the Schur-Agler class},
J. Math. Anal. Appl. {\bf 273} (2002), no. 2, 328--348.

\bibitem{bb4}
J.A.~Ball and V.~Bolotnikov, \emph{Realization and interpolation
for Schur-Agler-class functions on domains with matrix polynomial
defining function in ${\bf C}^n$}, J. Funct. Anal. {\bf 213} (2004),
45--87.


\bibitem{bb5}
J.A.~Ball and V.~Bolotnikov, \emph{Interpolation problems with
operator argument for contractive-valued functions on general domains in
${\bf C}^n$}, Preprint.

\bibitem{BGM1} J.A.~Ball, G.~Groenewald and T.~Malakorn,
\emph{Structured noncommutative multidimensional linear systems},
SIAM J.~Control and Optimization, to appear.

\bibitem{BGM2} J.A.~Ball, G.~Groenewald and T.~Malakorn, Conservative
structured noncommutative multidimensional linear systems, in
{\it Linear Operators and Systems} (Ed. D.~Alpay and
I.~Gohberg), OT volume, Birkh\"auser-Verlag, Basel-Boston-Berlin, to
appear.

\bibitem{BGR} J.A.~Ball, I.~Gohberg and L.~Rodman, \emph{Interpolation
of Rational Matrix Functions}, OT45 Birkh\"auser-Verlag, 1990.

\bibitem{BLTT}
J.A.~Ball, W.~S.~Li, D.~Timotin and T.~T.~Trent, {\it A
commutant lifting theorem on the polydisc},
Indiana University Math. J. {\bf 48} (1999), 653--675.

\bibitem{BallMal}  J.A.~Ball and T.~Malakorn, {\it Multidimensional
linear feedback control systems and interpolation problems for
multivariable holomorphic functions},
Multidimens. Systems and Signal Process. {\bf 15} (2004), 7--36.

\bibitem{BSV} J.A.~Ball, C.~Sadosky and V.~Vinnikov, {\it Scattering
systems with several evolutions and multidimensional
input/state/output systems},  Integral Equations and Operator
Theory, to appear.

\bibitem{bt}
J.A.~Ball and T.~Trent,
{\it Unitary colligations, reproducing kernel Hilbert
spaces and Nevanlinna--Pick interpolation in several variables},
J. Funct. Anal. {\bf 157} (1998), no.1, 1--61.


\bibitem{BTV}
J.A.~Ball, T.~T.~Trent and V.~Vinnikov, \emph{Interpolation and
commutant lifting for multipliers on reproducing kernels Hilbert
spaces}, in  {\emph Operator Theory and Analysis: The M.A. Kaashoek
Anniversary Volume (Workshop in Amsterdam, Nov. 1997)}, OT 122,
Birkh\"auser-Verlag, Basel-Boston, 2001, pp. 89--138.

\bibitem{BV} J.A.~Ball and V.~Vinnikov, \emph{Lax-Phillips scattering
and conservative linear systems:  A Cuntz-algebra multidimensional
setting}, AMS Memoir, to appear.

\bibitem{CH} R.E.~Curto and D.A.~Herrero, {\em On closures of joint
similarity orbits},  Integral Equations and Operator Theory
\textbf{8} (1985), 489--556.

\bibitem{CJ}  T.~Constantinescu and J.L.~Johnson, \emph{A note on
noncommutative interpolation}, Canadian Math.~Bull.~{\bf 46} (2003) no. 1
59--70.

\bibitem{davidsonpitts}
K.R.~Davidson and D.R.~Pitts, \emph{Nevanlinna--Pick interpolation for
non-commutative analytic Toeplitz algebras}, Integral
Equations Operator Theory \textbf{31} (1998), no. 3, 321--337.


\bibitem{Drury} S.W.~Drury, {\em A generalization of von Neumann's
inequality to the complex ball}. Proc. Amer. Math. Soc. {\bf 68} (1978),
300--304.

\bibitem{FFGK} C.~Foia\c{s}, A.~Frazho, I.~Gohberg and M.A.~Kaashoek,
{\em Metric Constrained Interpolation, Commutant Lifting and
Systems}, \textbf{OT100}, Birkh\"auser-Verlag, Boston-Basel, 1998.

\bibitem{FM}
E.~Fornasini and G.~Marchesini, {\it State-space
realization theory of two-dimensional filters}, IEEE Trans. Automat.
Contr. {\bf AC-21}, No. 4, 1976, 484--492.


\bibitem{HMV} J.W.~Helton, S.~McCullough and V.~Vinnikov,
\emph{Noncommutative convexity arises from linear matrix
inequalities}, preprint, 2005.

\bibitem{Kac}
T.~Kaczorek, {\it Two-Dimensional Linear Systems},
Lecture Notes in Control and Information Sciences {\bf 68}
Springer-Verlag, Berlin, 1985.

\bibitem{KV} D.S.~Kalyuzhny\u{\i}-Verbovetzki\u{\i} and V.~Vinnikov,
\emph{Non-commutative positive kernels and their matrix functions},
Proceedings of the American Mathematical Society, to appear.

\bibitem{mccull2}
S.~McCullough, \emph{The local de {B}ranges-{R}ovnyak construction and
complete {N}evanlinna-{P}ick kernels}, in
\emph{Algebraic methods in operator theory} (Ed. R.~Curto and
P.E.T.~Jorgensen), Birkh\"auser--Verlag, Boston, 1994, pp. 15--24.

\bibitem{Mal}  T.~Malakorn, {\em Multidimensional Linear Systems and
Robust Control}, Dissertation, Department of Electrical and Computer
Engineering, Virginia Tech (April, 2003).

\bibitem{MS} P.S.~Muhly and B.~Solel, \emph{Hardy algebras, $W^{*}$
correspondences and interpolation theory}, Math.~Ann.~ \textbf{330}
(2004), no. 2, 353--415.

\bibitem{Popescu-CLT1} G.~Popescu, Isometric dilations for infinite
sequences of noncommuting operators, {\it Trans. Amer. Math. Soc.}
316 (1989), 523--536.


\bibitem{Popescu-multi} G.~Popescu, \emph{Multi-analytic operators on
Fock spaces}, Math.~Ann.~\textbf{303} (1995), 31--46.

\bibitem{popescu1}
G.~Popescu, \emph{Interpolation problems in several variables},
J. Math. Anal. Appl., {\textbf 227} (1998), 227--250.

\bibitem{Popescu-Poisson}  G.~Popescu, {\it Poisson transforms on some
$C^{*}$-algebras generated by isometries},  J. Funct. Anal. {\bf 161}
(1999), 27--61.


\bibitem{quig}
P.~Quiggin, \emph{For which reproducing kernel {H}ilbert
spaces is {P}ick's theorem true?} Integral Equations Operator Theory
{\textbf 16} (1993), no. 2, 244--266.

\bibitem{Roesser}
R.P.~Roesser, {\it A discrete state-space model for
linear image processing}, IEEE Trans. Automat. Control {\bf AC-20} (1975),
no. 1, 1--10.

\bibitem{RR} M.~Rosenblum and J.~Rovnyak, \emph{Hardy classes and
operator theory}, Oxford Mathematical Monographs, Oxford University
Press, New York, 1985; Dover republication, New York, 1997. MR
\textbf{97j:47002}

\bibitem{Rudin-ball}
W.~Rudin, \emph{Function theory in the unit ball of ${\bf {C}}\sp{n}$},
Springer-Verlag, New York, 1980.

\bibitem{sarasonbook}
D.~Sarason, \emph{ Sub-{H}ardy {H}ilbert spaces in the unit disk},
John Wiley and Sons Inc., New York, 1994.

\bibitem{schur} I.~Schur, \emph{ \"Uber Potenzreihen die im Innern
des Einheitskreises  Beschr\"{a}nkt Sind},. J. Reine Angew. Math.{\bf
14}, 7 (1917),  205--232.


\bibitem{taylor1}
J.L.~Taylor, {\it A joint spectrum for several commuting operators}, J.
Funct. Anal. {\bf 6} (1970), 172--191.

\bibitem{taylor2}
J.L.~Taylor, {\it The analytic-functional calculus for several commuting
operators},  Acta Math.  {\bf 125} (1970) 1--38.

\bibitem{tomerlin}
A.T.~Tomerlin, \emph{Products of {N}evanlinna-{P}ick
kernels and operator
colligations}, Integral Equations Operator Theory {\bf
38} (2000), no. 3, 350--356.

\bibitem{vasilescu}
F.-H.~Vasilescu, {\it A Martinelli type formula for the
analytic functional calculus}, Rev. Roum. Math. Pures Appl.
{\bf 23} (1978), no. 10, 1587--1605.

\end{thebibliography}
\providecommand{\bysame}{\leavevmode\hbox to3em{\hrulefill}\thinspace}

\end{document}